\documentclass{article}
\usepackage[utf8]{inputenc}

\usepackage[margin=1in]{geometry}
\usepackage{amsfonts}
\usepackage{amsthm}
\usepackage{amsmath}
\usepackage{mathtools}
\usepackage{amssymb}
\usepackage{siunitx}
\usepackage{verbatim}
\usepackage{physics}
\usepackage{enumerate}
\usepackage{multicol}
\usepackage[noadjust]{cite}
\allowdisplaybreaks
\usepackage{caption}
\usepackage{subcaption}
\usepackage{graphicx}
\usepackage{tcolorbox}
\usepackage{dsfont}
\usepackage{xcolor}
\usepackage{tikz-cd}
\usetikzlibrary{backgrounds}
\usetikzlibrary{babel}
\usepackage{subfiles}
\usepackage{longtable}
\usepackage{hyperref}
\usepackage{bbm}
\makeatletter
\def\blfootnote{\gdef\@thefnmark{}\@footnotetext}
\makeatother
\newtheorem{theorem}{Theorem}[section]
\newtheorem*{theorem*}{Theorem}
\newtheorem{corollary}{Corollary}[theorem]

\newtheorem{lemma}[theorem]{Lemma}
\newtheorem{defn}[theorem]{Definition}
\newtheorem{propn}[theorem]{Proposition}

\newcommand{\B}[1]{\mathbb{#1}}
\newcommand{\C}[1]{\mathcal{#1}}
\newcommand{\oop}{\preccurlyeq}
\newcommand{\is}{\mathbbm{1}}
\DeclareMathOperator{\erfc}{erfc}

\usepackage[skip=7.5pt plus1pt, indent=20pt]{parskip}

\title{Two-Sided Free Boundary Problems Arising From Branching-Selection Particle Systems}
\author{Jacob Mercer\footnote{\texttt{jacob.mercer@maths.ox.ac.uk}, Department of Statistics, University of Oxford}}
\date{}  

\begin{document}
\maketitle

\begin{abstract}
    We introduce and analyse a two-sided branching-selection particle system which generalises the well-known $N$-particle branching Brownian motion ($N$-BBM) model, which we call the $(N,p)$-BBM, where either the leftmost or rightmost particle is deleted at each branching event according to a parameter $p\in(0,1)$. We establish that, as $N\to\infty$, the empirical distribution of the $(N,p)$-BBM converges to a deterministic hydrodynamic limit described by a free boundary problem on a finite interval with two moving boundaries, and Neumann and Dirichlet boundary conditions parametrized by $p$. Again, this generalises the one-sided free boundary problem which characterises the hydrodynamic limit of the $N$-BBM. Existence and regularity of the free boundary problem is also proved, by appealing to a connection with inverse first passage problems. We further prove that the asymptotic velocity $v_{N,p}$ of the $(N,p)$-BBM converges, as $N\to\infty$, to $v_p$, the unique travelling wave speed of the limiting free boundary problem. These results generalize previous one-sided models and connect to broader classes of free boundary problems found in evolutionary dynamics and flame propagation.
\end{abstract}

\section{Introduction} \blfootnote{This publication is based on work supported by the EPSRC Centre for Doctoral Training in Mathematics of Random Systems: Analysis, Modelling, and Simulation (EP/S023925/1)}
Survival-of-the-fittest is the mechanism by which, in an environment with limited resources, fitter members of a population out-compete less fit individuals. Over time, this `selective-pressure' drives species to evolve towards higher and higher fitness. A number of questions naturally arise about this selective mechanism. At what rate does fitness increase in a population due to survival-of-the-fittest? How does this increase in fitness depend on the population size? What is the variance or distribution of fitness levels in a population? 

Physicists Brunet and Derrida pioneered the study of so called `branching-selection particle systems' or `Brunet-Derrida systems' to understand these questions about survival-of-the-fittest in a mathematical framework (\cite{shiftCutoff},\cite{microModelsCutoff},\cite{effectOfMicroNoise}). These branching-selection systems are characterised by having a fixed population size, $N$, and a selection mechanism in order to enforce this fixed population size when a birth occurs. 

A well-studied branching-selection system is called the $N$-BBM. In the $N$-BBM, we have $N$ particles on the real line which move about as independent standard Brownian motions. Each particle, at rate $1$, splits into $2$ particles at the same location (the binary branching step) and simultaneously, we delete the particle with the leftmost position (the selection step). 

Much is known about the $N$-BBM system. Firstly, we have a `hydrodynamic limit' result (\cite{hydroNBBM}); as $N\to\infty$, the empirical density of particles in the $N$-BBM system with initial density $\rho$ (in the sense that each particle independently has initial density $\rho$) converges to a deterministic limit. Moreover, this limit is the solution of the free boundary problem:
\begin{align} \label{nbbmFBP}
\begin{cases}
    u_t=\frac{1}{2}u_{xx} + u  & x\in(L_t,\infty), t>0, \\
    u(x,t)=0 & x\notin (L_t,\infty), t>0, \\
    \int_{L_t}^\infty u(y,t)dy=1 & t>0, \\
    u \text{ continuous},\\
    u(x,0)=\rho(x).\\
\end{cases}
\end{align}

So the $\frac12u_{xx}$ term describes the diffusion of the Brownian motions, the $+u$ term describes the fact that each particle is branching at unit rate, and the $\int_{L_t}^\infty u(y,t)dy=1$ occurs due to the fact that the population size is constant (thus the empirical density truly remains a probability density for all time). The free boundary $L_t$ represents the location of the leftmost particle, left of which we trivially have no mass. Note that because this is a free boundary problem, a solution is not just a function $u$, but a pair $(u,L)$ satisfying \eqref{nbbmFBP}. An important question to clarify is whether, or for which $\rho$, a solution exists. For the free boundary problem \eqref{nbbmFBP}, this is shown in \cite{globalSolsNBBM}.

We also know what the asymptotic speed of the particle system is; that is, the asymptotic speed at which fitness increases in the population due to the selection mechanism. If $X_i^N(t)$ denotes the position of the $i$\textsuperscript{th} leftmost particle of the $N$-BBM at time $t$, then it is known (see Appendix to \cite{me}) that
$$\lim_{t\to\infty}X_i^N(t)/t = v_N = \sqrt{2}-\frac{\pi^2}{\sqrt{2}\log^2 N}+ o((\log N)^{-2})$$
for any $i=1,2,\ldots,N$. Now the FBP (\ref{nbbmFBP}) has infinitely many travelling wave solutions, one for each speed $c\in [\sqrt{2},\infty)$. Therefore we can observe that as $N\to\infty$, the $N$-BBM `selects' the minimal travelling wave speed, since $v_N \to \sqrt{2}$. This is called the \textit{weak selection principle}. A stronger version of this result has also been attained by Berestycki and Tough \cite{beresTough}, in which they show that not only does the asymptotic velocity select the minimal travelling wave speed, but that in fact the stationary empirical measure of the $N$-particle system converges to the minimal travelling wave; this is called the \textit{selection principle}.

The $N$-BBM has been generalised in a number of ways. For example, one can consider the $N$-BBM in which when a particle branches, its offspring is a random distance from itself (see \cite{durrettRemenik} or \cite{nonLocalNBBM}). Accordingly, the hydrodynamic limit of these systems is described by a non-local PDE with a single free boundary. 

In this paper we will study a generalisation of the $N$-BBM which we call the $(N,p)$-BBM. In the $(N,p)$-BBM, we again have $N$ particles moving about as independent standard Brownian motions on the real line. However in the $(N,p)$-BBM, when a branching event occurs, we may delete either the leftmost or rightmost particle. We kill the leftmost with probability $p$ and otherwise we kill the rightmost. So the $(N,p)$-BBM is a generalisation, with $p=1$ corresponding to $N$-BBM.

The motivation for studying this system is that, where the hydrodynamic limit of the $N$-BBM is a one-sided free boundary problem, the hydrodynamic limit of the $(N,p)$-BBM (for $p\in (0,1)$) is a two-sided free boundary problem. Two-sided free boundary problems appear in a number of different contexts. Most notably, two-sided free boundary problems have been considered to study predator-prey evolution (\cite{predPreyFBP}, Section 5 of \cite{spreadVanDichot}) or to study flame propagation (\cite{guoHuKoh03},\cite{guoHu06}). 

Recall that in the one-sided case, the $N$-BBM is a toy model in which $\B{R}$ is a fitness landscape; location on $\B{R}$ denoting relative fitness. In the two-sided case of the $(N,p)$-BBM, we could consider location on $\B{R}$ being a 1-dimensional summary of some trait which does not correspond to fitness. Then killing on the right and left may correspond to different causes of death which affect individuals at either extreme of the trait with different strengths. 

The main result of this paper is to prove that the hydrodynamic limit of the $(N,p)$-BBM particle system with initial density $\rho$, under suitable conditions, is the solution to the free boundary problem
\begin{align} \label{fbpMain}
\begin{cases}
    u_t=\frac{1}{2}u_{xx} + u  & x\in(L_t,R_t), t>0, \\
    u(x,t)=0 &  x\notin (L_t,R_t), t>0, \\
    u_x(L_t,t)=2p, u_x(R_t,t)=2(p-1) & t>0, \\
    u \text{  continuous}, \\
    u(x,0)=\rho(x).
\end{cases}
\end{align}

Of course, it is natural to ask whether the free boundary problem has a solution and whether it is unique. We will prove the affirmative; that the equation does in fact have a unique solution and further that any solution has highly regular boundaries. 

The proof connects this work into another area of problems called \textit{inverse first passage problems}. In general, inverse first passage problems ask whether, given a probability distribution on $[0,\infty)$, it is possible to construct a boundary such that the first passage time of Brownian motion to this boundary is distributed according to the given distribution. Just as with branching-selection particle systems and free boundary problems, these have been studied in the one-sided case \cite{ccs22}\cite{cccs11}\cite{cccs06}\cite{ekstromJanson} and two-sided symmetric case \cite{anulova} on $\B{R}$ (corresponding to the $N$-BBM \cite{hydroNBBM} and Brownian bees \cite{afree}\cite{beesBBNP} respectively). In our case, we will see that if $(u,L,R)$ is a solution to \eqref{fbpMain}, then for Brownian motion $W$ with initial density $\rho$, the boundaries $L$ and $R$ are such that the first exit time of $(L_s,R_s)$ by $W_s$ is exponentially distributed, and the probability that $W$ exits on the left (resp. right) is $p$ (resp. $1-p$). We draw on the extensive literature on the inverse first passage problem to help us show existence and regularity of the solution of \eqref{fbpMain}. 

\section{Statement of results} \label{results}
Consider the free boundary problem (\ref{fbpMain}). We prove in Section \ref{probRepOfSolutions} that there exists a solution of the FBP \eqref{fbpMain} with highly regular free boundaries $L,R\in C^\infty((0,\infty))$, and, in Section \ref{comparison}, by sandwiching any solution of the FBP, that the solution is unique. 

\vspace{1em}

\begin{theorem} \label{fbpExistence}
    For any probability density function $\rho$ on $\B{R}$, there exists a unique triplet $(u,L,R)$ solving the FBP \eqref{fbpMain}, and moreover $L,R\in C^\infty((0,\infty))$. Moreover, the boundaries $(L,R)$ are such that if $W$ has initial density and $\tau:=\inf\{t>0:W_t\notin (L_t,R_t)\}$, then $\B{P}(\tau>t)=e^{-t}$ and $\B{P}(W_\tau=L_\tau)=p=1-\B{P}(W_\tau=R_\tau)$; that is $(L,R)$ solves the asymmetric inverse first passage problem. 
\end{theorem}

Given, then, that a solution to the desired FBP exists, we can state the hydrodynamic limit theorem for the $(N,p)$-BBM. We prove that the empirical distribution of the $(N,p)$-BBM converges, and describe the rate of convergence.

\vspace{1em}

\begin{theorem} \label{hydroLimitThm}
    Consider the $(N,p)$-BBM process $(X^{N,p}(t))_{t\geq 0}$ such that, at time $t=0$, each particle independently is distributed with density $\rho$. Let $\pi^{N,p}_{\rho,t}$ be the empirical distribution of the $N$ particles of $(X^{N,p}(t))_{t\geq 0}$ at time $t$. Fix $t>0$, and let $u(x,t)$ be the solution to the FBP \eqref{fbpMain} with initial condition $\rho$. Then there exists a constant $C=C(t)$ such that for any $a\in \B{R}$ and $B,\beta>0$, for $N$ sufficiently large
    $$\B{P}\Bigg( \Bigg|\int_a^\infty \pi^N_{\rho,t}(dr) - \int_a^\infty u(r,t)dr\Bigg|>BN^{-\beta}\Bigg)<B^{-9}CN^{9\beta-2}.$$
\end{theorem}

Not only does the empirical distribution converge to $u$, be we also prove (in Section \ref{leftmostrightmost}) the convergence of the left and rightmost particles of $X^{N,p}$ to the left and right barriers, $L$ and $R$, of \eqref{fbpMain}. 

\vspace{1em}

\begin{theorem} \label{leftmostRightmostTheorem}
    Let $(u(x,t),L_t,R_t)_{t\geq 0,x\in \B{R}}$ be the solution to the free boundary problem (\ref{fbpMain}). For fixed $t$, $X_1^N(t)$ and $X_N^N(t)$ converge as $N\to\infty$ to $L_t$ and $R_t$ respectively, almost surely and in expectation.,
\end{theorem}

The strategy for proving Theorem \ref{hydroLimitThm} follows the ideas that De Masi et al. used for proving the hydrodynamic limit of the $N$-BBM, and can be summarised by the following diagram. 

\begin{center}
       \begin{tikzcd}[row sep=huge, column sep=tiny, inner frame sep=5mm]
           X^{N,p,\delta,-}\arrow[d]&&\preccurlyeq && X^{N,p}\arrow[d,dashed,color=red] &&\preccurlyeq && X^{N,p,\delta,+} \arrow[d] \\
            u^{\delta,-}(x,t) && \xrightarrow[\leq]{\delta\to 0}
          && u(x,t)&& \xleftarrow[\leq]{\delta \to 0} && u^{\delta,+}(x,t)
     \end{tikzcd}
\end{center}

The aim is to prove the convergence described by the red arrow; the convergence of the empirical density of $X^{N,p}$ to the solution $u$ of the free boundary problem as $N\to\infty$. To avoid dealing with deletion of particles at random times, we consider instead upper and lower boundaries with deletion only at discrete times $k\delta$ for $k\in \B{N}$, $\delta>0$. That is, we couple $X^{N,p}$ to particle systems $X^{N,p,\delta,\pm}$ so that
$$X^{N,p,\delta,-}(0)\oop X^{N,p}(0) \oop X^{N,p,\delta,+}(0) \implies X^{N,p,\delta,-}(k\delta)\oop X^{N,p}(k\delta) \oop X^{N,p,\delta,+}(k\delta),$$
for all $k\in \B{N}$ and such that $X^{N,p,\delta,\pm}$ are branching Brownian motions with selection only at discrete times $k\delta$, $k\in \B{N}$. For small $\delta>0$, these closely approximate $X^{N,p}$, and their hydrodynamic limits are easier to determine. This coupling is constructed and the ordering proven in Section \ref{twoSidedCoupling}. 

We then show in Section \ref{convOfCoupling} that the empirical distributions of the upper and lower particle systems $X^{N,p,\delta,+}$ and $X^{N,p,\delta,-}$ converge respectively to deterministic limits $u^{\delta,+}$ and $u^{\delta,-}$ as $N\to\infty$. Moreover, we prove in Section \ref{detBarrierSameLimit} that as $\delta \downarrow 0$, the functions $u^{\delta,+}$ and $u^{\delta,-}$ become arbitrarily close and converge to the same limit. Finally, in order to prove Theorem \ref{hydroLimitThm}, we prove that there is an ordering of the deterministic functions $u^{\delta,-}$ and $u^{\delta,+}$ relative to the solution of the free boundary problem (\ref{fbpMain}); that is, given the same initial condition, then
$$u^{\delta,-}(x,t) \leq u(x,t) \leq u^{\delta,+}(x,t) \quad \forall t\geq 0.$$
This is the done in Section \ref{comparison}, using the probabilistic representation of solutions $u$ which we show in Section \ref{probRepOfSolutions}. Section \ref{comparison} concludes with a proof of Theorems \ref{fbpExistence} and \ref{hydroLimitThm}.

Next, we show that this free boundary problem (\ref{fbpMain}) has a unique travelling wave solution with a speed which we can explicitly calculate.

\vspace{1em}

\begin{propn} \label{travellingWave}
    For $p\in (0,1)$, the free boundary problem (\ref{fbpMain}) has a unique (up to translation) travelling wave solution
    $$(u(x,t),L_t,R_t)=(w(x-ct),ct,R_0+ct)=\left(\frac{2p}{\sqrt{2-c^2}}e^{-c(x-ct)}\sin((x-ct)\sqrt{2-c^2}),ct,\frac{\pi}{\sqrt{2-c^2}}+ct\right),$$
    with unique wave speed $c=\sqrt{\frac{2\log^2\left(p/(1-p)\right)}{\log^2\left(p/(1-p)\right)+\pi^2}}$ when $p\geq 1/2$ and $c=-\sqrt{\frac{2\log^2\left(p/(1-p)\right)}{\log^2\left(p/(1-p)\right)+\pi^2}}$ when $p\leq 1/2$
\end{propn}

\begin{proof}
    Fix $p$ and suppose that $(w(x-ct),ct,R_0+ct)$ is a travelling wave solution with speed $c$, with $w(x)\geq 0\iff x\in (0,R_0)$. Plugging $u(x,t)=w(x-ct)$ into $u_t = \frac12 u_{xx} + u$, we get $\frac12 w'' + cw' + w = 0$, a second order linear ODE with solution $w(x)=Ae^{-(c-\sqrt{c^2-2})x}+Be^{-(c+\sqrt{c^2-2})x}$ for some $A,B\in \B{C}$. 

    Now let us apply the boundary conditions to $w$. The fact that $u$ is continuous gives that $w(0)=w(R_0)=0$, and $w(0)=0$ implies that $A=-B$. Suppose for contradiction that there is a solution with wave speed $c > \sqrt{2}$.  Then $w(R_0)=0$ implies that $Ae^{-(c-\sqrt{c^2-2})R_0}-Ae^{-(c+\sqrt{c^2-2})R_0}=0$ which means that we must have $A=0$ or $R_0\in \{0,\infty\}$. Trivially, we cannot have $A=0$ or $R_0=0$, and $R_0=\infty$ leads to the contradiction $2(p-1)=w'(R_0)=w'(\infty)=0$, therefore $c\leq \sqrt{2}$, and the exponents become complex. Using $e^{ix}=\cos(x)+i\sin(x)$ gives that the solution must be of the form $w(x)=\tilde{A}e^{-cx}\sin(x\sqrt{2-c^2})$. 

    Then $w(R_0)=0$ and $w\geq 0$ gives that $R_0=\pi/\sqrt{2-c^2}$, and $w'(0)=2p$ gives that $\tilde{A}=2p/\sqrt{2-c^2}$. Finally $w'(R_0)=2(p-1)$ gives that $2(p-1)=-2p\exp(-c\pi/\sqrt{2-c^2})$, which we solve to get that $c=\pm\sqrt{\frac{2\log^2\left(p/(1-p)\right)}{\log^2\left(p/(1-p)\right)+\pi^2}}$, as required.
\end{proof}

Note that in the special cases $p=0$ or $p=1$, the $(N,p)$-BBM is either an $N$-BBM system (which has asymptotic velocity converging to $\sqrt{2}$) or a reflection of $N$-BBM (which has asymptotic velocity $-\sqrt{2}$), and this agrees with taking limits $p\uparrow 1$ and $p\downarrow 0$ respectively. However the $(N,p)$-BBM differs from the $N$-BBM in that there is a \textit{unique} travelling wave solution for the hydrodynamic limit of the $(N,p)$-BBM, since the solution has compact support, whereas the hydrodynamic limit of the $N$-BBM has infinitely many travelling wave solutions. The support of the solution of (\ref{nbbmFBP}) is semi-infinite, and different travelling wave solutions correspond to different rates of decay at $+\infty$. In Section \ref{speed}, we prove that the limiting asymptotic velocity of the $(N,p)$-BBM coincides with the unique travelling wave speed $c$ of \eqref{fbpMain} which we found in Proposition \ref{travellingWave}.

\vspace{1em}

\begin{theorem} \label{speedThm}
    Let $X^{N,p}$ be an $(N,p)$-BBM process with initial condition $\rho$ which has bounded support. Then for every fixed $N$, there exists a constant $v_N$ such that
    $$\lim_{t\to\infty}X_1^{N,p}(t)/t = \lim_{t\to\infty}X_N^{N,p}(t)/t=v_N,$$
    where $X_1^{N,p}$ and $X_N^{N,p}$ are the leftmost and rightmost particles of the system, respectively, and moreover,
    $\lim_{N\to\infty}v_N = c$, where $c$ is the unique travelling wave speed of the free boundary problem (\ref{fbpMain}).
\end{theorem}

Although it is not proven here, we expect that $v_N - \lim_{N\to\infty}v_N$ is on the order of $(\log N)^{-2}$, which is what we see for the $N$-BBM and $N$-BRW. 

This limiting asymptotic speed also agrees with another heuristic argument using results about branching Brownian motion in a compact interval with absorption at the boundaries (see \cite{strip} for a study of this system in criticality). Consider branching Brownian motion with drift $-\mu$ in a compact interval $(0,K)$ in which a particle is killed when it exits $(0,K)$. It is known that the critical width of the interval for survival is $K=\pi /\sqrt{2-\mu^2}$ (see for example Proposition 1 of \cite{strip}). At this critical width, the expected population size is roughly constant, and the density of particles at large times is approximately (Lemma 5 of \cite{bbs1}):
$$ f(x):=Ze^{-\mu x}\sin(\frac{\pi x}{K})=Ze^{-\mu x}\sin(\sqrt{2-\mu^2}x),$$
where $Z$ is a normalising constant independent of $x$. Suppose we were to choose $\mu$ so that roughly a proportion $p$ of particles exit $(0,K)$ at $0$ and a proportion $1-p$ at $K$, in order that the system behaves roughly like an $(N,p)$-BBM. In order to see $f'(0)=2p$ and $f'(K)=2(p-1)$, we must choose $\mu$ so that
$$
\frac{p}{p-1}=\frac{f'(0)}{f'(K)}=\frac{Z\sqrt{2-\mu^2}}{-Z\exp\left(-\frac{\pi \mu}{\sqrt{2-\mu^2}}\right)\sqrt{2-\mu^2}} \implies \mu = \pm\sqrt{\frac{2\log^2(p/(1-p))}{\log^2(p/(1-p)) + \pi^2}},
$$
which agrees with the calculation from Proposition \ref{travellingWave}. Finally, in Section \ref{speed}, using Proposition \ref{travellingWave}, we prove that the limiting asymptotic speed of the $(N,p)$-BBM is in fact this travelling wave speed $\mu=c=\pm\sqrt{\frac{2\log^2(p/(1-p))}{\log^2(p/(1-p))+\pi^2}}$.

\section{Some properties of the $(N,p)$-BBM} \label{properties}
To prove certain properties of the $(N,p)$-BBM it will be helpful to give a specific construction in terms of Brownian motions, a Poisson process, and uniform and Bernoulli random variables. From this construction it will be easy to prove subsequent properties. In what follows, we will write the initial condition of the $(N,p)$-BBM as a vector $\nu=(\nu_1,\ldots,\nu_N)\in \B{R}^N$, describing the positions of the $N$ particles on $\B{R}$. However we will often think of this instead as a probability distribution which is the normalised sum of atoms $\frac1N \sum_{i=1}^N \delta_{\nu_i}$.

Therefore let $B=((B^j(t))_{t\geq 0})_{1\leq j\leq N}$ be a family of $N$ i.i.d. Brownian motions and $Q=(Q(t))_{t\geq 0}$ a Poisson process of rate $N$. Let $I=(I_i)_{i\in \B{N}}$ be an i.i.d. sequence of random variables which are uniform on $\{1,\ldots,N\}$ and $P=(P_i)_{i\in \B{N}}$ an i.i.d. sequence of Bernoulli random variables each with parameter $p$. Then for any initial configuration $\nu \in \B{R}^N$ we can write the $(N,p)$-BBM $(X^{N,p}(t))_{t\geq 0}$ as a function
$$(X^{N,p}(t))_{t\geq 0}=\Theta^N \left(\Xi(\nu,B,Q,I,P,t)_{t\geq 0}\right),$$ 
where $\Theta^N:\B{R}^N \to \B{R}^N$ is the function which orders components of $\B{R}^N$-valued vectors.

This is constructed as follows. Firstly, we define $\Xi(\nu,B,Q,I,P,0)=\nu$. Now define $\tau_0=0$, and let $\tau_1<\tau_2<\cdots$ be the discontinuities of the Poisson process $Q$, which will be the successive branching times of the process. Suppose for induction that $\Xi$ is defined up to the time $\tau_i$. Thence we will define $\Xi$ for $t\in (\tau_i,\tau_{i+1}]$. Let $\Xi_j(\nu,B,Q,I,P,\tau_i)$ denote the $j$\textsuperscript{th} smallest particle of $\Xi(\nu,B,Q,I,P,\tau_i)$, breaking ties arbitrarily. Then for $t\in [\tau_i,\tau_{i+1})$, drive the particle at this location by the Brownian motion $(B^j(t))_{t\in [\tau_i,\tau_{i+1})}$. So the particles of $\Xi(\nu,B,Q,I,P,t)$ are at locations $\{\Xi_j(\nu,B,Q,I,P,\tau_i)+B^j(t)-B^j(\tau_i):j=1,2,\ldots,N\}$. Now we define the operator $k:\B{R}^N\times \{1,\ldots,N\}\times \{0,1\}\to\B{R}^N$ by:
$$k(v,i,q)=\begin{cases}
        (v_2,v_3,\ldots,v_{i-1},v_i,v_i,v_{i+1},\ldots,v_N) & \text{ if }q=1, \\
        (v_1,v_2,\ldots, v_{i-1},v_i,v_i,v_{i+1},\ldots,v_{N-1}) & \text{ if }q=0.
\end{cases}$$

Then note that if $i$ is a uniform random variable and $q$ is a $Ber(p)$ random variable, then this operator describes the branching/selection mechanism in the $(N,p)$-BBM process; a uniformly chosen particle branches into two and either the leftmost or rightmost particle is simultaneously killed (depending on the Bernoulli random variable). Thus finally, having defined the function $\Xi$ up to the time $\tau_{i+1}-$, we define
$$\Xi(\nu,B,Q,I,P,\tau_{i+1})=k\left(\Theta^N(\Xi(\nu,B,Q,I,P,\tau_{i+1}-),I_{i+1},P_{i+1}\right) .$$

Continuing inductively, we can define $\Xi(\nu,B,Q,I,P,t)$ and thus $X^{N,p}(t)=\Theta^N(\Xi(\nu,B,Q,I,P,t))$ for all $t$. From this description of the particle system the monotonicity property will clearly follow. 

Before proving monotonicity, let us first define what we mean by \textit{order} for an $N$ particle system:

\vspace{1em}

\begin{defn}
    Let $A=\{A_1,\ldots,A_M\}\in \mathbb{R}^M$ and $B=\{B_1,\ldots,B_N\}\in \B{R}^N$. We say that $A\oop B$ if and only if for all $x\in \B{R}$:
    $$|\{i:A_i \geq x\}| \leq |\{j:B_j \geq x\}|,$$ 
   where $|S|$ represents the cardinality of set $S$.
\end{defn}

From this definition, it naturally follows that if $A=\{A_1,\ldots, A_N\}$ and $B=\{B_1,\ldots,B_N\}$, then $A\oop B \iff -B\oop -A$, where $-A=\{-A_1,\ldots,-A_N\}$.

\vspace{1em}

\begin{propn} \label{monotonicity}
    \textit{(Monotonicity)} The $(N,p)$-BBM is monotonic in the sense that if $(X^{N,p}(t))_{t\geq 0}$ and $(\tilde{X}^{N,p}(t))_{t\geq 0}$ are both $(N,p)$-BBM processes with $X^{N,p}(0)=\nu \oop \tilde{\nu} =\tilde{X}^{N,p}(0)$, then there exists a coupling of the processes such that $X^{N,p}(t)\oop \tilde{X}^{N,p}(t)$ for all $t\geq 0$ almost surely.
\end{propn}

\begin{proof}
    Define $(X^{N,p}(t))_{t\geq 0}=(\Theta^N(\Xi(\nu,B,Q,I,P,t)))_{t\geq 0}$ and $(\tilde{X}^{N,p}(t))_{t\geq 0}=(\Theta^N(\Xi(\tilde{\nu},B,Q,I,P,t)))_{t\geq 0}$. Now we are given that $\nu\oop \tilde{\nu}$, so $\Xi(\nu, B, Q, I, P, \tau_0)\oop \Xi(\tilde{\nu},B,Q,I,P,\tau_0)$, and suppose for induction that
    $$\Xi(\nu,B,Q,I,P,\tau_i)\oop \Xi(\tilde{\nu},B,Q,I,P,\tau_i).$$

    This implies that $\Xi_j(\nu,B,Q,I,P,\tau_i)\leq \Xi_j(\tilde{\nu},B,Q,I,P,\tau_i)$ for $j=1,2,\ldots,N$ and therefore as the particles of $X^{N,p}$ and $\tilde{X}^{N,p}$ at these locations respectively are both driven by the same Brownian motion $B^{j}$, thus it is immediate that $\Xi_j(\nu,B,Q,I,P,t)\leq \Xi_j(\tilde{\nu},B,Q,I,P,t)$ for $j=1,2,\ldots,N$ and $t<\tau_{i+1}$. Then the result follows from the simple observation that when $v\oop \tilde{v}$ are ordered vectors in $\B{R}^N$, then $k(v,i,q)\oop k(\tilde{v},i,q)$ for any $i=1,2,\ldots,N$ and $q\in \{0,1\}$, as this implies that $\Xi(\nu,B,Q,I,P,\tau_{i+1})\oop \Xi(\tilde{\nu},B,Q,I,P,\tau_{i+1})$ and therefore by induction this ordering remains true for all $t\geq 0$. 
\end{proof}

Next we prove that any time discretisation of the $X^{N,p}$, as viewed from its leftmost particle, is positive Harris recurrent. First we recall the definition of positive Harris recurrence. 

\vspace{1em}

\begin{defn} (Definition 2.2 of \cite{athreyaNey})
    A discrete time stochastic process $(X_n)_{n\geq 0}$ on a measurable space $(S, \mathfrak{S})$ is called \textbf{Harris recurrent} if there exists an element $A\in \mathfrak{S}$ and a probability measure $\phi$ on $A$ such that $\B{P}(X_n \in A \text{ for some }n\geq 1|X_0=x)=1$ for all $x\in S$, and $\B{P}(X_1\in E|X_0=x)\geq \lambda \phi(E)$ for all $x\in A$ and $E\subseteq A$. Further, it is called \textbf{positive Harris recurrent} if
    $$\sup_{x\in A}\B{E}_x[\tau_A]<\infty,$$ where $\tau_A:=\inf\{n\geq 1:X_n\in A\}$ is the hitting time of the set $A$.
\end{defn}

Now let $X^{N,p}$ be an $(N,p)$-BBM process, and let $X_1^{N,p}(t)$ be the location of the leftmost particle of $X^{N,p}$ at time $t$. Then we define the process $(\bar{X}^{N,p}(t))_{t\geq 0}:=(X^{N,p}(t)-X^{N,p}_1(t)\underline{1})_{t\geq 0}$ to be the process $X^{N,p}$ as viewed from its leftmost particle. 

\vspace{1em}

\begin{propn} \label{harrisRecc}
    Fix $\delta>0$. Then the stochastic process $(\bar{X}^{N,p}(n\delta ))_{n\in \B{N}}$ is positive Harris recurrent, and consequently has a unique stationary distribution to which it converges. 
\end{propn}

\begin{proof}
    Let $A=[0,1]^N$. We will first show that the random process $\bar{X}^{N,p}$ returns to $A$ in finite expected time. Consider our construction
    $$(X^{N,p}(t))_{t\geq 0}=\Theta^N\left(\left(\Xi\left(\nu, B,Q,I,P,t\right)\right)\right)_{t\geq 0}$$
    of $X^{N,p}$ as a function of random variables. Now consider the event $G_n=E_n\cap F_n$ where we define the events
    \begin{align*}
        E_n:&=\{Q(n\delta)-Q((n-1)\delta)=N-1\}
        \cap\{I_{i}=1\text{ and }P_i=0\text{ for }Q((n-1)\delta)< i\leq Q(n\delta)\} \\
        F_n:&=\left\{\sup_{(n-1)\delta\leq s<t\leq n\delta}|B^j(s)-B^j(t)|\leq 1/4N\text{ for }j=1,2,\ldots N\right\}
    \end{align*}
    So $E_n$ is the event that in the interval $[(n-1)\delta,n\delta]$, the leftmost particle of $X^{N,p}(t)$ branches $N-1$ times and at each of these branching events the rightmost particle is killed, and $F_n$ is the event that the biggest displacement difference between two points of a Brownian motion $B^j$ in the time interval $[(n-1)\delta,n\delta]$ is at most $1/2N$.
    
    Now we note that by our construction, under the event $E_n$, every particle trajectory of a particle of $X^{N,p}$ is made up of at most $N$ segments of Brownian paths $(B^j(u))_{s\leq u\leq t}$, since between discontinuity times of $Q$, a particle is driven by a single Brownian motion $B^j$, and under the event $E_n$ there are $N-1$ discontinuities of $Q$ in $[(n-1)\delta,n\delta]$. Therefore the event $G_n$ ensures that the distance between $X^{N,p}_i(n\delta)$ and $X_1^{N,p}(n\delta)$ is at most $1/4N \times N \times 2 = 1/2$. Thus on the event $G_n$, all particles in $X^{N,p}(n\delta)$ are within distance $1$ of each other, and therefore $\bar{X}^{N,p}(n\delta)\in A$. Therefore as the event $G_n$ is independent of the configuration of particles at time $(n-1)\delta$, thus $\B{P}(G_n)=q_{\text{hit}}>0$ is the same non-zero constant for each $n$. Thus the hitting time of $A$, $\tau_A:=\inf\{t:X^{N,p}(t)\in A\}$, from any configuration is bounded above by a geometric random variable with parameter $q_{\text{hit}}$, and hence $\B{E}_x[\tau_A]\leq q_{\text{hit}}^{-1}<\infty$ for all $x\in \C{S}$. Thus to prove that $(X(n\delta))_{n\geq 0}$ is positive Harris recurrent, it just remains to prove that there exists a constant $\lambda$ and a probability measure $\phi$ such that $\B{P}_x(\bar{X}^{N,p}(\delta)\in S)\geq \lambda \phi(S)$ for any $x=(0,x_2,x_3,\ldots,x_N)\in A$ and $S=\{0\}\times S_2 \times \cdots \times S_N \subseteq A$. 
    
    We will prove this inequality by taking $\phi$ to be the Lebesgue measure on the cube $\{0\}\times [0,1]^{N-1}$, so that, for example, $Leb(S)=\prod_{i=2}^N \int_{S_i}dy_i$. Now consider the event $E:=\{Q(\delta)=0, B^1(\delta)-B^1(0)\in [-1,0]\}$ that during $[0,\delta]$ no particle branches and that the leftmost particle moves at most $1$ unit to the left. This event occurs with a non-zero probability, say $q_E$, independent of the initial condition $x$. Thus we have that
    \begin{align}
        \nonumber \B{P}_x(\bar{X}^{N,p}(\delta)\in S) &\geq q_E \B{P}_x(\bar{X}^{N,p}(\delta)\in S|E) \\
        \nonumber &= q_E \B{P}_x(X^{N,p}_i(\delta) - X^{N,p}_1(\delta) \in S_i \text{ for } i=2,3,\ldots,N|E) \\
        \nonumber &\geq q_E \inf_{x_1 \in [-1,0]}\B{P}_x(X^{N,p}_i(\delta)-x_1 \in S_i \text{ for } i=2,3,\ldots,N|E) \\
        \nonumber &\geq q_E \inf_{x_1 \in [-1,0]}\prod_{i=2}^N \int_{S_i}(2\pi \delta)^{-1/2}e^{-\frac{(y_i - x_i+x_1)^2}{2\delta}}dy_i \\
        \label{harrisExpBound} &\geq q_E \inf_{x_1\in [-1,0]}\prod_{i=2}^N \int_{S_i} (2\pi \delta)^{-1/2}e^{-2/\delta}dy_i \\
        \nonumber &=q_E (2\pi \delta)^{-N/2}e^{-2N/\delta}\inf_{x_1\in [-1,0]} Leb(S) = q_E (2\pi \delta)^{-N/2}e^{-2N/\delta} Leb(S),
    \end{align}
    where the inequality (\ref{harrisExpBound}) follows from the fact that $x_1\in [-1,0]$, $x_i\in [0,1]$, and $y_i\in S_i\subseteq [0,1]$, so $(y_i-x_i+x_1)^2 \leq 4$. Therefore we have shown that the exists a constant $\lambda$ and measure $\phi$ such that $\B{P}_x(X^{N,p}(\delta)\in S)\geq \lambda \phi(S)$ for any $x\in A$ and $S\subseteq A$. Thus $(\bar{X}^{N,p}(n\delta))_{n\in\B{N}}$ is positive Harris recurrent. Therefore by Theorems 4.1 and 6.1 of \cite{athreyaNey}, $(\bar{X}^{N,p}(n\delta))_{n\in \B{N}}$ has a unique invariant distribution to which it converges. \end{proof}

Next we will use Kingman's subadditive ergodic theorem to prove that the $(N,p)$-BBM has an asymptotic velocity. Let us first recall the following definitions. The first is from Kingman's survey of subadditive processes \cite{sipKingman}:

\vspace{1em}

\begin{defn}
    A family $\C{Y}=(Y(s,t):s,t\in [0,\infty), s<t)$ is called a \textit{continuous-parameter subadditive process} if it satisfies the following properties:
    \begin{enumerate}[({A}1)]
        \item $Y(s,u)\leq Y(s,t) + Y(t,u)$ for all $s<t<u$.
        \item For all $T\in [0,\infty)$, the distribution of $Y(s+T,t+T)$ is the same as the distribution of $Y(s,t)$. 
        \item For all $t\in [0,\infty)$, the expectation $\B{E}[Y(0,t)]$ exists and satisfies $\B{E}[Y(0,t)]\geq - \gamma_0 t$ for some $\gamma_0 > - \infty$.
    \end{enumerate}
\end{defn}

The second comes from the textbook of Todorovic \cite{todorovic}, Definition 1.7.1:

\vspace{1em}

\begin{defn}
    Let $T\subset \B{R}$ be an interval. A stochastic process $\{\xi(t):t\in T\}$  is \textit{separable} if there exists a countable dense subset $D\subset T$ and null set $\Lambda\subset \Omega$ such that for any open  interval $I$ and closed set $C$
    $$\{\omega:\xi(t,\omega)\in C, t\in I\cap D\}\setminus \{\omega:\xi(t,\omega)\in C, t\in I\}\subset \Lambda.$$
\end{defn}

\vspace{1em}

\begin{theorem}[Kingman, Theorem 4, Section 1.4 of \cite{sipKingman}] \label{kingmanSub}
    If $\C{Y}$ is a continuous-parameter subadditive process then the limit $\gamma:=\lim_{t\to\infty}\B{E}[Y(0,t)]/t$ exists. If, in addition, $Y$ is cadlag (and hence separable) and $\B{E}[\sup_{s<t,s,t\in I}|Y(s,t)|]<\infty$ for any bounded interval $I\subseteq [0,\infty)$, then $\lim_{t\to\infty}Y(0,t)/t$ exists almost surely and in $L^1$ and equals $\gamma$. 
\end{theorem}

\vspace{1em}

\begin{propn} \label{kingmanSpeed}
    For any bounded initial configuration $\nu$, the limits $\lim_{t\to\infty}\B{E}[X_1^{N,p}(t)]/t$ and $\lim_{t\to\infty}\B{E}[X_N^{N,p}(t)]/t$ exist and are equal, and the limits $\lim_{t \to \infty}X_1^{N,p}(t)/t$ and $\lim_{t\to\infty}X_N^{N,p}(t)/t$ exist and are equal almost surely and in $L^1$.  
\end{propn}

\begin{proof}
    In order to prove that the limit $X_N^{N,p}(t)/t$ exists, we will construct a continuous parameter subadditive process $\C{Y}$ such that $Y(0,t)=X^{N,p}_N(t)$. Initially, we will suppose that the initial configuration is a constant vector $\nu(a):=(a,a,\ldots,a)$, a condition which we will relax later. Then for any $t\geq s\geq 0$ let us define:
    $$Z(s,t):=\Xi\left(\nu,((B^j(s+u))_{u\geq 0})_{1\leq j\leq N}, (Q(s+u))_{u\geq 0}, (I_i)_{i>Q(s)}, (P_i)_{i>Q(s)}, t-s\right).$$

This defines a stochastic flow of the $(N,p)$-BBM process. Simply, this is the $(N,p)$-BBM constructed so that for $s < t < u$, there is a natural coupling between the processes $(Z(s,w))_{w\geq u}$ and $(Z(t,w))_{w\geq u}$. Now let $Z_N(s,t)$ be the location of the rightmost particle of $Z(s,t)$, and consider $Y(s,t):=Z_N(s,t)-Z_N(s,s)$. Therefore by definition, both $Y(s,t)$ and $Y(s+T,t+T)$ describe the position of the rightmost particle of an $(N,p)$-BBM started from initial configuration $(0,0,\ldots,0)$ at $t-s$ units of time ago, and therefore have the same law. Thus $\C{Y}$ satisfies property \textit{(A2)}.

    Secondly, property \textit{(A3)} is easy to prove. Clearly the $(N,p)$-BBM can be embedded inside a branching Brownian motion (BBM) started from $N$ particles with configuration $\nu$. Now if $(\{X_u(t), u\in \C{N}(t)\})_{t\geq 0}$ is a branching Brownian motion with particles at locations $X_u(t)$ and particle set $\C{N}(t)$ at time $t$ started from a single particle at the origin, then it is well known (see \cite{lalleySellke} for example) that $D(t):=\sum_{u\in \C{N}(t)}\exp(\sqrt{2}X_u(t) - 2t)$ is a martingale (called the derivative martingale). Then by Jensen's inequality
    $$\exp(\sqrt{2} \B{E}[\max_u X_u(t)] - 2t) \leq \B{E}\left[\exp(\sqrt{2}\max_u X_u(t) - 2t)\right] \leq \B{E}[D(t)]=\B{E}[D(0)]=1,$$
    and therefore $\B{E}[\max_u X_u(t)] \leq \sqrt{2}t$ and by symmetry $\B{E}[\min_u X_u(t)] \geq -\sqrt{2}t$. Therefore we can lower bound the expected position of the leftmost particle of a BBM started from $N$ particles at the origin by $-N\sqrt{2}t$, therefore proving the property (\textit{A3}).

    To prove the sub-additivity property (\textit{A1}), we note that trivially $Z(s,t)\oop \nu(Z_N(s,t))$ (recalling that $\nu(a)=(a,\ldots,a)$). Then by the definition of $Z$, we know that $(Z(s,t+v))_{v\geq 0}$ and $(Z(t,t+v)-Z(t,t)+\nu(Z_N(s,t)))_{v\geq 0}$ are two coupled $(N,p)$-BBM processes, which have initial configurations $Z(s,t)$ and $\nu(Z_N(s,t))$ respectively. Therefore by monotonicity (Proposition \ref{monotonicity}), for any $v\geq 0$, $Z(s,t+v)\oop Z(t,t+v)-Z(t,t)+\nu(Z_N(s,t))$. Hence for $v=u-t$
    $$Y(s,u)=Z_N(s,u)-Z_N(s,s)\leq Z_N(t,u)-Z_N(t,t)+Z_N(s,t)-Z_N(s,s)=Y(t,u)+Y(s,t),$$
    thus proving the sub-additivity property \textit{(A1)}. Therefore $\C{Y}=(Y(s,t):s,t\in [0,\infty),s<t)$ is a continuous parameter subadditive process.

    Finally in order to apply Theorem \ref{kingmanSub}, it just remains to show that for any bounded interval $I\subseteq [0,\infty)$ we have $E[\sup_{s<t,s,t\in I} |Y(s,t)|]<\infty$. So fix $I=[a,b]\subseteq [0,\infty)$. Again, recall that the $(N,p)$-BBM can be considered as being embedded into $N$ independent BBMs, each started from a single particle at the origin, and hence $\sup_{s<t,s,t,\in I}|Y(s,t)|$ is bounded by the maximum up to time $b-a$ of $N$ independent BBMs, all started from $0$ (since $Y(s,s)=(0,\ldots,0)$ when $\nu = (a,\ldots,a)$). Thus by the many-to-one lemma (see for example Lemma 2.1 in \cite{kim}), we certainly have that:
    \begin{align}\label{boundForKingman}\B{E}\left[\sup_{s<t,s,t\in I}|Y(s,t)|\right]\leq N\B{E}\left[\sup_{u\in \C{N}(b-a)}\sup_{0\leq t\leq b-a}|X_u(t)|\right]=Ne^{b-a}\B{E}\left[\sup_{0\leq t\leq b-a}|B(t)|\right]<\infty,\end{align}
    where $(B(t))_{t\geq 0}$ is a standard Brownian motion and $\{X_u(t):u\in \C{N}(t)\}_{t\geq 0}$ is a branching Brownian motion started with a single particle at the origin.

    Therefore we can apply Theorem \ref{kingmanSub} to the process $\C{Y}$ to conclude that $\lim_{t\to\infty}\B{E}[Y(0,t)]/t$ exists, and $\lim_{t\to\infty}Y(0,t)/t$ exists and is equal to $\lim_{t\to\infty}\B{E}[Y(0,t)]/t$ almost surely and in $L^1$. Therefore we have that for an $(N,p)$-BBM process $(X^{N,p}(t))_{t\geq 0}$ with initial condition $\nu(a)=(a,\ldots,a)$, the limit $\gamma^{N,p}_N:=\lim_{t\to\infty}\B{E}[X_N^{N,p}(t)]/t$ exists and the limit $\lim_{t\to\infty}X_N^{N,p}(t)/t$ exists and is equal to $\gamma_N^{N,p}$ almost surely and in $L^1$. Then to generalise the result to any bounded initial condition $\tilde{\nu}$, say with each element of $\nu$ bounded between $m$ and $M$, we can couple the process $X^{N,p}(t)$ with initial condition $\tilde{\nu}$ by monotonicty so that 
    $$(\Xi(\nu(m),B,Q,I,P,t))_{t\geq 0}\oop(\Xi(\tilde{\nu},B,Q,I,P,t))_{t\geq 0}\oop (\Xi(\nu(M),B,Q,I,P,t))_{t\geq 0},$$
    thus proving that the asymptotic velocity of $X^{N,p}_N(t)$ exists for any initial condition with bounded support.    

    By symmetry, we can similarly apply Theorem \ref{kingmanSub} to the continuous-parameter subadditive process $\tilde{\C{Y}}=(-Z_1(s,t)+Z_1(s,s):s,t\in [0,\infty),s<t)$, and therefore $\gamma_1^{N,p}:=\lim_{t\to\infty}\B{E}[X_1^{N,p}(t)]/t$ exists and $\lim_{t\to\infty}X_1^{N,p}(t)/t$ exists and equals $\gamma_1^{N,p}$ almost surely and in $L^1$. 

    It just remains to show that $\gamma_1^{N,p}=\gamma_N^{N,p}$. Recall the event $G_n$ that we defined in the proof of Proposition \ref{harrisRecc}, and define the random times $\tau_0=0$ and subsequently $\tau_{i+1}:=\inf\{n>\tau_i:n\in \B{N},\; G_n\text{ occurs}\}$. Recall that, at the stopping times $\tau_i$, we have $|X_N^{N,p}(\tau_i)-X_1^{N,p}(\tau_i)|\leq 1$. Note that as $Q$ has independent increments, and $I$ and $P$ are i.i.d families of random variables, thus $\tau_{i+1}-\tau_i$ are i.i.d and have finite mean. 
    
    Therefore since $\tau_i\to\infty$ as $i\to\infty$, it follows that $$0=\lim_{i\to\infty}\B{E}[X_N^{N,p}(\tau_i)-X_1^{N,p}(\tau_i)]/\tau_i = \lim_{t\to\infty}\B{E}[X_N^{N,p}(t)-X_1^{N,p}(t)]/t =\gamma_N^{N,p}-\gamma_1^{N,p},$$ 
    therefore $\gamma_N^{N,p}=\gamma_1^{N,p}$ as required. 
\end{proof}

\section{Coupling particle systems with two-sided selection} \label{twoSidedCoupling}
The aim of this section is to construct upper and lower bounding particle systems which can be coupled to $X^{N,p}$. In particular, for any $\delta>0$, we want to construct particle systems $X^{N,p,\delta,+}$ and $X^{N,p,\delta,-}$ such that we have the order $X^{N,p,\delta,-}(k\delta) \oop X^{N,p}(k\delta) \oop X^{N,p,\delta,+}(k\delta)$ at all discrete times $k\delta$, $k\in \B{N}$. Moreover, we will construct the upper and lower bounding systems such that $X^{N,p,\delta,\pm}$ behave as branching Brownian motions for all intervals $(k\delta, k\delta + \delta)$.

Our couplings build on the ideas used in \cite{hydroNBBM}, which constructs couplings for the $N$-BBM by discretising time to $\delta \B{N}$ and performing selection only at the discrete times. To do this, we will firstly construct the $(N,p)$-BBM as a process embedded in a branching Brownian motion (BBM) which we describe by colouring. 

\paragraph{Construction of $X$:} Consider a branching Brownian motion with deletion, $\C{B}^{N,\nu}=(\C{B}^{N,\nu}(t))_{t\geq 0}=\{B_u(t):u\in \C{N}^{N,\nu}(t)\}_{t\geq 0}$, starting from $N$ particles, whose initial locations are given by the vector $\nu\in \B{R}^N$. We define $X=X^{N,p}$ as a function $\Gamma$ of the BBM $\C{B}^{N,\nu}$ and a sequence of i.i.d. $\text{Bernoulli}(p)$ random variables, $P=(P_i)_{i\in \B{N}}$. We initially set $X(0)=\Gamma(\C{B}^{N,\nu},P,0)=\Theta^N(\nu)$, where $\Theta^N:\B{R}^N\to\B{R}^N$ is the function which puts the $N$ elements of a vector into increasing order. 

The particles of the BBM will be coloured either blue or yellow, and initially we set the colour of all particles to yellow. Each time a particle branches, it gives birth to two particles of the same colour. Additionally, at the $i$\textsuperscript{th} branching event of a yellow particle, if $P_i=1$, we simultaneously recolour the leftmost yellow particle blue, and if $P_i=0$, we simultaneously delete the rightmost yellow particle. Thus in the coloured process $\C{B}^{N,\nu}$ there are at all times exactly $N$ yellow particles. For $t>0$, we define $X(t)=\Gamma(\C{B}^{N,\nu},P,t)$ to be the vector of the locations of the $N$ yellow particles of $\C{B}^{N,\nu}(t)$ in increasing order. 

Now we will consider coupling $\C{B}^{N,\nu}$ to a branching Brownian motion starting from a different initial condition and with fewer than $N$ particles. So let $\tilde{\nu}\in \B{R}^N$ be such that $\nu\oop \tilde{\nu}$. Then let $\C{B}^{m,\tilde{\nu}}=(\C{B}^{m,\tilde{\nu}}
(t))_{t\geq 0}$ be a branching Brownian motion starting with $m\leq N$ particles, at the locations of the leftmost $m$ particles of $\tilde{\nu}$. Again, we initially colour the $m$ particles of $\C{B}^{m,\tilde{\nu}}$ yellow, although at later times we will also colour some particles of $\C{B}^{m,\tilde{\nu}}$ blue. Now for a coloured BBM $\C{B}$, let us define $LY(\C{B}(t),k)$ (resp. $RY(\C{B}(t),k)$) to be the vector of locations of the $k$ leftmost (resp. rightmost) yellow particles of $\C{B}(t)$. Therefore $\nu\oop \tilde{\nu}$ implies that we have initially $LY(\C{B}^{N,\nu}(0),m)\oop RY(\C{B}^{m,\tilde{\nu}}(0),m)$. 

We will couple $\C{B}^{N,\nu}$ to $\C{B}^{m,\tilde{\nu}}$ so that $LY(\C{B}^{N,\nu}(t),k)\oop RY(\C{B}^{m,\tilde{\nu}}(t),k)$ for all $t\geq 0$, for a suitable choice of $k$. Our coupling is constructed as follows. Define $\tau_0=0$ and let $\tau_1<\tau_2<\cdots$ be the branching times of yellow particles in the process $\C{B}^{N,\nu}$. Consider the time interval $[\tau_i,\tau_{i+1}]$. Let $y^m(t)$ denote the number of yellow particles in $\C{B}^{m,\tilde{\nu}}(t)$, and suppose that $y^m(\tau_i)\leq N-1$. Suppose for induction that $LY(\C{B}^{N,\nu}(\tau_i),y^m(\tau_i))\oop RY(\C{B}^{m,\tilde{\nu}}(\tau_i),y^m(\tau_i))$. Let $\C{B}^{N,\nu}_j(\tau_i)$ denote the position of the $j$\textsuperscript{th} leftmost yellow particle of $\C{B}^{N,\nu}$, and call $(B^{i,j}(t))_{t\geq 0}$ the Brownian motion driving this particle over the time interval $[\tau_i,\tau_{i+1}]$. Then we will drive the $j$\textsuperscript{th} leftmost yellow particle of $\C{B}^{m,\tilde{\nu}}(\tau_i)$ by the same $(B^{i,j}(t))_{t\geq 0}$ on $[\tau_{i},\tau_{i+1}]$. Moreover, on $(\tau_i,\tau_{i+1})$, we do no recolouring of particles. Therefore it is clear that:
\begin{align*}LY(\C{B}^{N,\nu}(\tau_i),y^m(\tau_i))&\oop RY(\C{B}^{m,\tilde{\nu}}(\tau_i),y^m(\tau_i))\implies \{\C{B}^{N,\nu}_j(\tau_i): j\leq y^m(\tau_i)\}\oop \{\C{B}^{m,\tilde{\nu}}_j(\tau_i):j\leq y^m(\tau_i)\} \\ & \implies \{\C{B}^{N,\nu}_j(\tau_i)+B^{i,j}(t-\tau_i):j\leq y^m(\tau_i)\}\oop \{\C{B}^{m,\tilde{\nu}}_j(\tau_i)+B^{i,j}(t-\tau_i):j\leq y^m(\tau_i)\},\end{align*}
and so $LY(\C{B}^{N,\nu}(t),y^m(\tau_i))\oop RY(\C{B}^{m,\tilde{\nu}}(t),y^m(\tau_i))$ for $t\in [\tau_i,\tau_{i+1})$. Now at the branching time $\tau_{i+1}$, if the particle at location $\C{B}^{N,\nu}_j(\tau_{i+1}-)$ branches, then we also branch the particle at location $\C{B}^{m,\tilde{\nu}}_j(\tau_{i+1}-)$ if $j\leq y^m(\tau_i)$, and otherwise no branching of particles occurs in $\C{B}^{m,\tilde{\nu}}$ at time $\tau_{i+1}$. 
Recall that in $\C{B}^{N,\nu}(\tau_{i+1})$, if $P_{i+1}=1$, we recolour the leftmost yellow particle blue, and if $P_{i+1}=0$, we delete the rightmost yellow particle. Then for $\C{B}^{m,\tilde{\nu}}$, if $P_{i+1}=1$, we will recolour the leftmost yellow particle blue, and if $P_{i+1}=0$, we will not change $\C{B}^{m,\tilde{\nu}}$. We claim that, in any case, $LY(\C{B}^{N,\nu}(\tau_{i+1}),y^m(\tau_{i+1}))\oop RY(\C{B}^{m,\tilde{\nu}}(\tau_{i+1}),y^m(\tau_{i+1}))$

Let us consider what happens in each case. Firstly, note that $$LY(\C{B}^{N,\nu}(\tau_{i+1}-),y^m(\tau_{i+1}-))\oop RY(\C{B}^{m,\tilde{\nu}}(\tau_{i+1}-),y^m(\tau_{i+1}-))$$ implies that $\C{B}^{N,\nu}_j(\tau_{i+1}-)\leq \C{B}^{m,\tilde{\nu}}_j(\tau_{i+1}-)$ for each $j=1,2,\ldots,y^m(\tau_{i+1}-)$. Therefore if the $i$\textsuperscript{th} leftmost yellow particles of $\C{B}^{N,\nu}$ and $\C{B}^{m,\tilde{\nu}}$ both branch, then (before deletion and recolouring), the leftmost $y^m(\tau_{i+1}-)+1$ yellow particles of $\C{B}^{N,\nu}$ still lie to the left of the $y^m(\tau_{i+1}-)+1$ yellow particles of $\C{B}^{m,\tilde{\nu}}$, and if no particle of $\C{B}^{m,\tilde{\nu}}$ branches, then (before deletion and recolouring), the leftmost $y^m(\tau_{i+1}-)$ yellow particles of $\C{B}^{N,\nu}$ trivially must lie to the left of the $y^m(\tau_{i+1}-)$ yellow particles of $\C{B}^{m,\tilde{\nu}}$. Now we consider how deletion and recolouring affect the ordering. Define $h$ to be the function which takes as input a vector of any length and removes the smallest element (for example $h((4,1,3,2))=(4,3,2)$). Note that if $P_{i+1}=1$, then recolouring the leftmost yellow particle of $\C{B}^{N,\nu}(\tau_{i+1}-)$ or $\C{B}^{m,\tilde{\nu}}(\tau_{i+1}-)$ blue is equivalent to applying $h$ to the set of yellow particles in $\C{B}^{N,\nu}(\tau_{i+1}-)$ or $\C{B}^{m,\tilde{\nu}}(\tau_{i+1}-)$ respectively. Then if $\underline{u},\underline{v}\in \B{R}^n$ and $\underline{u}\oop \underline{v}$, we have $h(\underline{u})\oop h(\underline{v})$. Therefore in the case that the leftmost yellow particles of $\C{B}^{N,\nu}(\tau_{i+1}-)$ and $\C{B}^{m,\tilde{\nu}}(\tau_{i+1}-)$ are coloured blue (ie. if $P_{i+1}=1)$, it follows that
\begin{align*}LY(\C{B}^{N,\nu}(\tau_{i+1}),y^m(\tau_{i+1}))&=h(LY(\C{B}^{N,\nu}(\tau_{i+1}-),y^m(\tau_{i+1}-)))\\
&\oop h(RY(\C{B}^{m,\tilde{\nu}}(\tau_{i+1}-),y^m(\tau_{i+1}-)))=RY(\C{B}^{m,\tilde{\nu}}(\tau_{i+1}),y^m(\tau_{i+1})).\end{align*}
On the other hand, if $P_{i+1}=0$, and we delete the rightmost yellow particle of $\C{B}^{N,\nu}(\tau_{i+1})$, then as $y^m(\tau_{i+1})\leq N$ and $\C{B}^{N,\nu}$ always has exactly $N$ yellow particles, thus $LY(\C{B}^{N,\nu}(\tau_{i+1}),y^m(\tau_{i+1}))$ is unchanged by the deletion. Therefore in the case $P_{i+1}=0$, it immediately follows that $LY(\C{B}^{N,\nu}(\tau_{i+1}),y^m(\tau_{i+1}))\oop RY(\C{B}^{m,\tilde{\nu}}(\tau_{i+1}),y^m(\tau_{i+1}))$. Therefore, by induction, we have proven that, up until and including the first branching event at which we have $y^m(\tau_\ell)=N$, there is a coupling of $\C{B}^N$ and $\C{B}^m$ such that $LY(\C{B}^{N,\nu}(\tau_\ell),y^m(\tau_\ell))\oop RY(\C{B}^{m,\tilde{\nu}}(\tau_\ell),y^m(\tau_\ell))$. 

Of course, there are also blue particles in $\C{B}^{N,\nu}$ and $\C{B}^{m,\tilde{\nu}}$ whose motion and branching we have not yet described. We will couple the blue particles of $\C{B}^{N,\nu}$ and $\C{B}^{m,\tilde{\nu}}$ in the following way: if a blue particle of $\C{B}^{N,\nu}$ branches, and at the time of branching is the $j$\textsuperscript{th} leftmost blue particle, we simultaneously branch the $j$\textsuperscript{th} leftmost blue particle of $\C{B}^{m,\tilde{\nu}}$. In this way, the number of blue particles in $\C{B}^{N,\nu}(t)$ is always at least the number of blue particles in $\C{B}^{m,\tilde{\nu}}(t)$. In fact, we can observe that the only instance in which we have \textit{strictly more} blue particles in $\C{B}^{N,\nu}$ than in $\C{B}^{m,\tilde{\nu}}$ is if there are no remaining yellow particles in $\C{B}^{m,\tilde{\nu}}$ to turn blue; that is, if \textit{all} yellow particles have been turned blue. In this instance, $\C{B}^{m,\tilde{\nu}}$ can thereafter only contain blue particles. Therefore if there are any remaining yellow particles in $\C{B}^{m,\tilde{\nu}}$, we know that the number of blue particles in $\C{B}^{N,\nu}$ and $\C{B}^{m,\tilde{\nu}}$ are equal. It is not necessary to couple the Brownian motions driving the blue particles; it is sufficient just to say that they move as independent Brownian motions.

Let us take a moment here to observe that, according to our construction, the number of yellow particles in $\C{B}^{m,\tilde{\nu}}$ may decrease. This can happen at the branching time $\tau_i$ if the $j$\textsuperscript{th} leftmost yellow particle of $\C{B}^{N,\nu}$ branches with $j>y^m(\tau_i)$ (so that no branching occurs in $\C{B}^{m,\tilde{\nu}}$ at time $\tau_i$) and simultaneously the leftmost yellow particle of $\C{B}^{N,\nu}(\tau_i-)$ and $\C{B}^{m,\tilde{\nu}}(\tau_i-)$ are coloured blue (i.e. $P_i=1$). In such a case $y^m(t)$ decreases by $1$ at the time $\tau_i$. However, we can observe that once $y^m(t)\geq N$, then $y^m(t)$ can no longer decrease, since when $y^m(t)\geq N$, every time a yellow particle of $\C{B}^{N,\nu}$ branches, a yellow particle of $\C{B}^{m,\tilde{\nu}}$ also branches. 

So let $\tau_\ell$ be the first branching event at which we reached $y^m(\tau_\ell)=N$. We will now inductively define $\C{B}^{m,\tilde{\nu}}$ for times $t\geq \tau_\ell$, coupled to $\C{B}^{N,\nu}$, so that $LY(\C{B}^{N,\nu}(t),N)\oop RY(\C{B}^{m,\tilde{\nu}}(t),N)$ for $t\geq \tau_\ell$. Consider the time interval $[\tau_i,\tau_{i+1}]$, for $i\geq \ell$, and suppose for induction that $LY(\C{B}^{N,\nu}(\tau_i),N)\oop RY(\C{B}^{m,\tilde{\nu}}(\tau_i),N)$. Call $(B^{i,j}(t))_{t\geq 0}$ the Brownian motion driving the $j$\textsuperscript{th} leftmost yellow particle of $\C{B}^{N,\nu}(\tau_i)$, $\C{B}^{N,\nu}_j(\tau_i)$, on $[\tau_i,\tau_{i+1})$. Then we will also drive by $B^{i,j}(t)$ the $y^m(\tau_i)-N+j$\textsuperscript{th} leftmost yellow particle of $\C{B}^{m,\tilde{\nu}}(\tau_i)$, $\C{B}^m_{y^m(\tau_i)-N+j}(\tau_i)$ (i.e. the $j$\textsuperscript{th} leftmost out of the rightmost $N$ yellow particles) on $[\tau_i,\tau_{i+1})$. Therefore it is clear that:
\begin{align*}
    \{\C{B}^{N,\nu}_j(\tau_i):&1\leq j\leq N\}\oop \{\C{B}^{m,\tilde{\nu}}_{y^m(\tau_i)-N+j}(\tau_i):1\leq j\leq N\} \\
    & \implies \{\C{B}^{N,\nu}_j(\tau_i)+B^{i,j}(t-\tau_i):1\leq j\leq N\} \oop \{\C{B}^{m,\tilde{\nu}}_{y^m(\tau_i)-N+j}(\tau_i)+B^{i,j}(t-\tau_i): 1\leq j\leq N\},
\end{align*}
and so $LY(\C{B}^{N,\nu}(t),N)\oop RY(\C{B}^{m,\tilde{\nu}}(t),N)$ for $t\in [\tau_i,\tau_{i+1})$. Then at time $\tau_{i+1}$, if the $j$\textsuperscript{th} leftmost particle of $\C{B}^{N,\nu}(\tau_{i+1})$ branches, then we branch the $y^m(\tau_i)-N+j$\textsuperscript{th} leftmost particle of $\C{B}^{m,\tilde{\nu}}(\tau_i)$. As before, since $\C{B}^{N,\nu}_j(\tau_{i+1}-)\leq \C{B}^{m,\tilde{\nu}}_{y^m(\tau_i)-N+j}(\tau_{i+1}-)$, thus the ordering $LY(\C{B}^{N,\nu}(t),N)\oop RY(\C{B}^{m,\tilde{\nu}}(t),N)$ is preserved under branching. Simultaneously to branching, if $P_{i+1}=1$, we recolour the leftmost yellow particle of both $\C{B}^{N,\nu}(\tau_{i+1}-)$ and $\C{B}^{m,\tilde{\nu}}(\tau_{i+1}-)$ blue, and if $P_{i+1}=0$, we delete the rightmost yellow particle of $\C{B}^{N,\nu}(\tau_{i+1}-)$. By identical reasoning to before, in either case ($P_{i+1}=0$ or $P_{i+1}=1$), the ordering is preserved, and $LY(\C{B}^{N,\nu}(\tau_{i+1}),N)\oop RY(\C{B}^{m,\tilde{\nu}}(\tau_{i+1}),N)$. So by induction, it follows that $LY(\C{B}^{N,\tilde{\nu}}(t),N)\oop RY(\C{B}^{m,\tilde{\nu}}(t),N)$ for all times $t>\tau_\ell$. 

As with times $t\leq \tau_{\ell}$, when $t\geq \tau_\ell$ there are also blue particles of $\C{B}^{N,\nu}$ and blue and yellow particles of $\C{B}^{m,\tilde{\nu}}$ which we have not described the behaviour of. Again, all we require is that when the $j$\textsuperscript{th} leftmost blue particle of $\C{B}^{N,\nu}$ branches, the $j$\textsuperscript{th} leftmost blue particle of $\C{B}^{m,\tilde{\nu}}$ simultaneously branches. Since we have assumed that for $t\geq \tau_{\ell}$ we have at least $N$ yellow particles in $\C{B}^{m,\tilde{\nu}}$, thus by our earlier remark, the number of blue particles in $\C{B}^{N,\nu}(t)$ is equal to the number of blue particles in $\C{B}^{m,\tilde{\nu}}(t)$. It is not necessary to couple the driving Brownian motions of the blue particles or the yellow particles which are not among the $N$ rightmost as branching times, but just to assume that they are driven by independent Brownian motions and branch at unit rate. 

Now for times $t\geq 0$, we define $\Gamma^+(\C{B}^{m,\tilde{\nu}},P,t)$ to be the vector of positions of the particles (of any colour) of the BBM $\C{B}^{m,\tilde{\nu}}$ which starts from the leftmost $m$ particles of the initial configuration $\tilde{\nu}$ and is coupled to $\C{B}^{N,\nu}$ in the above way. Therefore we can observe that when $y^m(t)\geq N$ (that is, there are at least $N$ yellow particles in the BBM $\C{B}^{m,\tilde{\nu}}$), we have:
\begin{align} \label{moreThanNIneq}
    \Gamma(\C{B}^{N,\nu},P,t)=LY(\C{B}^{N,\nu}(t),N)\oop RY(\C{B}^{m,\tilde{\nu}}(t),N) \oop \Gamma^+(\C{B}^{m,\tilde{\nu}},P,t).
\end{align}

The equality is trivial since the process $\Gamma(\C{B}^{N,\nu},P,t)$ is precisely the vector of locations of the $N$ yellow particles in $\C{B}^{N,\nu}$. The first inequality $\oop$ we proved by induction. The second inequality $\oop$ follows from the fact that the yellow particles considered by $RY(\C{B}^{m,\tilde{\nu}}(t),N)$ are a subset of the particles of $\Gamma^+(\C{B}^{m,\tilde{\nu}},P,t)$. When defining $\Gamma^+$, we consider both blue and yellow particles because this allows us to more easily find the hydrodynamic limit of the upper bounding system, since $\Gamma^+$ is then simply a branching Brownian motion, whose hydrodynamic limit is more easily attained. 

\paragraph{Construction of $X^+$:} Fix some small parameter $\delta$ and initial conditions $\nu\oop \tilde{\nu}$, and consider the time interval $[0,\delta]$. Let the process $X$ on the interval $[0,\delta]$ be defined through the function $X(t)=\Gamma(\C{B}^{N,\nu},P,t)$, and suppose that there are $b^N(\delta)$ blue particles in the BBM $\C{B}^{N,\nu}$ at time delta. Define $X^+(0)=\tilde{\nu}$. Now consider the smallest integer $\hat{m}$ such that $|\C{N}^{\hat{m},\tilde{\nu}}(\delta)|\geq N+b^N(\delta)$. Since, for $\hat{m}\leq N$, the number of blue particles in $\C{B}^{\hat{m},\tilde{\nu}}(\delta)$ is at most $b^N(\delta)$, thus $|\C{N}^{\hat{m},\tilde{\nu}}(\delta)|\geq N+b^N(\delta)$ implies that $\C{B}^{\hat{m},\tilde{\nu}}(\delta)$ has at least $N$ yellow particles, and therefore (\ref{moreThanNIneq}) yields that $X(\delta)\oop \Gamma^+(\C{B}^{\hat{m},\tilde{\nu}},P,\delta)$.

Then for $t\in (0,\delta)$, define $X^+(t):=\Gamma^+(\C{B}^{\hat{m},\tilde{\nu}},P,t)$. Finally, define $X^+(\delta)$ to be the $N$ rightmost particles (of any colour) of $\Gamma^+(\C{B}^{\hat{m},\tilde{\nu}},P,\delta-)$. Therefore $RY(\C{B}^{\hat{m},\tilde{\nu}}(\delta),N)\oop X^+(\delta)$, and thus $X(\delta)\oop X^+(\delta)$. Now suppose for induction that we have defined the processes $X$ and $X^+$ up to time $k\delta$ such that $X(i\delta)\oop X^+(i\delta)$ for $i\leq k$. Then using the above construction, with initial conditions $X(k\delta)$ and $X^+(k\delta)$, we can thus construct the upper bounding process for times $[k\delta,(k+1)\delta]$. Hence by induction, we can construct $X^+$ for all $t\geq 0$ such that $X(i\delta)\oop X^+(i\delta)$ for all $i\in \B{N}$.

\paragraph{Construction of $X^-$:} Consider now the process $-X$. Observe that the process $-X$ is an $(N,1-p)$-BBM process. Using the construction above, we can construct an upper bound $X^{N,1-p,\delta,+}$ for the $(N,1-p)$-BBM process $-X$ so that $-X(k\delta)\oop X^{N,1-p,\delta,+}(k\delta)$ for $k\in \B{N}$. Therefore as $X^{N,1-p,\delta,+}$ and $-X$ are vectors of the same size (namely $N$) at each time $k\delta$, we have that $-X^{N,1-p,\delta,+}(k\delta)\oop X(k\delta)$ for all $k$. Therefore we define the lower bound $X^-=X^{N,p,\delta,-}$ to be $X^-(t):=-X^{N,1-p,\delta,+}(t)$ for $t\geq 0$. In words, at the start of each interval $(k\delta,k\delta+\delta)$, we delete from the process $X^-$ a number of its leftmost particles. Then over $(k\delta,k\delta+\delta)$, $X^-$ grows and diffuses so that at time $k\delta+\delta$ it has more than $N$ yellow particles, and finally at time $k\delta+\delta$ we delete all but the leftmost $N$ particles of $X^-$. 

Thus we have constructed upper and lower bounding processes for $X$ such that at all discrete times $k\delta$, $k\in \B{N}$, we have $X^-(k\delta)\oop X(k\delta)\oop X^+(k\delta)$. Let us state this as the following Proposition:

\vspace{1em}

\begin{propn} \label{upDownCoupling}
    Let $X$ be an $(N,p)$-BBM process with initial configuration $\nu$, let $\nu^\pm$ be initial configurations such that $\nu^- \oop \nu \oop \nu^+$, and fix $\delta>0$. Then there exist upper and lower bounding processes $X^\pm$ with initial configurations $\nu^\pm$ respectively such that $X^-(k\delta)\oop X(k\delta)\oop X^+(k\delta)$ for all $k\in \B{N}$, and on each interval $((k-1)\delta,k\delta)$, the processes $X^-$ and $X^+$ are branching Brownian motion processes.
\end{propn}

\section{Convergence of the coupled upper and lower bounds} \label{convOfCoupling}
The reason that we defined the upper and lower bounds, $X^\pm$, as we did is because it is easier to describe the dynamics of selection/deletion of particles at deterministic discrete times rather than random times, and therefore stating the limits of the systems $X^\pm = X^{N,p,\delta,\pm}$ is considerably easier. We will describe the limiting behaviour of $X^\pm$  in terms of the Gaussian diffusion operator 
$$G_t f(x):=\int_{\B{R}} (2\pi t)^{-1/2} e^{-(x-y)^2/2t}f(y)dy,$$ and the cut operators 
$$C_k^L f(x):= f(x) \is_{\{\int_x^\infty f(y)dy < k\}} \text{ and } C_k^R f(x):= f(x) \is_{\{\int_{-\infty}^x f(y)dy < k\}}.$$ So $C_k^L$ is the operator which cuts mass on the left and leaves a total mass $k$, and $C_k^R$ is the operator which cuts mass on the left and leaves a total mass $k$. Then we claim that the limiting distributions of $X^-$ and $X^+$ are:
$$S^{\delta,-}_{k\delta}\rho = (C^R_1 e^\delta G_\delta C^L_{1-p(1-e^{-\delta})})^k \rho \quad \text{and}\quad S^{\delta,+}_{k\delta}\rho = (C^L_1 e^\delta G_\delta C^R_{1-(1-p)(1-e^{-\delta})})^k \rho,$$
respectively, in the following sense:

\vspace{1em}

\begin{propn} \label{couplingConv} Fix $\delta>0$ and let $\pi^{N,\delta,\pm}_{\rho,t}$ be the empirical probability measure of the particle system $X^{N,p,\delta,\pm}$ at time $t$ with initial density $\rho$ in the sense that each particle's initial position independently has density $\rho$. Then there exists a constant $C$ such that for any $B,\beta>0$, $k\in \B{N}$, for $N$ sufficiently large:
\begin{align}\label{convProbBound}
\sup_{a\in \B{R}}\B{P}\left(\left|\int_a^\infty \pi^{N,\delta,\pm}_{\rho,k\delta}(dr)-\int_a^\infty S_{k\delta}^{\delta,\pm}\rho(r)dr\right|>BN^{-\beta}\right)\leq B^{-4}Ce^{4k\delta}k^5 e^{4\delta}N^{4\beta-2}\end{align}
\end{propn}

To prove this probabilistic bound, we will first state a number of results which will be useful in the course of the proof, and give various scenarios in which a probabilistic bound of the above form appears. These may be considered `moderate deviations' results. They are simple technical results, so their proofs are omitted here and can be found in the appendix.  

\vspace{1em}

\begin{lemma} \label{llnProbBound}
    Suppose that $A_1,A_2,A_3,\ldots$ is a sequence of random variables with means $\mu_1,\mu_2,\mu_3,\ldots$ such that the sequences $N^{-1}\sum_{i=1}^N \B{E}[(A_i-\mu_i)^4]$ and $N^{-2}\sum_{i\neq j}^N\B{E}[(A_i-\mu_i)^2]\B{E}[(A_j-\mu_j)^2]$ converge, say with $N^{-2}\sum_{i\neq j}^N\B{E}[(A_i-\mu_i)^2]\B{E}[(A_j-\mu_j)^2]\leq k$ for $N$ sufficiently large. Then there exists a constant $C$ independent of the distributions of $A_i$ such that for any $B,\beta>0$, for sufficiently large $N$:
    $$\B{P}\left(\left|\frac1N \sum_{i=1}^N A_i - \frac1N \sum_{i=1}^N \mu_i \right|>BN^{-\beta}\right)\leq B^{-4}CkN^{4\beta -2}.$$
\end{lemma}

\vspace{1em}

\begin{corollary} \label{probBoundAddition}
    Suppose that $A_1,A_2,A_3,\ldots$ and $A'_1,A'_2,A'_3,\ldots$ are sequences of random variables such that there exist constants $C_A$ and $C_{A'}$ such that for any $B,\beta>0$, for $N$ sufficiently large, $\B{P}(|A_N-a|>BN^{-\beta})<B^{-4}C_AN^{4\beta-2}$ and $\B{P}(|A'_N-a'|>BN^{-\beta})<B^{-4}C'_AN^{4\beta-2}$. Then for any $B,\beta>0$, for $N$ sufficiently large,
    $$\B{P}(|(A_N+A'_N)-(a+a')|>BN^{-\beta})<16B^{-4}(C_A+C'_A)N^{4\beta-2}$$
\end{corollary}

\vspace{1em}

\begin{lemma} \label{compoundPoissonLemma}
    For $i=1,2,\ldots,n$, let $A_{i,1},A_{i,2},A_{i,3},\ldots$ be a sequence of i.i.d. random variables with with mean $\mu_i$ such that $\B{E}[A_{i,1}^4]<\infty$. Let $M^1_N,M^2_N,\ldots$ be i.i.d. Poisson random variables with parameter $\lambda N$. Then there exists a constant $C$ independent of $\lambda$, $n$, and $A_{i,j}$ such that for any $B,\beta>0$, for $N$ sufficiently large:
    $$\B{P}\left(\left|\sum_{i=1}^n \frac1N \sum_{j=1}^{M^i_N} A_{i,j} - \lambda \sum_{i=1}^n \mu_i\right|>BN^{-\beta}\right)<B^{-4}C\lambda^2 \left(\sum_{i=1}^n\B{E}[A_{i,1}^2]^2+6\sum_{i,j=1,i\neq j}^n \B{E}[A_{i,1}^2]\B{E}[A_{j,1}^2]\right)N^{4\beta-2}.$$
\end{lemma}

\vspace{1em}

\begin{propn} \label{preKillProp}
    Let $|X^{N,p,\delta,\pm}(t)|$ denote the number of particles in the process $X^{N,p,\delta,\pm}$ at time $t$. Then there exists a constant $C$ such that for any $B,\beta,\delta>0$, for $N$ sufficiently large
    \begin{align}\label{lemmaMinusAmounts}\B{P}\left(\left||X^{N,p,\delta,-}(\delta-)|/N-(e^\delta(1-p)+p)\right|>BN^{-\beta}\right)< B^{-4}C\delta^2 e^{4\delta}N^{4\beta-2},\end{align}
    and
    \begin{align}\label{lemmaPlusAmounts}\B{P}\left(\left||X^{N,p,\delta,+}(\delta-)|/N-(pe^\delta+(1-p))\right|>BN^{-\beta}\right)< B^{-4}C\delta^2 e^{4\delta}N^{4\beta-2}.\end{align}
\end{propn}

\begin{proof}
Consider the process $X^+=X^{N,p,\delta,+}$ defined by the construction in Section \ref{twoSidedCoupling} as a function of $\C{B}^{N,\nu}$ and $P$. Recall that $X^{+}(0+)$ is constructed by taking the minimum number $\hat{m}$ of leftmost particles of the BBM $\C{B}^{N,\nu}$ at time $0$ so that there are at least $N+b^N(\delta)$ particles in $\C{B}^{\hat{m},\nu}$ alive at time $\delta$, where $b^N(\delta)$ is the number of blue particles in $\C{B}^{N,\nu}(\delta)$. 

Now let $K_i(\delta)$ denote the number of offspring (of any colour) at time $\delta$ of the $i$\textsuperscript{th} leftmost particle of $\C{B}^N(0)$ for $i=1,2,\ldots,N$. So for each $i=1,2,\ldots,N$, $K_i(\delta)$ are i.i.d. random variables with geometric distribution of parameter $e^\delta$ (since $K_i$ counts the number of individuals in a Yule processes). Thus $|X^+(\delta-)|$ differs from $N+b^N(\delta)$ by at most $\max_{i=1}^N K_i(\delta)$. Then since $K_i(\delta)$ has finite moments of all orders, Corollary 7 of \cite{downey} states that $\B{E}\left[\max_{i=1}^N K_i(\delta)\right]=o(N^\gamma)$ for any $\gamma>0$. Therefore by Markov's inequality:
$$\B{P}\left(N^{-1}\max_{i=1}^N K_i(\delta)>BN^{-\beta}\right)=o\left(N^{\gamma}/N^{2-2\beta}\right)$$
for any $\gamma>0$. So choosing $\gamma$ such that $\gamma < 2\beta$, we can use the union bound to yield:
\begin{align*}
    \B{P}(||X^+(\delta-)|/N- (pe^\delta + 1-p)|&>BN^{-\beta})\\
    \leq \B{P}(&1+b^N(\delta) + \max_{i=1}^N K_i(\delta)/N - (pe^\delta + 1-p)>BN^{-\beta}) \\
    &+ \B{P}(1+b^N(\delta)/N - (pe^\delta + 1-p)<-BN^{-\beta}) \\
    \leq \B{P}(&\max_{i=1}^N K_i(\delta) > BN^{1-\beta}/2) + \B{P}(|1+b^N(\delta)/N - (pe^\delta + 1-p)|>BN^{-\beta}/2) \\
    \leq o(&N^{2\beta+\gamma - 2}) + \B{P}(|1+b^N(\delta)/N - (pe^\delta + 1-p)|>BN^{-\beta}/2) 
\end{align*}

Now let us consider how to count $b^N(\delta)$, the number of blue particles in $\C{B}^{N,\nu}$. By construction, $b^N(t)$ counts the number of individuals in a population with initially zero individuals, with immigration at rate $Np$ and branching at rate $1$ independently for each individual. Immigration at rate $1$ occurs because yellow particles are turned blue at rate $Np$, and then subsequently each blue particle branches independently at rate $1$. Now fix $n\in \B{N}$ and for $i\leq n$, let $M^{i,n}$ be the number of particles which turn from yellow to blue in $\C{B}^{N,\nu}$ during $\left[ \frac{(i-1)\delta}{n},\frac{i\delta}{n}\right)$, which is a Poisson random variable of parameter $Np\delta/m$. By the strong law of large numbers for Poisson processes, $M^{i,n}/N\to p\delta/n$ almost surely as $N\to\infty$. 

Now, for $i=1,2,\ldots,n$ let $A^i$ count the number of offspring, after $\delta-\frac{i\delta}{n}$ units of time, in a Yule process with branching rate $1$ which starts from a single particle. Let $A_1^i,A_2^i,\ldots$ be i.i.d copies of $A^i$. So $\B{E}A_j^i = e^{\delta-\frac{i\delta}{n}}$. Then the number of descendants at time $\delta$ of a particle which immigrates into the population during $[\frac{(i-1)\delta}{n}, \frac{i\delta}{n})$ is stochastically bounded above and below by $A_j^{i-1}$ and $A_j^i$ respectively. So by the strong law of large numbers:
\begin{align} \label{noParticlesUpper}
U_{N,n} :=\frac1N \sum_{i=1}^n \sum_{j=1}^{M^{i,n}}A_j^{i-1} = \sum_{i=1}^n \frac{M^{i,n}}{N}\frac{1}{M^{i,n}}\sum_{j=1}^{M^{i,n}}A_j^{i-1} \xrightarrow[N\to\infty]{a.s.,L^1}\sum_{i=1}^n \frac{p\delta}{n}e^{\delta-\frac{(i-1)\delta}{n}}=:u(n),\end{align}
and
\begin{align} \label{noParticlesLower}
L_{N,n}:=\frac1N \sum_{i=1}^n \sum_{j=1}^{M^{i,n}}A_j^{i} = \sum_{i=1}^n \frac{M^{i,n}}{N}\frac{1}{M^{i,n}}\sum_{j=1}^{M^{i,n}}A_j^{i} \xrightarrow[N\to\infty]{a.s.,L^1}\sum_{i=1}^n \frac{p\delta}{n}e^{\delta-\frac{i\delta}{n}}=:\ell(n).\end{align}
So $U_{N,n}$ and $L_{N,n}$ provide stochastic upper and lower bounds for $b^N(t)/N$. Furthermore, they converge to $u(n)$ and $\ell(n)$ respectively, which are upper and lower bounding Reimann sums approximating $\int_0^t pe^x dx$, and converging to $\int_0^t pe^x dx = p(e^t-1)$ as $n$ goes to infinity. 

Observe that for any $i,n$, $M^{i,n}$ is a Poisson random variable with parameter $p\delta N/n$ and that  for any $i,j$, $A^i_j$ has bounded fourth and second moments. More precisely, $A_j^i$ has second moment $e^{2(\delta-i\delta/n)}\leq e^{2\delta}$ and therefore by Lemma
\ref{compoundPoissonLemma}, there exists constants $C,C'$ such that for any $B,\beta>0$ and any $n\in \B{N}$:
\begin{align} \label{upperLowerllnBnd}
    \B{P}(|U_{N,n}-u(n)|>BN^{-\beta})\leq B^{-4}C'\left(\sum_{i=1}^n \big(\frac{p\delta}{n}\big)^2 e^{4\delta} + 6\sum_{i,j=1,i\neq j}^n \big(\frac{p\delta}{n}\big)^2 e^{4\delta}\right)N^{4\beta-2}
    \leq B^{-4}Ce^{4\delta}\delta^2 N^{4\beta-2},
\end{align}
and similarly this bound holds for $\B{P}(|L_{N,n}-\ell(n)|>BN^{-\beta})$. By the union bound
$$\B{P}(|1+b^N(\delta)/N - (1+pe^\delta-p)|>BN^{-\beta})\leq \B{P}(b^N(\delta)/N - (pe^\delta-p)>BN^{-\beta}) + \B{P}(b^N(\delta)/N - (pe^\delta-p)<-BN^{-\beta}).$$
By the triangle inequality, we have
\begin{align*}\B{P}(b^N(\delta)/N - (pe^\delta-p)>&BN^{-\beta}) \leq \B{P}(U_{N,n} - (pe^\delta-p)>BN^{-\beta}) \\
&\leq \B{P}(U_{N,n}-u(n)>BN^{-\beta}/2) + \B{P}(u(n)-(pe^{\delta}-p)>BN^{-\beta}/2) \\
&\leq 16B^{-4}CN^{4\beta-2}\delta^2 e^{4\delta}+\is_{\{BN^{-\beta}/2<u(n)-(1+pe^\delta-p)\}},\end{align*}
therefore for all $N$, since $u(n)\xrightarrow[n\to\infty]{}1+pe^\delta - p$, we can choose $n$ sufficiently large to prove that $$\B{P}(|b^N(\delta)/N - (pe^\delta - p)|>BN^{-\beta}) \leq B^{-4}CN^{4\beta-2}\delta^2 e^{4\delta}.$$
Similarly, by using the lower bounds $L_{N,n}$ and $\ell(n)$, we have that $\B{P}(b^N(\delta)/N - (pe^\delta-p)<-BN^{-\beta})\leq B^{-4}Ce^{4\delta}\delta^2N^{4\beta-2}$. Putting these two bounds together gives the result.
\end{proof}

In order to prove Proposition \ref{couplingConv}, we first prove the following Lemma, from which the result will follow by an inductive argument. 

\vspace{1em}

\begin{lemma}
    \label{couplingConvLemma}
    Fix $\delta > 0$, and let $(\rho^N)_{N\geq 1}$ be a sequence of probability measures describing initial configurations of $N$ particles such that no two particles are in the same location (i.e. $\rho^N$ is a sum of $N$ atoms of weight $1/N$ at distinct locations). Let $\pi^{N,\delta,\pm}_{\rho^N,t}$ be the empirical probability measure of the particle system $X^{N,p,\delta,\pm}$ at time $t$ with the particles initial configurations described by $\rho^N$. Then there exists a constant $C$ independent of $(\rho^N)_{N\geq 0}$ such that for any $B,\beta,\delta>0$, we have that for $N$ sufficiently large:
    \begin{align}
        \sup_{a\in \B{R}}\B{P}\Bigg(\Big|\int_a^\infty \pi^{N,\delta,\pm}_{\rho^N,\delta}(dr)-\int_a^\infty S_\delta^{\delta,\pm}\rho^N(dr)\Big|>BN^{-\beta}\Bigg) \leq B^{-4}Ce^{4\delta}N^{4\beta-2}.
    \end{align}
\end{lemma}

Before proving the above Proposition, let us remark that the locations of the $N$ atoms in the distribution $\rho^N$ may be deterministic or drawn as i.i.d. samples from a probability distribution on $\B{R}$, provided that with probability $1$, no two particles are in the same location.

\begin{proof}
    We will focus here on proving the case of $X^{N,p,\delta,-}$, since the case $X^{N,p,\delta,+}$ then follows immediately by symmetry. Let us begin by defining the cutpoints $\tilde{L}^{N,\delta}$ and $\tilde{R}^{N,\delta}$. These will be the point at which the cutting operators $C^R_1$ and $C^L_{1-p(1-e^{-\delta})}$ which make up the operator $S_\delta^{\delta,-}$ cut, when applied to the initial configuration $\rho^N$. So:
    $$\tilde{L}^{N,\delta}:=\sup\Big\{a:\int_a^\infty\rho^N(dr) \geq 1-p(1-e^{-\delta})\Big\}\text{ and }\tilde{R}^{N,\delta}:=\inf\Big\{a:\int_{-\infty}^a e^\delta G_\delta C^L_{1-p(1-e^{-\delta})} \rho^N(dr)\geq 1\Big\},$$
    so in particular, $\tilde{L}^{N,\delta}$ is the point left of which we delete mass from the initial configuration $\rho^N$ at the start of the interval $[0,\delta]$ due to $C^L_{1-p(1-e^{-\delta})}$, and $\tilde{R}^{N,\delta}$ is the point right of which we delete mass from the function $e^\delta G_\delta C^L_{1-p(1-e^{-\delta})}\rho^N$ at the end of the time interval $[0,\delta]$ due to $C^R_1$. Let $x^N_1,\ldots,x^N_N$ be the locations of the $N$ atoms of $\rho^N$, and consider $N$ independent branching Brownian motions (BBMs) starting from the locations $x^N_1,\ldots,x^N_N$. Say that the $i$\textsuperscript{th} BBM has, at time $t$, $N^i_t$ particles, which we label $B^{i,1},\ldots,B^{i,N^i_t}$. Then define the empirical distribution $\tilde{\pi}^{N,\delta}_{\rho^N,\delta}$ by:
    $$\tilde{\pi}^{N,\delta}_{\rho^N,\delta}(A):=\frac1N \sum_{i=1}^N \sum_{j=1}^{N^i_\delta}\is_{\{B^{i,j}_\delta\in A\}}\is_{\{B_0^{i,j}\geq \tilde{L}^{N,\delta}\}}\is_{\{B_\delta^{i,j}\leq \tilde{R}^{N,\delta}\}}.$$

    So this is the normalised empirical distribution of particles in a BBM with deletion of particles left of $\tilde{L}^{N,\delta}$ at time $0$ and particles right of $\tilde{R}^{N,\delta}$ at time $\delta$. Essentially, we are using $\tilde{\pi}$ to approximate $\pi$; using deterministic cutoff points $\tilde{L}^{N,\delta}$ and $\tilde{R}^{N,\delta}$ to approximate the random cutoff points in the process $X^+$.
    
    Now consider the random variables \sloppy$\sum_{j=1}^{N_\delta^i}\is_{\{B^{i,j}_\delta\in A\}}\is_{\{B_0^{i,j}\geq \tilde{L}^{N,\delta}\}}\is_{\{B_\delta^{i,j}\leq \tilde{R}^{N,\delta}\}}$, which have finite moments of all orders and second moment on the order of $e^{2\delta}$. Then notice that $\B{E}\Big[\int_a^\infty \tilde{\pi}^{N,\delta,-}_{\rho^N,\delta}(dr)\Big] = \int_a^\infty S_\delta^{\delta,-}\rho^N(dr)$. Therefore Lemma (\ref{llnProbBound}) gives that there exists $C_1$ such that for $B,\beta>0$ and $a\in \B{R}$, then for $N$ sufficiently large:
    \begin{align} \label{piTildeSLink}
        \B{P}\Bigg(\Big|\int_a^\infty \tilde{\pi}^{N,\delta}_{\rho^N,\delta}(dr)-\int_a^\infty S_\delta^{\delta,-}\rho^N(dr)\Big|>BN^{-\beta}\Bigg) \leq B^{-4}C_1e^{4\delta}N^{4\beta-2}.
    \end{align}
    
    In the case that $\rho^N$ is random and describe an initial configuration in which $N$ locations are i.i.d. samples of a distribution $\rho$, then this can be strengthened to:
    \begin{align*}
        \B{P}\Bigg(\Big|\int_a^\infty \tilde{\pi}^{N,\delta}_{\rho^N,\delta}(dr)-\int_a^\infty S_\delta^{\delta,-}\rho(dr)\Big|>BN^{-\beta}\Bigg) \leq B^{-4}C_1e^{4\delta}N^{4\beta-2}.
    \end{align*}

    We can similarly express the empirical density of $X^{N,p,\delta,-}$ with initial configuration $\rho^N$, $\pi^{N,\delta,-}_{\rho^N,\delta}$, by considering the branching Brownian motions and cutting points. So let $L^{N,\delta}$ and $R^{N,\delta}$ be the cutpoints such that for any Borel measurable set $A$,
    \begin{align*}
        \pi^{N,\delta}_{\rho^N,\delta}(A):=\frac1N \sum_{i=1}^N \sum_{j=1}^{N_\delta^i} \is_{\{B^{i,j}_\delta \in A\}}\is_{\{B_0^{i,j}\geq L^{N,\delta}\}}\is_{\{B^{i,j}_\delta\leq R^{N,\delta}\}}.
    \end{align*}

    Therefore by the triangle inequality, the inequality (\ref{piTildeSLink}), and Corollary \ref{probBoundAddition}, to prove the desired result, it remains to prove a bound for 
    \begin{align*}
    \B{P}\left(\left|\int_a^\infty \pi_{\rho^N,\delta}^{N,\delta,-}(dr)-\int_a^\infty \tilde{\pi}_{\rho^N,\delta}^{N,\delta,-}(dr)\right| \geq BN^{-\beta}\right).\end{align*}
    
    Then using the triangle inequality:
    \begin{align}
    \nonumber \Big|\int_a^\infty \pi_{\rho^N,\delta}^{N,\delta,-}(dr) - &\int_a^\infty \tilde{\pi}^{N,\delta,-}_{\rho^N,\delta}(dr)\Big| \leq \frac1N \sum_{i=1}^N \sum_{j=1}^{N_\delta^i}\Big| \is_{\{B_0^{i,j}\geq L^{N,\delta}\}}\is_{\{B_\delta^{i,j}\leq R^{N,\delta}\}}-\is_{\{B_0^{i,j}\geq \tilde{L}^{N,\delta}\}}\is_{\{B_\delta^{i,j}\leq \tilde{R}^{N,\delta}\}}\Big| \\
    \label{absolIn} \leq \frac1N &\sum_{i=1}^N \sum_{j=1}^{N_\delta^i}\Big| \is_{\{B_0^{i,j}\geq L^{N,\delta}\}}\is_{\{B_\delta^{i,j}\leq R^{N,\delta}\}}-\is_{\{B_0^{i,j}\geq L^{N,\delta}\}}\is_{\{B_\delta^{i,j}\leq \tilde{R}^{N,\delta}\}}\Big| \\
    \nonumber &+ \frac1N \sum_{i=1}^N \sum_{j=1}^{N_\delta^i}\Big| \is_{\{B_0^{i,j}\geq L^{N,\delta}\}}\is_{\{B_\delta^{i,j}\leq \tilde{R}^{N,\delta}\}}-\is_{\{B_0^{i,j}\geq \tilde{L}^{N,\delta}\}}\is_{\{B_\delta^{i,j}\leq \tilde{R}^{N,\delta}\}}\Big| \\
    \label{absolOut} \leq \Big| \frac1N &\sum_{i=1}^N \sum_{j=1}^{N_\delta^i} \is_{\{B_0^{i,j}\geq L^{N,\delta}\}}\is_{\{B_\delta^{i,j}\leq R^{N,\delta}\}}-\frac1N \sum_{i=1}^N \sum_{j=1}^{N_\delta^i} \is_{\{B_0^{i,j}\geq L^{N,\delta}\}}\is_{\{B_\delta^{i,j}\leq \tilde{R}^{N,\delta}\}}\Big| \\
    &+ \frac1N \sum_{i=1}^N \sum_{j=1}^{N_\delta^i} \big| \is_{\{B_0^{i,j}\geq L^{N,\delta}\}}- \is_{\{B_0^{i,j}\geq \tilde{L}^{N,\delta}\}}\big| \\
    \label{couplingConvMiddleBound} \leq \Big| \frac1N &\sum_{i=1}^N \sum_{j=1}^{N_\delta^i} \is_{\{B_0^{i,j}\geq L^{N,\delta}\}}\is_{\{B_\delta^{i,j}\leq R^{N,\delta}\}}-\frac1N \sum_{i=1}^N \sum_{j=1}^{N_\delta^i} \is_{\{B_0^{i,j}\geq \tilde{L}^{N,\delta}\}}\is_{\{B_\delta^{i,j}\leq \tilde{R}^{N,\delta}\}}\Big| \\
    \nonumber + \Big| &\frac1N \sum_{i=1}^N \sum_{j=1}^{N_\delta^i} \is_{\{B_0^{i,j}\geq \tilde{L}^{N,\delta}\}}\is_{\{B_\delta^{i,j}\leq\tilde{R}^{N,\delta}\}}-\frac1N \sum_{i=1}^N \sum_{j=1}^{N_\delta^i} \is_{\{B_0^{i,j}\geq L^{N,\delta}\}}\is_{\{B_\delta^{i,j}\leq \tilde{R}^{N,\delta}\}}\Big| \\
    &+ \frac1N \sum_{i=1}^N \sum_{j=1}^{N_\delta^i} \big| \is_{\{B_0^{i,j}\geq L^{N,\delta}\}}- \is_{\{B_0^{i,j}\geq >\tilde{L}^{N,\delta}\}}\big| \\
    \label{couplingConvFinalBound}
    \leq \Big| \frac1N &\sum_{i=1}^N \sum_{j=1}^{N_\delta^i} \is_{\{B_0^{i,j}\geq L^{N,\delta}\}}\is_{\{B_\delta^{i,j}\leq R^{N,\delta}\}}-\frac1N \sum_{i=1}^N \sum_{j=1}^{N_\delta^i} \is_{\{B_0^{i,j}\geq \tilde{L}^{N,\delta}\}}\is_{\{B_\delta^{i,j}\leq \tilde{R}^{N,\delta}\}}\Big|\\
    & \label{couplingConvFinalBound2} + \frac2N \sum_{i=1}^N \sum_{j=1}^{N_\delta^i} \big| \is_{\{B_0^{i,j}\geq L^{N,\delta}\}}- \is_{\{B_0^{i,j} \geq \tilde{L}^{N,\delta}\}}\big|   \\
    \label{couplingConvFinalFinalBound}
    \leq \Big| \frac1N &\sum_{i=1}^N \sum_{j=1}^{N_\delta^i} \is_{\{B_0^{i,j} \geq L^{N,\delta}\}}\is_{\{B_\delta^{i,j}\leq R^{N,\delta}\}}-\frac1N \sum_{i=1}^N \sum_{j=1}^{N_\delta^i} \is_{\{B_0^{i,j}\geq \tilde{L}^{N,\delta}\}}\is_{\{B_\delta^{i,j}\leq \tilde{R}^{N,\delta}\}}\Big| \\
    \nonumber
    & + \Bigg| \frac2N \sum_{i=1}^N \sum_{j=1}^{N_\delta^i} \is_{\{B_0^{i,j}\geq L^{N,\delta}\}}- (e^\delta - pe^\delta + p) \Bigg| +  \Bigg| \frac2N \sum_{i=1}^N \sum_{j=1}^{N_\delta^i} \is_{\{B_0^{i,j}\geq \tilde{L}^{N,\delta}\}}- (e^\delta - pe^\delta + p) \Bigg|
\end{align}
where the inequalities (\ref{absolIn}),(\ref{couplingConvMiddleBound}), (\ref{couplingConvFinalBound}) follow from the triangle inequality, the inequality (\ref{absolOut}) follows from the fact that every term in the first sum of (\ref{absolIn}) has the same sign (namely the sign of $R^{N,\delta}-\tilde{R}^{N,\delta}$), and the inequality (\ref{couplingConvFinalFinalBound}) follows from the fact that every term in (\ref{couplingConvFinalBound2}) has the same sign (namely the sign of $\tilde{L}^{N,\delta}-L^{N,\delta}$) and the triangle inequality. 

Let us first deal with the first term of (\ref{couplingConvFinalFinalBound}). By definition of $L^{N,\delta}$ and $R^{N,\delta}$, we certainly have that by $\frac1N \sum_{i=1}^N \sum_{j=1}^{N_\delta^i}\is_{\{B_0^{i,j}\geq L^{N,\delta}\}}\is_{\{B_\delta^{i,j}\leq R^{N,\delta}\}}= 1$. Next, observe that since $e^\delta G_\delta C^L_{1-p(1-e^{-\delta})}\rho^N(x)$ is a smooth function, thus $\int_{-\infty}^{\tilde{R}^{N,\delta}}e^\delta G_\delta C^L_{1-p(1-e^{-\delta})}\rho^N(dr)$ (and therefore $\int_{\B{R}}S^{\delta,-}_{\delta}\rho^N(dr)$) is exactly equal to $1$. Therefore applying equation (\ref{piTildeSLink}) with $a=-\infty$, we get that for $N$ sufficiently large:
\begin{align}\label{firstTermBound} \B{P}\Bigg(\Bigg|\frac1N \sum_{i=1}^N \sum_{j=1}^{N_\delta^i}\is_{\{B_0^{i,j}\geq \tilde{L}^{N,\delta}\}}\is_{\{B_\delta^{i,j}\leq \tilde{R}^{N,\delta}\}} - 1 \Bigg| >BN^{-\beta}\Bigg)< B^{-4}C_2e^{4\delta}N^{4\beta-2},
\end{align}
where the factor of $e^{4\delta}$ comes from the fact that the second moment of the random variable $\sum_{j=1}^{N^i_{\delta}}\is_{\{B_0^{i,j}\geq \tilde{L}^{N,\delta}\}}\is_{\{B_\delta^j \leq \tilde{R}^{N,\delta}\}}$ is bounded above by $\B{E}[(N_\delta^i)^2]\leq 2e^{2\delta}$. 

We now turn our attention to the second and third terms of (\ref{couplingConvFinalFinalBound}). Now by Proposition \ref{preKillProp}, since $\sum_{i=1}^N \sum_{j=1}^{N_\delta^i}\is_{\{B_0^{i,j}>L^{N,\delta}\}}=|X^{N,p,\delta,-}(\delta-)|$, thus there exists $C_3$ such that for any $B,\beta,\delta>0$ and $N$ sufficiently large, 
\begin{align}
    \label{firstPartOfSecondTermBound} \B{P}\left(\left|\frac1N \sum_{i=1}^N \sum_{j=1}^{N_\delta^i} \is_{\{B_0^{i,j}>L^{N,\delta}\}} - (e^\delta(1-p)+p)\right|>BN^{-\beta}\right) < B^{-4}C_3\delta^2 e^{4\delta}N^{4\beta-2},
\end{align}
Now consider the cumulative distribution of the initial configuration $\rho^N$. Since $\rho^N$ is the sum of $N$ distinct atoms, each of mass at most $1/N$, thus $\sum_{i=1}^N \is_{\{B^{i,j}_\delta \geq \tilde{L}^{N,\delta}\}}-(1-p(1-e^{-\delta})) \leq 1/N$. Therefore $\B{E}\big[\sum_{i=1}^N\sum_{j=1}^{N_\delta^i}\is_{\{B_0^{i,j}>L^\delta\}}\big]-e^\delta(1-p(1-e^{-\delta}))\leq e^\delta/N$, so by Lemma \ref{llnProbBound}, there exists $C_4$ such that for any $B\in (0,1),\beta>0$, then for $N$ sufficiently large we have
by Lemma \ref{llnProbBound} and the union bound, there exists $C_4$ such that for any $B,\beta,\delta>0$, for $N$ sufficiently large
\begin{align}
    \label{secondPartOfSecondTermBound} \B{P}\Bigg(\Bigg|\frac1N \sum_{i=1}^N \sum_{j=1}^{N_\delta^i} \is_{\{B_0^{i,j}>\tilde{L}^{N,\delta}\}} - &e^\delta (1-p(1-e^{-\delta}))\Bigg|>BN^{-\beta}\Bigg) \\
    \nonumber\leq \B{P}\Bigg(\Bigg|\frac1N \sum_{i=1}^N&\sum_{j=1}^{N_\delta^i} \is_{\{B_0^{i,j}>\tilde{L}^{N,\delta}\}}-\B{E}\Big[\frac1N \sum_{i=1}^N\sum_{j=1}^{N_\delta^i} \is_{\{B_0^{i,j}>\tilde{L}^{N,\delta}\}}\Big]\Bigg|>BN^{-\beta}/2\Bigg)\\
    \nonumber &+\B{P}(e^\delta/N > BN^{-\beta}/2)\\
    \nonumber \leq B^{-4}C_4 e^{4\delta}&N^{4\beta-2},
\end{align}
where again, the factor of $e^{4\delta}$ comes from the fact that the second moment of the random variable $\sum_{j=1}^{N_\delta^i} \is_{\{B_0^{i,j}\geq \tilde{L}^{N,\delta}\}}\is_{\{B_\delta^i\leq \tilde{R}^{N,\delta}\}}$ is bounded above by $\B{E}[(N_\delta^i)^2]\leq 2e^{2\delta}$. Therefore, by Corollary \ref{probBoundAddition}, putting together (\ref{firstTermBound}), (\ref{firstPartOfSecondTermBound}), and (\ref{secondPartOfSecondTermBound}), it follows that there exists a constant $C_5$ such that for $N$ sufficiently large, 
\begin{align}
    \B{P}\left(\left|\int_a^\infty \pi_{\rho^N,\delta}^{N,\delta,-}(dr)-\int_a^\infty \tilde{\pi}_{\rho^N,\delta}^{N,\delta,-}(dr)\right|\geq BN^{-\beta}\right)\leq B^{-4}C_5 e^{4\delta}N^{4\beta-2}
\end{align}
Then putting this together with (\ref{piTildeSLink}), we get that there exists a $C$ such that for $B,\beta>0$, for $N$ sufficiently large:
\begin{align*}\label{uniformHydroBound} \sup_{a\in \B{R}}\B{P}\left(\left| \int_a^\infty \pi^{N,\delta,-}_{\rho^N,\delta}(dr) - \int_a^\infty S_\delta^{\delta,-}\rho^N(r)dr\right|>BN^{-\beta}\right)<B^{-4}Ce^{4\delta}N^{4\beta-2}.\end{align*} \end{proof}

The proof of Proposition \ref{couplingConv} then follows by an inductive argument.

\begin{proof} \textit{(Of Proposition \ref{couplingConv})} 
Consider the empirical measure-valued process $(\pi^{N,\delta,-}_{\rho,k\delta})_{k\geq 1}$, and for ease of notation, write $\rho^N_i:=\pi^{N,\delta}_{\rho,i\delta}$, with $\rho^N_0=\rho$. Note that with this definition, $\pi^{N,\delta}_{\rho,i\delta}=\pi^{N,\delta}_{\rho^N_{i-1},\delta}$. Therefore we can write $\int_a^\infty \pi^{N,\delta,-}_{\rho, k\delta}- \int_a^\infty S_{k\delta}^{\delta,-}\rho(dr)$ as a telescoping sum:
\begin{align*} \int_a^\infty \pi^{N,\delta,-}_{\rho^N_{k-1},\delta}(dr) - \int_a^\infty S_\delta^{\delta,-}\rho^N_{k-1}(r)dr+\int_a^\infty S_\delta^{\delta,-}\rho^N_{k-1}(r)dr - \int_a^\infty S_{2\delta}^{\delta,-}\rho^N_{k-2}(r)dr +\cdots \\
- \int_a^\infty S_{(k-1)\delta}^{\delta,-}\rho^N_1(r)dr + \int_a^\infty S_{(k-1)\delta}^{\delta,-}\rho^N_1(r)dr - \int_a^\infty S_{k\delta}^{\delta,-}\rho_0^N(r)dr,\end{align*}

Then we can make use of the fact (see Lemma \ref{operatorProps}, properties (a)\&(c) below) that for any $a\in \B{R}$:
$$\Bigg|\int_a^\infty S_{\delta}^{\delta,-}f(x)dx - \int_a^\infty S_{\delta}^{\delta,-}g(x)dx\Bigg| \leq e^\delta \int_{-\infty}^\infty \big|f(x)dx - g(x)\big|dx.$$ 

Therefore as the result of Lemma \ref{couplingConvLemma} holds, with $C$ independent of $\rho^N$, we can apply it to every sequence of initial conditions $(\rho^N_i)_{N\geq 1}$ for $i=0,1,\ldots,k$. Observe that the first sequence (i.e. $(\rho_0^N)_{N\geq 1}$) is random, and each subsequent sequence $(\rho^N_i)_{N\geq 1}$ is deterministic in the sense that it is measurable with respect to the natural filtration of the process $X^{N,\delta,-}$ up to time $i\delta$. Let us define $$U_\ell^N := \left|\int_{a}^\infty \pi_{\delta}^{N,\delta,-}\rho^N_\ell(dr) - \int_a^\infty S_{\delta}^{\delta,-}\rho^N_\ell(r)dr \right| .$$
Therefore by the triangle inequality:
$$\Bigg| \int_a^\infty \pi_{\rho,k\delta}^{N,\delta,-}(dr) - \int_a^\infty S_{k\delta}^{\delta,-}\rho(r)dr\Bigg| \leq U^N_{k-1}+ e^\delta U^N_{k-2} + \cdots + e^{(k-2)\delta}U^N_1+ e^{(k-1)\delta}U^N_0.$$

So then by the union bound, for any $a\in \B{R}$ and $B,\beta>0$, then for $N$ sufficiently large:
\begin{align*}\sup_{a\in \B{R}}\B{P}\Big(\Big| \int_a^\infty \pi_{\rho,k\delta}^{N,\delta,-}(dr)-\int_a^\infty S_{k\delta}^{\delta,-}\rho(r)dr\Big|>BN^{-\beta}) \leq \sum_{\ell=0}^{k-1} \B{P}(e^{(k-1-\ell)\delta}|U_\ell|>BN^{-\beta}/k) \\
\leq \sum_{\ell=0}^{k-1} \left(\frac{B}{ke^{(k-1-\ell)\delta}}\right)^{-4}Ce^{4\delta}N^{4\beta-2}\leq B^{-4}Ck^5e^{4k\delta}e^{4\delta}N^{4\beta-2},\end{align*}
which completes the proof.
\end{proof}

\section{Convergence of the deterministic limits}
\label{detBarrierSameLimit}
In this section, we will prove that both the upper and lower deterministic bounds, $S_t^{\delta,-}\rho$ and $S_t^{\delta,+}\rho$, converge to the same limit in the appropriate sense as $\delta\downarrow 0$. We employ a similar strategy to that found in Section 5 of \cite{hydroNBBM}, with some added difficulty coming from the fact that $S_t^{\delta,\pm}$ both have cutting of mass on both the left and the right. 

For the remainder of this section, it will be easiest if we define the operators $D_m^L,D_m^R$ to be the operators which cut a mass of $m$ on the left and right hand side respectively. This is in contrast to the operators $C^L_m$ and $C^R_m$ in which $m$ denotes the amount of mass \textit{remaining after cutting} rather than the amount of mass being cut. So $D_m^L f(x)=\is_{x\geq d}f(x)$ where $\int_d^\infty f(x)dx=\int_{-\infty}^\infty f(x)dx - m$. For function $u$ with total mass $\int_{\B{R}}u(y)dy = 1$, it is easily calculated that the operators:
$$S_{k\delta}^{\delta,-}\rho = (C_1^Re^\delta G_\delta C_{1-p(1-e^{-\delta})}^L)^k\rho \quad \text{and} \quad S_{k\delta}^{\delta,+}\rho = (C_1^L e^\delta G_\delta C_{1-(1-p)(1-e^{-\delta})}^R)^k \rho$$
can be rewritten as 
$$S^{\delta,-}_{k\delta}\rho = (D^R_{(1-p)(e^\delta-1)} e^\delta G_\delta D^L_{p(1-e^{-\delta})})^k \rho \quad \text{and} \quad S^{\delta,+}_{k\delta}\rho = (D^L_{p(e^\delta-1)} e^\delta G_\delta D^R_{(1-p)(1-e^{-\delta})})^k \rho.$$

\vspace{1em}

\begin{propn} \label{detBoundDist} For any $\delta > 0$ and $k\in \B{N}$,
    $$\norm{ S_{k\delta}^{\delta,+}\rho(x) - S_{k\delta}^{\delta,-}\rho(x)}_{1} < 2(e^{\delta}-1)e^{k\delta},$$
    where the $\|\cdot\|_1$ is defined by $\|f\|_1 = \int_\B{R}|f(s)|ds$
\end{propn}

Notice that this immediately implies that $\big|\int_x^\infty S_{k\delta}^{\delta,+}\rho(s)ds - \int_x^\infty S_{k\delta}^{\delta,-}\rho(s)ds\big|<2(e^\delta-1)e^{k\delta}$ for all $x\in \B{R}$. Let us start by proving the following properties of the operators $D^L$, $D^R$, and $G$. 

\vspace{1em}

\begin{lemma} \label{operatorProps}
    For positive function $u$ we have the following properties:
    \begin{enumerate}[(a)]
    \setlength{\itemindent}{.5in}
        \item $\|D^{L/R}_a u- D_b^{L/R}u\|_1\leq |a-b|$; 
        \item $\|D^{L/R}_a  - D_a^{L/R}v\|_1 \leq \|u-v\|_1$;
        \item $\|G_t u - G_t v\|_1 \leq \|u-v\|_1$;
        \item $D_a^L D_b^R u = D_b^R D_a^L u$;
        \item $D_a^{L}D_b^{L}u = D^{L}_{a+b}u$ and $D_a^RD_b^Ru = D_{a+b}^Ru$;
        \item For scalar $c>0$, $D^{L/R}_a (cu) = cD^{L/R}_{a/c}u$.
    \end{enumerate}
\end{lemma}

\begin{proof}Properties \textit{(b),(c)} are proven in \cite{hydroNBBM} (properties \textit{(d),(e)} of Proposition 12), and properties \textit{(d),(e)} follow immediately from the definition. We will prove property \textit{(a)} for $\|D^L_a u - D^L_b u\|_1$, and the case for $\|D^R_a u - D^R_b u\|_1$ follows by symmetry. Suppose without loss of generality that $a>b$. Then:
$$\|D^L_a u - D^L_b u \|_1 = \int_\B{R}|\is_{\{x\geq d_a\}}u(x) - \is_{\{x\geq d_b\}}u(x)|dx$$
where $d_a$ and $d_b$ are such that $\int_{d_a}^\infty u(x)dx = \int_\B{R}u(x)dx - a$ and $\int_{d_b}^\infty u(x)dx = \int_\B{R} u(x)dx - b$. Then every term $\is_{\{x\geq d_a\}}u(x) - \is_{\{x\geq d_b\}}u(x)$ has the same sign, so we can take the norm outside the integral so that
$$\|D_a^L u - D_b^L u \|_1 = \Big|\int_{d_a}^{\infty}u(x)dx - \int_{d_b}^\infty u(x)dx\Big| = \Big|\int_\B{R}u(x)dx - a - \int_\B{R}u(x)dx + b\Big| = |a-b|.$$

For constant $c>0$, note that $D^{L}_a(cu) = cu(x)\is_{x\geq d_a}$, where $\int_{d_a}^\infty cu(x)dx = \int_{\B{R}}cu(x)dx - a$ which implies that $\int_{d_a}^\infty u(x)dx = \int_{\B{R}}u(x)dx - a/c$. So $D^L_a(cu)=cD^L_{a/c}u$. \end{proof}

\begin{proof} (of Proposition \ref{detBoundDist}) Define
\begin{align*}u_k &= \left(e^\delta G_\delta D^L_{p(1-e^{-\delta})}D^R_{(1-p)(e^\delta -1)}\right)^{k-1} e^\delta G_\delta D^L_{p(1-e^{-\delta})}\rho \\
v_k &= \left(e^\delta G_\delta D^R_{(1-p)(1-e^{-\delta})}D^L_{p(e^\delta -1)}\right)^{k-1} e^\delta G_\delta D^R_{(1-p)(1-e^{-\delta})}\rho,\end{align*}
so we can write:
$$S_{k\delta}^{\delta,-}\rho = D^R_{(1-p)(e^\delta-1)}u_k \quad ; \quad S_{k\delta}^{\delta,+}\rho = D^L_{p(e^\delta - 1)}v_k.$$

Then we have
\begin{align*}
    \| u_k - v_k\|_1 &=\| e^\delta G_\delta D^L_{p(1-e^{-\delta})}D^R_{(1-p)(e^\delta -1)}u_{k-1} - e^\delta G_\delta D^R_{(1-p)(1-e^{-\delta})}D^L_{p(e^\delta-1)}v_{k-1}\|_1 && \text{ by definition}\\
    &\leq e^\delta \| D^L_{p(1-e^{-\delta})}D^R_{(1-p)(e^\delta-1)}u_{k-1} - D^R_{(1-p)(1-e^{-\delta})}D^L_{p(e^\delta-1)}v_{k-1}\|_1 && \text{ by \textit{(c)}}\\
    &\leq e^\delta \| D^L_{p(1-e^{-\delta})}D^R_{(1-p)(e^\delta-1)}u_{k-1} - D^L_{p(e^\delta-1)}D^R_{(1-p)(e^{\delta}-1)}u_{k-1}\|_1 \\
    &\quad \quad \quad \quad + e^\delta\| D^L_{p(e^\delta-1)}D^R_{(1-p)(e^{\delta}-1)}u_{k-1}- D^L_{p(e^\delta-1)}D^R_{(1-p)(1-e^{-\delta})}v_{k-1}\|_1 &&\text{ by \textit{(d)}, $\Delta$-ineq.}\\
    &\leq e^\delta\left( p(e^\delta - 1)-p(1-e^{-\delta})\right) + e^\delta \|D^R_{(1-p)(e^\delta - 1)}u_{k-1} - D^R_{(1-p)(1-e^{-\delta})}v_{k-1}\|_1 &&\text{ by \textit{(a),(b)}}\\
    &\leq p(e^\delta - 1)^2 + e^\delta \|D^R_{(1-p)(e^\delta - 1)}u_{k-1} - D^R_{(1-p)(1-e^{-\delta})}u_{k-1}\|_1 \\
    &\quad \quad \quad \quad + e^\delta \|D^R_{(1-p)(1-e^{-\delta})}u_{k-1} - D^R_{(1-p)(1-e^{-\delta})}v_{k-1}\|_1 &&\text{ by $\Delta$-ineq.}\\
    &\leq (e^\delta - 1)^2 + e^\delta \|u_{k-1} - v_{k-1}\|_1 &&\text{ by \textit{(a),(b)}}
\end{align*}
Therefore by recursion, since $(1+e^\delta + \cdots +e^{\ell \delta})=(e^{(\ell+1)\delta}-1)/(e^\delta-1)$, we can bound
$$\|u_k-v_k\|_1 \leq (e^\delta-1)(e^{(k-1)\delta}-1)+e^{(k-1)\delta}\|u_1-v_1\|_1$$
therefore
\begin{align*}
    \|S_{k\delta}^{\delta,-}\rho &- S_{k\delta}^{\delta,+}\rho\|_1 = \|D^R_{(1-p)(e^\delta - 1)}u_k - D^L_{p(e^\delta - 1)}v_k\|_1 \\
    &\leq \|D^R_{(1-p)(e^\delta - 1)}u_k - u_k\|_1 + \|u_k - v_k\|_1 + \|v_k - D^L_{p(e^\delta - 1)}v_k\|_1 && \text{ by $\Delta$-ineq.}\\
    &\leq (1-p)(e^\delta - 1) + \Big((e^\delta-1)(e^{(k-1)\delta}-1)+e^{(k-1)\delta}\|u_1-v_1\|_1\Big) + p(e^\delta - 1) && \text{ by \textit{(a)}} \\
    & = e^{(k-1)\delta}(e^\delta - 1) + e^{(k-1)\delta} \|e^\delta G_\delta D^L_{p(1-e^{-\delta})}\rho - e^\delta G_\delta D^R_{(1-p)(1-e^{-\delta})}\rho\|_1 &&\text{ by def. of $u_1,v_1$}\\
    & \leq e^{(k-1)\delta}(e^\delta - 1) + e^{k\delta} \|D^L_{p(1-e^{-\delta})}\rho - \rho\|_1 + e^{k\delta} \|\rho - D^R_{(1-p)(1-e^{-\delta})}\rho\|_1 &&\text{ by \textit{(c)} \& $\Delta$-ineq.}\\
    & = e^{(k-1)\delta}(e^\delta - 1) + e^{(k-1)\delta}(e^\delta - 1) = 2e^{(k-1)\delta}(e^\delta - 1) &&\text{ by \textit{(a)}}
\end{align*}
as required. 
\end{proof}

We will now prove the following Proposition, which is the same as Proposition 13 from \cite{hydroNBBM}, adapted to our setting. Analogously to our ordering $\oop$ for sets of real numbers (i.e. $A\oop B$ if and only if $|\{a\in A:a\geq x\}|\leq |\{b\in B:b\geq x\}|$ for all $x\in\B{R}$) we make a similar definition for positive and integrable real valued functions. We say $f\oop g$ if and only if $\int_x^\infty f(y)dy \leq \int_x^\infty g(y)dy$  for all $x\in \B{R}$. Then we can prove that the following ordering (in the sense of $\oop$) holds for the functions $S_{k\delta}^{\delta,\pm}\rho$:

\vspace{1em}

\begin{propn} \label{halfTime} For every $\delta > 0, k\in \B{N}$, we have
$$S_{k\delta}^{\delta,-}\rho\overset{(i)}{\oop} S_{k\delta}^{\delta/2,-}\rho \overset{(ii)}{\oop} S_{k\delta}^{\delta/2,+}\rho \overset{(iii)}{\oop} S_{k\delta}^{\delta,+}\rho $$
\end{propn}

Before the proof, we need the following Lemma:

\vspace{1em}

\begin{lemma} \label{orderSwitching} For any $0\leq a \leq \int_{\B{R}}w(s)ds$ and $\delta >0$:
    $$G_\delta D_a^L w \oop D_a^L G_\delta w\quad \quad D_a^R G_\delta w \oop G_\delta D_a^R w$$
\end{lemma}

\begin{proof}
    Let $q$ be the cut point such that $D^L_a G_\delta w(x)=\is_{x\geq q}G_\delta w(x)$. So then for $y\leq q$ we have
    $$\int_y^\infty D_a^L G_\delta w(s)ds = \int_{\B{R}}w(s)ds - a = \int_{\B{R}}G_\delta D_a^L w(s)ds \geq \int_y^\infty G_\delta D_a^L w(s)ds,$$
    and for $y\geq q$, since $D_a^L w(s)\leq w(s)$ point-wise, we have
    $$\int_y^\infty D_a^L G_\delta w(s)ds = \int_y^\infty G_\delta w(s)ds \geq \int_y^\infty G_\delta D_a^L w(s)ds.$$
    Now note that by symmetry 
    $$\int_y^\infty G_\delta D^L_a w(s)ds \leq \int_y^\infty D^L_a G_\delta w(s)ds \; \forall y \Leftrightarrow \int_{-\infty}^y G_\delta D^R_a w(s)ds \leq \int_{-\infty}^y D^R_a G_\delta w(s)ds\; \forall y,$$
    therefore since $\int_{-\infty}^\infty D^R_a G_\delta w = \int_{-\infty}^\infty G_\delta D^R_a w$, we have
    $$\int_y^\infty D_a^R G_\delta w = \int_{-\infty}^\infty D_a^R G_\delta w - \int_{-\infty}^y D_a^R G_\delta w \leq \int_{-\infty}^\infty G_\delta D_a^R w - \int_{-\infty}^y G_\delta D^R_a w =\int_y^\infty G_\delta D^R_a w.$$
    Therefore we have $G_\delta D^L_a w \oop D^L_a G_\delta w$ and $D^R_a G_\delta w \oop G_\delta D^R_a w$ as required. 
\end{proof}

\begin{proof} \textit{(Of Proposition \ref{halfTime}}) First we note that the comparison $\overset{(ii)}{\oop}$ follows immediately from Propositions \ref{upDownCoupling} and \ref{couplingConv}

To prove $\overset{(i)}{\oop}$, define $H_t\rho := \left((D^R_{(1-p)(e^t-1)}e^t G_t D^L_{p(1-e^{-t})}\right)\rho$ and note that
\begin{align*}
    G_{\frac{\delta}{2}}&D^L_{p(e^{-\frac{\delta}{2}}-e^{-\delta})}D^L_{p(1-e^{-\frac{\delta}{2}})}\rho \oop D^L_{p(e^{-\frac{\delta}{2}}-e^{-\delta})}G_{\frac{\delta}{2}}D^L_{p(1-e^{-\frac{\delta}{2}})}\rho &&\text{ by \textit{(e)} \& Lemma \ref{orderSwitching}}\\
    &\Rightarrow e^{\frac{\delta}{2}}G_{\frac{\delta}{2}}D^L_{p(1-e^{-\delta})}\rho \oop D^L_{p(1-e^{-\frac{\delta}{2}})}e^{\frac{\delta}{2}}G_{\frac{\delta}{2}}D^L_{p(1-e^{-\frac{\delta}{2}})}\rho &&\text{ by \textit{(e)}\&\textit{(f)}}\\
    &\Rightarrow D^R_{(1-p)(e^{\frac{\delta}{2}}-1)}G_{\frac{\delta}{2}}e^{\frac{\delta}{2}}G_{\frac{\delta}{2}}D^L_{p(1-e^{-\delta})}\rho \oop G_{\frac{\delta}{2}}D^R_{(1-p)(e^{\frac{\delta}{2}}-1)}D^L_{p(1-e^{-\frac{\delta}{2}})}e^{\frac{\delta}{2}}G_{\frac{\delta}{2}}D^L_{p(1-e^{-\frac{\delta}{2}})}\rho &&\text{ by Lemma \ref{orderSwitching}}\\
    &\Rightarrow D^R_{(1-p)(e^{\frac{\delta}{2}}-1)}e^{\frac{\delta}{2}}G_{\delta}D^L_{p(1-e^{-\delta})}\rho \oop G_{\frac{\delta}{2}}D^L_{p(1-e^{-\frac{\delta}{2}})}H_{\frac{\delta}{2}}\rho &&\text{ by \textit{(d)}}\\
    &\Rightarrow D^R_{(1-p)(e^{\frac{\delta}{2}}-1)}e^{\frac{\delta}{2}}D^R_{(1-p)(e^{\frac{\delta}{2}}-1)}e^{\frac{\delta}{2}}G_{\delta}D^L_{p(1-e^{-\delta})}\rho \oop D^R_{(1-p)(e^{\frac{\delta}{2}}-1)}e^{\frac{\delta}{2}}G_{\frac{\delta}{2}}D^L_{p(1-e^{-\frac{\delta}{2}})}H_{\frac{\delta}{2}}\rho &&\text{ by \textit{(f)}}\\
    &\Rightarrow H_\delta \rho \oop (H_{\frac{\delta}{2}})^2\rho &&\text{ by \textit{(e)}\&\textit{(f)}}
\end{align*}

Thus by induction it holds that
$$S_{k\delta}^{\delta,-}\rho =\left(H_\delta\right)^k \rho \oop \left(\left(H_{\delta/2}\right)^2\right)^k \rho = S^{\delta/2,-}_{k\delta} \rho,$$
and by the same argument it holds that $S^{\delta/2,+}_{k\delta}\rho \oop S^{\delta,+}_{k\delta}\rho$. Thus we have proven $\overset{(i)}{\oop}$. Then $\overset{(iii)}{\oop}$ follows by symmetry.
\end{proof}

 Next we give a simple analytical result about uniform convergence of cumulative distribution functions. The proof is simple and can be found in the appendix. 

\vspace{1em}

\begin{lemma} \label{functionLims} Let $(f_n(x))_{n\geq 1}$ be a sequence of non-negative functions such that $\int_{\B{R}}f_n(s)ds=1\; \forall n$. Suppose that the sequence is also monotonic increasing (resp. decreasing) in the sense that $m < n \implies \int_x^\infty f_m(s)ds \leq \int_x^\infty f_n(s)ds$ for all $x$ (resp. $\int_x^\infty f_n(s) \leq \int_x^\infty f_m(s)ds$ for all $x$), and that the sequence is bounded above (resp. below) by a bounded function $b$ (in the sense that $f_n \oop b$) such that $b$ is non-negative and $\int_{\B{R}}b(s)ds =1$. Then there exists a function $\Psi$ such that
$$\sup_{x\in \B{R}}\left| \int_x^\infty f_n(s)ds - \Psi(x)\right|\xrightarrow[]{n\to\infty}0.$$
\end{lemma}

This allows us to prove the following.

\vspace{1em}

\begin{propn} \label{sameLimitPropn} $\int_x^\infty S_{t}^{t/2^n,-}\rho(y)dy$ and $\int_x^\infty S_{t}^{t/2^n,+}\rho(y)dy$ both converge to the same limit as $n\to\infty$.
\end{propn}

\begin{proof} By Proposition \ref{halfTime}, the sequence $(S_{t}^{t/2^n, -}\rho)_{n\in \B{N}}$ is monotonic increasing and bounded above by $S^{t,+}_{t}\rho$, and $(S_{t}^{t/2^{n}, +}\rho)_{n\in \B{N}}$ is monotonic decreasing and bounded below by $S^{t,-}_{t}\rho$. Therefore by Lemma \ref{functionLims}, both sequences converge uniformly to limits, say $\Psi^-$ and $\Psi^+$ respectively. Moreover, since $\int_x^\infty S_t^{t/2^n,-}\rho(y)dy$ is monotonic increasing in $n$ and $\int_x^\infty S_t^{t/2^n,+}\rho(y)dy$ is monotonic decreasing in $n$, thus $\|\Psi^+(x)-\Psi^-(x)\|_1$ is bounded above by $|\int_x^\infty S_t^{t/2^n,+}\rho(y)dy - \int_x^\infty S_t^{t/2^n,-}\rho(y)dy|$ for all $n$. By Proposition \ref{detBoundDist}, this upper bound converges to $0$ for fixed $t$ as $n\to\infty$, therefore $\|\Psi^+(x)-\Psi^-(x)\|_1=0$ and hence $\Psi^+(x)\equiv \Psi^-(x)$, thus confirming that $\int_x^\infty S_t^{t/2^n,-}\rho(y)dy$ and $\int_x^\infty S_t^{t/2^n,+}\rho(y)dy$ converge to the same limit. 
\end{proof}

With this result, if we can show that any solution $u$ of the FBP \eqref{fbpMain} is bounded $S_{t}^{t/n,-}\rho(x) \oop u(x,t)\oop S_{t}^{t/n,+}\rho(x)$, then $u$ must be the limit of $S_{t}^{t/n,-}$ and $S_t^{t/n,+}$ as $n\to\infty$. We will show this in Section \ref{comparison}, but before then we must show that the FBP does in fact have \textit{at least one} solution $u$. This we do in the Section \ref{probRepOfSolutions}.

\section{Probabilistic representation of free boundary problem solutions} \label{probRepOfSolutions}
In this section, we prove that there exists a solution of the 2-sided FBP \eqref{fbpMain}. We do so by appealing to a connection with the \textit{inverse first passage problem}.

Suppose you are given a probability distribution $\mu$ on $[0,\infty)$. The (one-sided) inverse first passage problem asks whether one can find a boundary $L:[0,\infty)\to\B{R}$ such that $\tau^L:=\inf\{t>0:W_t\leq L_t\}\sim \mu$. It is known that for $\mu=Exp(1)$ there does in fact exist such an $L$ (Theorem 1, \cite{cccs11}, applying $p(t)=e^{-t}$). Moreover, it was proven by Chen et al. (Theorem 1 of \cite{ccs22}) that the solution $L$ is $C^\infty((0,\infty))$. Therefore given an initial density $\rho$ of the Brownian motion $W$, consider finding boundaries $L$ and $R$ such that: 
\begin{equation}\label{tauLtauR}\begin{split} \tau^L&:=\inf\{t>0:W_t\leq L_t\}\sim Exp(p),\\
\tau^R&:=\inf\{t>0:W_t\geq R_t\}\sim Exp(1-p).\end{split}\end{equation}
If we define $\tau:=\tau^L\wedge \tau^R$ then $\B{P}_\rho(\tau>t)=e^{-t}$ and $\B{P}_\rho(W_\tau=L_\tau)=p=1-\B{P}_\rho(W_\tau=R_\tau)$, or equivalently $\B{P}_\rho(\tau\leq t,W_\tau=L_\tau)=p(1-e^{-t})$ and $\B{P}_\rho(\tau\leq t,W_\tau=R_\tau)=(1-p)(1-e^{-t})$. We may say that $(L,R)$ thus solves the asymmetric two-sided inverse first passage problem. From this solution $(L,R)$ we may construct a solution of \eqref{fbpMain}. 

\vspace{1em}

\begin{propn} \label{repThem}
    Let $\rho$ be a probability density function and let $\tau_L$, $\tau_R$, $L$, and $R$ be defined by \eqref{tauLtauR}. Define $\tau=\tau^L\wedge \tau^R$ and  
    $$u(x,t)=e^t \frac{\partial}{\partial x}\B{P}_\rho(W_t\leq x,\tau>t).$$
    In other words, $u(\cdot,t)$ is the conditional density of a Brownian motion with intial density $\rho$ conditioned not to exit $(L_s,R_s)$ before time $t$. Then $u$ is a solution to the FBP \eqref{fbpMain}. 
\end{propn}

First, we prove the following Lemma, which shows that if a Brownian motion has initial density $f$ with support $(L_0,R_0)$ and has left and right derivatives at $R_0$ and $L_0$ respectively, then the probability of the Brownian motion exiting on the left or right is proportional to $f'$ at that boundary. We will subsequently apply this Lemma with $f(x)=u(x,t)$ for $t>0$, since Theorem 1 of \cite{ccs22} tells us that the boundaries $L$ and $R$ are $C^\infty((0,\infty))$, and hence Lipschitz, around $t$ for any $t>0$. 

\vspace{1em}

\begin{lemma} \label{erfcPropn} Let $s\mapsto L_s$ and $s\mapsto R_s$ be barriers which do not touch and are Lipschitz continuous in some neighbourhood of $0$. Let $f(x)$ be a continuous function with support $(L_0,R_0)$. Furthermore, suppose that $f$ has a right derivative at $L_0$ and a left derivative at $R_0$. For Brownian motion $W$ with initial density $f$, let $\tau:=\inf\{s\geq 0:W_s \notin (L_s,R_s)\}$ be the first exit time of $(L_s,R_s)$ by $W_s$. Then:
\begin{align*}\lim_{\delta\downarrow 0}\frac{1}{\delta}\int_{L_0}^{R_0}f(y)\B{P}_y(\tau<\delta, W_{\tau}=L_{\tau})dy=\frac12\lim_{x\downarrow L_0}f'(x)=: f^L\\
\lim_{\delta\downarrow 0}\frac{1}{\delta}\int_{L_0}^{R_0}f(y)\B{P}_y(\tau<\delta, W_{\tau}=R_{\tau})dy=-\frac12\lim_{x\uparrow R_0}f'(x)=: f^R\end{align*}
\end{lemma}

\begin{proof}
We begin by noting that as $L,R$ are Lipschitz continuous at $0$, there exists $C>0$ such that for $\delta$ sufficiently small, $|L_\delta-L_0|<C\delta$ and $|R_\delta-R_0|<C\delta$, then we have:
$$\B{P}_y(\tau<\delta,W_{\tau}=L_{\tau})\leq \B{P}_y\left(\min_{0\leq s\leq \delta}W_s < L_0 + C\delta\right),$$
and:
\begin{align*}\B{P}_y(\tau<\delta,W_{\tau}=L_{\tau})&\geq \B{P}_y\left(\max_{0\leq s\leq \delta}W_s - R_s < 0, \min_{0\leq s\leq \delta}W_s  - L_s <0\right) \\ &\geq \B{P}_y\left(\min_{0\leq s\leq \delta}W_s < L_0 - C\delta\right)-\B{P}_y\left(\max_{0\leq s\leq \delta}W_s > R_0 - C\delta\right),
\end{align*}
therefore putting these bounds together and using the reflection principle, we have:
$$\B{P}_y(\tau<\delta, W_{\tau}=L_{\tau})\in \left[\max\left\{0,\erfc\left(\frac{y-L_0+C\delta}{\sqrt{2\delta}}\right)-\erfc\left(\frac{R_0-C\delta-y}{\sqrt{2\delta}}\right)\right\}, \erfc\left(\frac{y-L_0-C\delta}{\sqrt{2\delta}}\right)\right],$$
where $\erfc(z)$ is the complementary error function defined by $\erfc(z)=\frac{2}{\sqrt{\pi}}\int_z^\infty e^{-s^2}ds$. Let us consider the integral $\frac{1}{\delta}\int_{L_0}^{R_0}f(y)\erfc\left(\frac{y-L_0\pm C\delta}{\sqrt{2\delta}}\right)dy$ first. We may split up the integral as:
\begin{align}\nonumber \frac{1}{\delta}\int_{L_0}^{R_0} f(y)&\erfc\left(\frac{y-L_0 \pm C\delta}{\sqrt{2\delta}}\right)dy\\
\label{upperBoundIntegralSplitup} &=\frac{1}{\delta}\int_{L_0}^{L_0+S(\delta)}f(y)\erfc\left(\frac{y-L_0\pm C\delta}{\sqrt{2\delta}}\right)dy + \frac{1}{\delta}\int_{L_0+S(\delta)}^{R_0} f(y)\erfc\left(\frac{y-L_0\pm C\delta}{\sqrt{2\delta}}\right)dy,\end{align}
where $S(\delta)$ is some function which we are yet to determine which we think of as being small. Consider the first term of (\ref{upperBoundIntegralSplitup}). Changing variables so that $z=(2\delta)^{-1/2}(y-L_0\pm C\delta)$, and Taylor expanding $f$ about $L_0$, we have:
$$\frac{1}{\delta}\int_{L_0}^{L_0+S(\delta)}f(y)\erfc\left(\frac{y-L_0\pm C\delta}{\sqrt{2\delta}}\right)dy=\sqrt{\frac{2}{\delta}}\int_{\pm \frac{C\delta}{\sqrt{2\delta}}}^{\frac{S(\delta)\pm C\delta}{\sqrt{2\delta}}}\left(\mp f^L C\delta+f^L\sqrt{2\delta}z+O((C\delta+\sqrt{2\delta}z)^2)\right)\erfc(z)dz.$$
Then as $\int_0^\infty \erfc(x)dx<\infty$, and $\int_0^\infty x\erfc(x)dx=1/4$, then certainly if $S(\delta)/\sqrt{2\delta}\to\infty$, we have: 
$$\sqrt{\frac{2}{\delta}}\int_{\pm\frac{C\delta}{\sqrt{2\delta}}}^{\frac{S(\delta)\pm C\delta}{\sqrt{2\delta}}}(\mp f^LC\delta +f^L\sqrt{2\delta}x)\erfc(x)dx\to \frac12f^L.$$ 
as $\delta\to 0$. Next we deal with the $O((C\delta+\sqrt{2\delta}z)^2)\erfc(z)$ term of the integral. For $z\in \left(\pm C\delta/\sqrt{2\delta},(S(\delta)\pm C\delta)/\sqrt{2\delta}\right)$, we have $z=O\left(\max\{\delta^{1/2},S(\delta)/\delta^{1/2}\right\})$, therefore $\delta z=O(\max\{\delta^{3/2},\delta^{1/2}S(\delta)\})$ and $\delta^{1/2}z^2 = O(\max\{\delta^{3/2},S(\delta)^2/\delta^{1/2}\})$. Therefore if $S(\delta)^2/\sqrt{2\delta} \to 0$, then: 
\begin{align*}
\sqrt{\frac{2}{\delta}}\int_{\pm\frac{C\delta}{\sqrt{2\delta}}}^{\frac{S(\delta)\pm C\delta}{\sqrt{2\delta}}}&O((C\delta+\sqrt{2\delta}z)^2)\erfc(x)dx=\delta^{-1/2}O\left(\max\{\delta^2, \delta^{3/2}z, \delta z^2\}\right) \\
=&O\left(\max\{\delta^{3/2},\delta z, \delta^{1/2}z^2\}\right)=O\left(\max\{\delta^{3/2},\delta^{1/2}S(\delta),S(\delta)^2/\delta^{1/2}\}\right) \to 0.
\end{align*}
Therefore choosing $S(\delta)=\delta^{3/8}$ and putting the two limits above together, we get that
$$\frac{1}{\delta}\int_{L_0}^{L_0+S(\delta)} f(y)\erfc\left(\frac{y-L_0\pm C\delta}{\sqrt{2\delta}}\right)dy \to \frac12f^L,$$

Now consider the second term of (\ref{upperBoundIntegralSplitup}). Once again chainging variables, we can upper bound the integral for sufficiently small $\delta$, by \begin{align*}\frac{1}{\delta}\int_{L_0+S(\delta)}^{R_0}f(y)&\erfc\Big(\frac{y-L_0\pm C\delta}{\sqrt{2\delta}}\Big)dy \leq \frac{\sqrt{2}\|f\|_{\infty}}{\sqrt{\delta}}\int_{\frac{\delta^{-1/8}\pm C\delta^{1/2}}{2^{1/2}}}^{\infty}\erfc(u)du \\&= \frac{2\sqrt{2}\|f\|_{\infty}}{\sqrt{\pi\delta}}\int_{\frac{\delta^{-1/8}\pm C\delta^{1/2}}{2^{1/2}}}^{\infty}\int_u^\infty e^{-t^2}dt\,du
\leq \frac{2\sqrt{2}\|f\|_{\infty}}{\sqrt{\pi\delta}}\int_{\frac{\delta^{-1/8}\pm C\delta^{1/2}}{2^{1/2}}}^{\infty}\int_u^\infty te^{-t^2}dt\,du \\&= \frac{2\sqrt{2}\|f\|_{\infty}}{\sqrt{\pi\delta}}\int_{\frac{\delta^{-1/8}\pm C\delta^{1/2}}{2^{1/2}}}^{\infty}e^{-u^2}du \leq \frac{2\sqrt{2}\|f\|_{\infty}}{\sqrt{\pi\delta}}\int_{\frac{\delta^{-1/8}\pm C\delta^{1/2}}{2^{1/2}}}^{\infty}ue^{-u^2}du \\
&= O(\delta^{-1/2}e^{-\delta^{-1/4}/2})\xrightarrow[\delta\to 0]{}0,
\end{align*}
so putting together both terms of (\ref{upperBoundIntegralSplitup}), we have that $\frac{1}{\delta}\int_{L_0}^{R_0} f(y)\erfc\left(\frac{y-L_0 \pm C\delta}{\sqrt{2\delta}}\right)dy\xrightarrow[\delta\to 0]{}\frac12f^L$. We can similarly split up the integral:
\begin{align*}\frac{1}{\delta}\int_{L_0}^{R_0} f(y)&\erfc\left(\frac{R_0-C\delta - y}{\sqrt{2\delta}}\right)\\
&=\frac{1}{\delta}\int_{L_0}^{R_0-S(\delta)}f(y)\erfc\left(\frac{R_0-C\delta-y}{\sqrt{2\delta}}\right)dy + \frac{1}{\delta}\int_{R_0-S(\delta)}^{R_0} f(y)\erfc\left(\frac{R_0-C\delta - y}{\sqrt{2\delta}}\right)dy,\end{align*}
where the first term converges to $0$ as $\delta \to 0$ and the second term converges to $\frac12f^R$. But notice that when $y\geq (R_0-L_0)/2-C\delta$, then $\erfc\left(\frac{y-L_0+C\delta}{\sqrt{2\delta}}\right)<\erfc\left(\frac{R_0-C\delta-y}{\sqrt{2\delta}}\right)$, therefore the lower bound for $\B{P}_y(\tau<\delta, W_{\tau}=L_{\tau})$ is $0$. Therefore if $\delta$ is chosen sufficiently small so that $S(\delta)=\delta^{3/8}<(R_0-L_0)/2-C\delta$, then certainly:
$$\frac{1}{\delta}\int_{L_0}^{R_0} f(y)\max\left\{0,\erfc\left(\frac{y-L_0+C\delta}{\sqrt{2\delta}}\right)-\erfc\left(\frac{R_0-C\delta-y}{\sqrt{2\delta}}\right)\right\}dy \to \frac12f^L,$$
and thus by sandwiching the result holds.
\end{proof}

Using this, we now prove Proposition \ref{repThem}. 

\begin{proof}\textit{(of Proposition \ref{repThem})}
Let $L$ and $R$ be defined as in \eqref{tauLtauR}. By the Feynman-Kac representation theorem (see Chapter 7, Theorem 4.2 of \cite{yongZhou}, for example) we have that 
\begin{align*} u(x,t) = e^t\B{E}_x[\rho(W_t)\is_{\tau>t}]\end{align*}
solves $u_t = \frac12u_{xx}+u$ for $x\in(L_t,R_t)$ and $u=0$ for $x\notin (L_t,R_t)$. Let $\tilde{p}_t(y,z)$ be the sub-probability density of $\{W_t=z,\tau>t|W_0=y\}$; that is, $\tilde{p}_t(y,z)=\frac{\partial}{\partial x}\B{P}_y(W_t\leq x,\tau>t)|_{x=z}$, and note that by the reversibility of Brownian motion $\tilde{p}_t(y,z)=\tilde{p}_t(z,y)$. Therefore \begin{align}\label{feynmanKacRep}u(x,t)=e^t\int_{\B{R}}\tilde{p}_t(x,y)\rho(y)dy=e^t\int_{\B{R}}\rho(y)\tilde{p}_t(y,x)dy=e^t\frac{\partial}{\partial x}\int_{\B{R}}\rho(y)\B{P}_y(W_t\leq x,\tau>t)dy.\end{align}
Then conditioning on the value of $W_t$ and applying Lemma \ref{erfcPropn}, we have that
\begin{align} \label{connectU_zTop}
    \nonumber \lim_{\delta\downarrow 0}\frac1\delta \int_\B{R}\rho(y)\B{P}_y(\tau\in &[t,t+\delta),W_{\tau}=L_{\tau})dy\\
    &= \lim_{\delta\downarrow 0}\frac1\delta \int_{\B{R}}\B{P}_z(\tau<\delta,W_{\tau}=L_{\tau})\int_{\B{R}}\rho(y)\tilde{p}_t(y,z)dy dz \\
    \nonumber &= \lim_{\delta\downarrow 0}\frac1\delta \int_{\B{R}}\B{P}_z(\tau<\delta, W_{\tau}=L_{\tau})e^{-t}u(z,t)dz = \lim_{z\downarrow L_t}\frac12 u_x(z,t)e^{-t}.
\end{align}
Since 
\begin{align*}\B{P}_\rho(\tau\in [t,t+\delta),W_{\tau}=L_{\tau})&=\B{P}_\rho(\tau<t+\delta, W_{\tau}=L_{\tau})-\B{P}_\rho(\tau<t+\delta, W_{\tau}=L_{\tau})\\
&=p(1-e^{-(t+\delta)})-p(1-e^{-t})=-p(e^{-(t+\delta)}-e^{-t}),
\end{align*}
therefore 
$$\lim_{\delta\downarrow 0}\frac{1}{\delta}\int_\B{R}\rho(y)\B{P}_y(\tau\in [t,t+\delta),W_{\tau}=L_{\tau})dy=\lim_{\delta\downarrow 0}\frac{-p(e^{-(t+\delta)}-e^{-t})}{\delta},$$
which we can identify as $-p\frac{d}{dt}e^{-t}=pe^{-t}$. Therefore by \eqref{connectU_zTop}, $u_x(L_t,t)=2p$. An identical argument shows that $u_x(R_t,t)=2(p-1)$. Therefore $u$ solves $u_t=\frac12 u_{xx}+u$ for $x\in (L_t,R_t)$, $u=0$ for $x\notin (L_t,R_t)$, $u$ is continuous, $u(x,0)=\rho(x)$, and the boundary conditions $u_x(L_t,t)=2p$ an $u_x(R_t,t)=2(p-1)$ hold. Thus $u$ is a solution to the FBP \eqref{fbpMain}.
\end{proof}

\section{A comparison theorem for the free boundary problem} \label{comparison}
The aim of this section is to prove that for any $t=k\delta>0$, any solution $u(x,k\delta)$ to the free boundary problem \eqref{fbpMain} with initial condition $\rho$ is sandwiched in between the functions $S_{\delta}^{k\delta,-}\rho(x)$ and $S_{k\delta}^{\delta,-}\rho(x)$ in the sense that
\begin{align} \label{compSectionAim}
\int_a^\infty S_{t}^{t/n,-}\rho(x)dx \leq \int_a^\infty u(x,t)dx\leq \int_a^\infty S_{t}^{t/n,+}\rho(x)dx.\end{align}
With Proposition \ref{detBoundDist}, this will allow us to show that the upper and lower bounds $S_{\delta}^{k\delta,\pm}\rho(x)$ converge to the same limit $u$, and as a consequence that there is only one solution $u$ of the FBP \eqref{fbpMain}. 
 
We first prove this for the case that $k=1$ and for the lower bound. That is, we prove that, for any fixed $\delta>0$, we have \begin{align}\label{Sinequ}\int_a^\infty C^R_1e^\delta G_\delta C^L_{1-p(1-e^{-\delta})}\rho(dr)\leq \int_a^\infty u(x,\delta)dx.\end{align} 
We can then inductively reapply the proposition to show that for any $k$, then inequality holds for $\int_a^\infty S_{k\delta}^{\delta,-}\rho(dr)\leq \int_a^\infty u(dr,t)$. By symmetry, the upper bound similarly holds, so that inequality \eqref{compSectionAim} holds for any $t>0$ and $n\in \B{N}$.

\vspace{1em}

\begin{propn} \label{comparisonTheorem} Let $u(x,t)$ be a solution to the free boundary problem \eqref{fbpMain} with initial condition $\rho(x)$ such that the boundaries $L,R$ are locally Lipschitz continuous on $(0,\infty)$. Let $v(x,t)$ be the solution to the PDE which satisfies $v_t = \frac12 v_{xx} + v$ and has initial condition $\rho_0(x):=C^L_{1-p(1-e^{-\delta})}\rho(x)$. Then for any $r$:
\begin{equation} \label{npIneq}
    1 \wedge \int_{-\infty}^r v(x,\delta) dx \geq \int_{-\infty}^r u(x,\delta)dx
\end{equation}
\end{propn}

We will prove the above Proposition in a similar fashion to the Proof of 10.3.13 in \cite{carinci}, by constructing 2 simpler particle systems, finding their hydrodynamic limits, and constructing a coupling between the two systems. 

\begin{proof} We begin by noting that by Theorem \ref{repThem}, the inequality \eqref{npIneq} is equivalent to the inequality:
\begin{equation}
    1 \wedge e^\delta \int \rho_0(x) \B{P}_x(B_\delta \leq r) dx \geq e^\delta \int \rho(x) \B{P}_x(B_\delta \leq r, \tau^{LR} > \delta) dx
\end{equation}
where $\tau^{LR}:=\inf\{s>0:B_s \notin (L_s,R_s)\}$ and $L_s,R_s$ are locally Lipschitz continuous curves such that for all $t\in [0,\delta]$
\begin{equation}\label{leftBdyCdn} e^t \int \rho(x) \B{P}_x(\tau^{LR} < t, B_{\tau^{LR}}=L_{\tau^{LR}})dx = p(e^t-1).\end{equation}
\begin{equation}\label{rightBdyCdn} e^t \int \rho(x) \B{P}_x(\tau^{LR} < t, B_{\tau^{LR}}=R_{\tau^{LR}})dx = (1-p)(e^t-1).\end{equation}
Since $e^\delta\int \rho(x)\B{P}_x(\tau^{LR}>\delta)dx = 1$, it suffices to show:
$$e^\delta \int \rho_0(x) \B{P}_x(B_\delta \leq r) dx \geq e^\delta \int \rho(x) \B{P}_x(B_\delta \leq r, \tau^{LR} > \delta) dx$$
So fix $n$ and let $M:=\lceil p(e^\delta -1)n\rceil$ and $N:=\lfloor (1-p)(e^\delta-1)n \rfloor$. Let $B^{(1)}_i,i=M+1,...,M+N+n$ be $N+n$ standard Brownian motions, and independently for each $i$, distribute $B^{(1)}_i(0)$ with initial density $\rho_0/\int \rho_0$. Let $B^{(2)}_i, i=1,...,M+N+n$ be $M+N+n$ standard Brownian motions, which are removed when they hit $L_s$ or $R_s$. Couple the particles so that $B^{(1)}_i(0)=B^{(2)}_i(0)$ for $i=M+1,...,M+N+n$ and distribute the particles $B^{(2)}_i(0),i=1,...,M$ with density $(\rho-\rho_0)/\int(\rho-\rho_0)$. Thus as $n\to\infty$, it is simple to check that the empirical distribution of $B^{(2)}(0)$ converges to $\rho$. 

Now we will couple the particles over the time interval $[0,\delta]$ as follows. We will call the particles $B^{(1)}_i$ the $(1)$-particles, and $B^{(2)}_i$ the $(2)$-particles. Recall that $B^{(1)}_i(0)=B^{(2)}_i(0)$ for $i=M+1,\ldots,M+N+n$, and we will call $\{B^{(1)}_i(0),B^{(2)}_i(0)\}$ a `married couple'. Whilst 2 particles are `married', they will be driven by the same Brownian motion.
The remaining uncoupled particles will be called `singletons'. We will call a married couple $\{B^{(1)}_k(t),B^{(2)}_\ell(t)\}$ `good' if $B^{(1)}_k(t)\leq B^{(2)}_\ell(t)$ and `bad' otherwise.

So suppose that at time $\tau$ the particle $B^{(2)}_i$ exits $(L_\tau,R_\tau)$ and is removed from the system. Then we describe how the particles are recoupled.

\begin{enumerate}
    \item If $B^{(2)}_i(\tau-)$ is a singleton particle, then when $B^{(2)}_i(\tau)$ is removed, we do not change the coupling. The numbers of good and bad couples are unchanged.
    \item If $B^{(2)}_i(\tau-)$ is part of a married couple, say $\{B^{(1)}_k(\tau-),B^{(2)}_i(\tau-)\}$ is a married couple:
    \begin{enumerate}
        \item If there are still $\geq N+n$ $(2)$-particles remaining at time $\tau$:
        \begin{enumerate}
            \item If $B^{(2)}_i(\tau)$ is removed because it hits $L_\tau$, then take any remaining singleton $(2)$-particle, say $B_j^{(2)}(\tau)$, and couple it to $B^{(1)}_k(\tau)$ to form a married couple $\{B^{(1)}_k(\tau),B^{(2)}_j(\tau)\}$. Since $B_i^{(2)}(\tau-)<B_j^{(2)}(\tau-)$, we have either replaced a good couple by a good couple, or a bad couple with a good or bad couple. Therefore the number of bad couples does not increase.
            \item If $B^{(2)}_i(\tau)$ is removed because it hits $R_\tau$, then take any remaining singleton $(2)$-particle, say $B_j^{(2)}(\tau)$, and couple it to $B^{(1)}_k(\tau)$ to form a married couple $\{B^{(1)}_k(\tau),B^{(2)}_j(\tau)\}$. Therefore the number of bad couples may increase by at most 1.
        \end{enumerate}
        Thus in either case, whilst there are $\geq N+n$ $(2)$-particles remaining, all $(1)$-particles are married. 
        \item If there are $<N+n$ $(2)$-particles remaining
        \begin{enumerate}
            \item If $B^{(2)}_i(\tau)$ is removed because it hits $L_\tau$:
            \begin{itemize}
                \item If $B^{(2)}_i(\tau) < B^{(1)}_k(\tau)$ (ie. $B_i^{(2)}(\tau-)$ was in a bad couple), then leave $B^{(1)}_k$ as a singleton, thus decreasing the number of bad couples by 1. 
                \item If $B^{(1)}_k(\tau) \leq B^{(2)}_i(\tau)=L_\tau$, then consider any $(2)$-particle in a bad couple, and recouple it to $B^{(1)}_k$. Since all $(2)$-particles are in $(L_\tau,R_\tau)$ at time $\tau$, this must form a good couple, and therefore we have decreased the number of bad couples by 1.
            \end{itemize}
            In either case, the number of bad couples decreases by 1 (unless there are no bad couples). If there are no bad couples, then no recoupling happens at the time $\tau$.
            \item Suppose $B^{(2)}_i$ is removed because it hits $R_s$. Then when $B^{(2)}_i$ is removed, do not recouple. The number of married couples decreases by 1; the removed couple may be good or bad.
        \end{enumerate}
    \end{enumerate}
\end{enumerate}

Now define:
\begin{align*}
    \alpha_n(r) &:= \sum_{i=M+1}^{M+N+n} \is_{\{B_i^{(1)}(\delta)\leq r\}} \\
    \beta_n(r) &:= \sum_{i=1}^{M+N+n} \is_{\{\text{$B^{(2)}_i(\delta)$ has not been removed}\}}\is_{\{B_i^{(2)}(\delta)\leq r\}}
\end{align*}
Then by the law of large numbers
$$\B{E}\left[\frac{1}{N+n}\sum_{i=M+1}^{M+N+n}\is_{\{B_i^{(1)}(\delta)\leq r\}}\right] \xrightarrow{n\to\infty} \int_{\B{R}} \frac{\rho_0(x)}{\int \rho_0}\B{P}_x(B_\delta \leq r)dx .$$
Therefore since $(N+n)/n \to 1+(1-p)(e^\delta-1)=p+e^\delta-pe^\delta = e^\delta \int \rho_0$ as $n\to\infty$, we conclude that
$$\frac{1}{n}\B{E}[\alpha_n(r)] \xrightarrow[n\to\infty]{} e^\delta \int \rho_0(x)\B{P}_x(B_\delta \leq r)dx.$$
Again, by the strong law of large numbers, since $(M+N+n)/n \to e^\delta$ as $n\to\infty$, thus 
$$\frac{1}{n}\B{E}[\beta_n(r)] \xrightarrow[n\to\infty]{} e^\delta \int \rho(x) \B{P}_x(B_\delta \leq r, \tau^{LR} > \delta) dx $$

Therefore to prove the claim it suffices to show that $\lim_{n\to\infty}\frac1n \B{E}[\alpha_n(r)-\beta_n(r)]\geq 0$ almost surely for all $r$. Let $A$ be the number of $(2)$-particles at time $\delta$. Then by \eqref{leftBdyCdn}, \eqref{rightBdyCdn} $\frac{1}{n}\B{E}[A]=\frac{M+N+n}{n}(1-(1-e^{-\delta}))\xrightarrow{n\to\infty} 1$.

Notice that whilst there are $\geq N+n$ $(2)$-particles alive, the number of married couples is $N+n$, and the number of bad couples is at most the number of particles which hit $R_s$, which by (\ref{rightBdyCdn}) has expectation $(1-p)M$. Then when there are fewer than $N+n$ $(2)$-particles alive, we remove a bad couple when a $(2)$-particle hits $L_s$ (if any bad couples remain). By (\ref{leftBdyCdn}), the expected number of times this happens is $\B{E}[p(N+n-A)]$. Let $\mathfrak{B}$ be the number of bad couples at time $\delta$. Therefore the expected number of bad couples, $\B{E}\mathfrak{B}$ is at most $((1-p)M - p(N+n-A))$. Therefore since $(1-p)M/n$ and $pN/n$ both converge to $p(1-p)(e^\delta-1)$ as $n\to\infty$, thus $\frac1n\B{E}[\mathfrak{B}]\xrightarrow{n\to\infty} 0$. 

Now note that for each good couple $\{B_k^{(1)}(\delta),B_{\ell}^{(2)}(\delta)\}$, we have $\is_{\{B_k^{(1)}(\delta)\leq r\}}-\is_{\{B_l^{(2)}(\delta)\leq r\}}\geq 0$, and for every bad couple, $\is_{\{B_k^{(1)}(\delta)\leq r\}}-\is_{\{B_\ell^{(2)}(\delta)\leq r\}}\geq -1$. 
Therefore $$\alpha_n(r) - \beta_n(r)\geq \sum_{\text{married couples}}\left(\is_{\{B_k^{(1)}(\delta)\leq r\}}-\is_{\{B_l^{(2)}(\delta)\leq r\}}\right) - \frac1n \#\{\text{singleton $(2)$-particles}\}.$$

Let $S_2(n)$ be the number of singleton $(2)$-particles at time $\delta$. Now by (\ref{leftBdyCdn}) and (\ref{rightBdyCdn}), each $(2)$-particle is, independently, alive at time $\delta$ with probability $e^{-\delta}$. Therefore by the central limit theorem, $\frac{S_2(n) - n}{\sqrt{ne^{-\delta}(1-e^{-\delta})}} \xrightarrow[n\to\infty]{d}Z$, where $Z$ is a standard Normal random variable. Therefore 
$$\B{E}[S_2(n)]\leq M\B{P}(S_2(n)>N+n)=M\B{P}\left(\frac{S_2(n)-n}{\sqrt{ne^{-\delta}(1-e^{-\delta})}} > N\right)\approx p(e^\delta-1)n\B{P}(Z>c\sqrt{n})\xrightarrow[n\to\infty]{}0,$$
and so as $n\to\infty$, $\B{E}[\alpha_n(r)-\beta_n(r)]\geq -\frac1n\B{E}[\mathfrak{B}]-\frac1n \B{E}[S_2(n)]\xrightarrow[n\to\infty]{}0$,
and therefore as $n$ was arbitrary
$$e^\delta \int \rho_0(x) \B{P}_x(B_\delta \leq r)dx\geq e^\delta \int \rho(x)\B{P}_x(B_\delta \leq r, \tau^{LR} > \delta)dx,$$
as required.
\end{proof}

Now we can observe that if $\int_{\B{R}}g(x)dx=\int_{\B{R}}h(x)dx =1$, then $\int_a^\infty g(x)dx\leq \int_a^\infty h(x)dx$ is equivalent to $\int_{-\infty}^a g(x)dx \geq \int_{-\infty}^a h(x)dx$. Also $\int_{-\infty}^a C^R_1g(x)dx =1\wedge \int_{-\infty}^a g(x)dx$ for any $a$. Therefore Proposition \ref{comparisonTheorem} implies (\ref{Sinequ}). Thus by induction and symmetry, the ordering (\ref{compSectionAim}) holds. We finally have all the tools ready to prove Theorem \ref{fbpExistence} and the hydrodynamic limit result of Theorem \ref{hydroLimitThm}.

\begin{proof}[(Of Theorem \ref{fbpExistence})]
By Proposition \ref{repThem}, there exists a solution $(u(x,t),L_t,R_t)$ of the FBP \eqref{fbpMain} and the boundaries $L$ and $R$ solves the two-sided inverse first passage problem. By Proposition \ref{comparisonTheorem}, the ordering \eqref{compSectionAim} holds for any solution $u$ of \eqref{fbpMain}. By Proposition \ref{detBoundDist}, these upper an lower bounds converge to the same limit, $u(x,t)$, as $n\to\infty$, and hence there can be at most one solution $u$ of \eqref{fbpMain}.
\end{proof}

\begin{proof}[(Of Theorem \ref{hydroLimitThm})] Note that by the union bound we have that for any $a\in \B{R},B,\beta>0$, we have
\begin{align}\nonumber \B{P}\Big(\Big| &\int_a^\infty \pi^N_{\rho,t}(dr) - \int_a^\infty u(dr,t)\Big|>BN^{-\beta}\Big) \\\label{theorem2Split}
&\leq \B{P}\Bigg(\int_a^\infty \pi^N_{\rho,t}(dr) - \int_a^\infty u(r,t)dr >BN^{-\beta}\Bigg) + \B{P}\Bigg(\int_a^\infty \pi^N_{\rho,t}(dr) - \int_a^\infty u(r,t)dr<-BN^{-\beta}\Bigg).\end{align}
Then since $\int_a^\infty \pi^N_{\rho,t}(dr)\leq \int_a^\infty \pi_{\rho,t}^{N,t/n,+}(dr)$ for any $n$ (by Proposition \ref{upDownCoupling}), thus using the union bound again
\begin{align*}
\B{P}\Bigg(\int_a^\infty \pi^N_{\rho,t}(dr) - \int_a^\infty u(r,t)dr>BN^{-\beta}\Bigg) &\leq \B{P}\Bigg(\int_a^\infty \pi^{N,t/n,+}_{\rho,t}(dr) - \int_a^\infty u(r,t)dr>BN^{-\beta}\Bigg) \\
\leq \B{P}\Bigg(\int_a^\infty \pi^{N,t/n,+}_{\rho,t}(dr) &- \int_a^\infty S_t^{t/n,+}\rho(dr) >\frac{BN^{-\beta}}{2}\Bigg)\\
&+ \B{P}\Bigg(\int_a^\infty S_t^{t/n,+}\rho(dr) - \int_a^\infty u(r,t)dr >\frac{BN^{-\beta}}{2}\Bigg) \\
\leq 16B^{-4}Ce^{4t}n^5e^{4t/n}N^{4\beta-2} &+ \is_{\{2e^t(e^{t/n}-1)>BN^{-\beta}/2\}},
\end{align*}
where for the first term second inequality we have used the bound from Proposition \ref{couplingConv} and for the second term we use the fact that by Proposition \ref{detBoundDist} and Proposition \ref{comparisonTheorem}: $$\Bigg|\int_x^\infty S_t^{t/n,+}\rho(dr) - \int_x^\infty u(r,t)dr\Bigg|<\Bigg|\int_x^\infty S_t^{t/n,+}\rho(dr) - \int_x^\infty S_t^{t/n,-}\rho(dr)\Bigg|<2e^t(e^{t/n}-1).$$
So then as $e^{t/n}-1< 2t/n$ for $n$ sufficiently large (say $n\geq n_1$), we choose $N$ to be sufficiently large so that $n:=8B^{-1}te^tN^{\beta}>n_1$, so that $\is_{\{2e^t(e^{t/n}-1)>BN^{-\beta}/2\}}=0$. Thus for some constant $C'$ independent of $a,B,\beta,t$ we have
\begin{align*}\B{P}\Bigg(\int_a^\infty \pi^N_{\rho,t}(dr) - \int_a^\infty u(r,t)dr>BN^{-\beta}\Bigg)\leq 16B^{-4}Ce^{4t}(8B^{-1}te^t N^\beta)^5\exp\left(\frac{4t}{8B^{-1}te^tN^\beta}\right)N^{4\beta-2} \\= B^{-9}C't^5e^{9t}e^{B/2N^{\beta}e^{t}}N^{9\beta-2}.\end{align*}
By symmetry, the same argument holds for the second term of (\ref{theorem2Split}), and summing these two upper bounds completes the proof.  
\end{proof}

\section{Convergence of the Leftmost and Rightmost Particles} \label{leftmostrightmost}
As well as proving the convergence of the empirical density of particles to the solution of the free boundary problem, we also want to prove convergence of the leftmost and rightmost particles to the free boundaries $L_t$ and $R_t$ respectively, almost surely and in expectation. This is stated in Theorem \ref{leftmostRightmostTheorem} and this result will be helpful to prove the convergence of the asymptotic speed in the next section. 

\vspace{1em}

\begin{propn} \label{asBdyConverge}
    Consider the solution $(u(x,t),L_t,R_t)_{t\geq 0,x\in \B{R}}$ of the FBP \eqref{fbpMain} with initial condition $\rho$. Let $(X^N(t))_{t\geq 0}$ be the $(N,p)$-BBM with each particle initially being i.i.d. with density $\rho$. Then for fixed $t$, $X_1^N(t)\xrightarrow[N\to\infty]{a.s.}L_t$ and $X_N^N(t)\xrightarrow[N\to\infty]{a.s.}R_t$. 
\end{propn}

\begin{proof}
    \sloppy To prove this theorem we will use the first Borel-Cantelli lemma. So fix $\epsilon>0$ and let $K$ be the local Lipschitz constant of $L$ for a small interval around $t$. Then define the event $F_t^N:=\{X_1^N(t)<L_t-2\epsilon\}$. We want to prove that $\sum_{N\geq 1}\B{P}(F_t^N)<\infty$. We will do this by considering the following events such that, when they occur simultaneously (which happens with high probability), we have $X_1^N(t)\geq L_t -2\epsilon$.
    
    Firstly, we want to have $\pi^N_{\rho,t-N^{-1/40}}((-\infty,L_{t-N^{-1/40}}))<N^{-1/10}$, which is to say that there are at most $N^{9/10}$ particles left of $L_{t-N^{-1/40}}$ at time $t-N^{-1/40}$. Call this event $E_t^{N,1}$. Choosing $\beta = 1/10$ in Theorem \ref{hydroLimitThm} shows that there exists a constant $c_1$ such that $\B{P}(E_t^{N,1})\geq 1-c_1N^{-11/10}$ for $N$ sufficiently large. 
    
    Let us call the particles which are left of $L_{t-N^{-1/40}}$ at time $t-N^{-1/40}$ red particles, and say that in a branching event, a red particle branches into two red particles. Next we want that over the time interval $[t-N^{-1/40},t]$, the red particles (of which there are initially at most $N^{9/10}$) have at most $N^{19/20}$ total descendants by time $t$. Call this event $E_t^{N,2}$. Note that the red particles branch at rate $1$, and we may delete the leftmost particle at branching events, therefore the number of descendants of a single particle in this system is bounded above by a Yule process. Thus the number of descendants of a particle after $N^{-1/40}$ time is bounded above by a $\text{Geom}(e^{-N^{-1/40}})$ random variable, therefore the probability of $E_t^{N,2}$ is bounded below by:
    $$\B{P}(E_t^{N,2})\geq \B{P}\left(\sum_{i=1}^{N^{9/10}}G_i\leq N^{19/20}\right),$$
    where $G_i$ are i.i.d. $\text{Geom}\big(e^{-N^{-1/40}}\big)$ random variables. Thus by Chebyshev's inequality, this happens with probability at least $1-\text{Var}\big(\sum_{i=1}^{N^{9/10}}G_i\big)/\big(N^{19/20}-\B{E}\big[\sum_{i=1}^{N^{9/10}}G_i\big]\big)^2\sim 1-N^{-41/40}$, therefore for some constant $c_2$, we have $\B{P}(E_t^{N,2})\geq 1-c_2N^{-41/40}$ for $N$ sufficiently large. 

    Let us call the remaining particles green, and say that at each branching event, a green particle branches into two green particles. So next we want that, over the time interval $[t-N^{-1/40},t]$, the green particles have at least $N^{19/20}$ branching events in which the leftmost particle is deleted. Call this event $E_t^{N,3}$. Each green particle branches at unit rate, and at each branching event, the leftmost particle (red or green) is deleted with probability $p$. Note that, given the event $E_t^{N,2}$, there are, for all times in $[t-N^{-1/40},t]$, at least $N-N^{19/20}$ green particles alive (since there are at most $N^{19/20}$ red particles). Therefore we can bound the probability of $E_t^{N,3}$ given $E_t^{N,2}$ below by:
    $\B{P}(E_t^{N,3}|E_t^{N,2})\geq \B{P}(H \geq N^{19/20})$,
    where $H$ is a $\text{Pois}(pN^{-1/40}(N-N^{19/20}))$ random variable. Then by the Chernoff bound $\B{P}(\text{Pois}(\lambda)\leq a)\leq e^{a-\lambda}\lambda^aa^{-a}=\exp(-\lambda +a(1+\ln(\lambda/a)))$, we get that 
    \begin{align*}
        \B{P}(E_t^{N,3}|E_t^{N,2})&\geq 1-\exp\left(-pN^{39/40}+pN^{37/40}+N^{19/20}\left(1+\ln(pN^{1/40}-pN^{-1/40})\right)\right)\\
        &\sim 1-\exp(-pN^{39/40}),
    \end{align*}
    therefore for some constant $c_3$ we have $$\B{P}(E_t^{N,2}\cap E_t^{N,3})=\B{P}(E_t^{N,2})\B{P}(E_t^{N,3}|E_t^{N,2})>1-c_2N^{-41/40}-c_3\exp(-pN^{39/40}).$$
    
    Finally, consider the Brownian motions driving the $N$ particles. Consider that the $(N,p)$-BBM process can be embedded inside a branching Brownian motion with $N$ particles at time $t-N^{-1/40}$. Say the $i$\textsuperscript{th} BBM has particles $\C{B}_u$ for $u\in\C{N}^i(t)$ at time $t$. We want that the maximum displacement of these Brownian motions over the time interval $[t-N^{-1/40},t]$ to be at most $\epsilon/3$; that is, $$\max_{i=1}^N \max_{u\in\C{N}^i(t)}\sup_{s\in [t-N^{-1/40},t]}|B_u(s)-B_u(t-N^{-1/40})|\leq \epsilon/3.$$
    Call this event $E_t^{N,4}$, and note that we have by the union bound, the many-to-one lemma, and the reflection principle
    \begin{align*}\B{P}(E_t^{N,4})&\geq 1-N\B{P}\left(\sup_{s\in [t-N^{-1/40},t]}|B_u(s)-B_u(t-N^{-1/40})|>\epsilon/3\text{ for some }u\in \C{N}^i(t)\right) \\
    &\geq 1-Ne^{N^{-1/40}}\B{P}\left(\sup_{s\in[t-N^{-1/40},t]}|B_u(s)-B_u(t-N^{-1/40})|> \epsilon/3\right)\\
    &= 1- 4Ne^{N^{-1/40}}\B{P}(B(N^{-1/40})>\epsilon/3) \\
    &= 1-12Ne^{N^{-1/40}}\exp(-\epsilon^2 N^{1/40}/18)/(\epsilon N^{1/80}\sqrt{2\pi}),\end{align*}
    where $B$ is a standard Brownian motion. 
        
    We will now show that on the event $E_t^{N,1}\cap E_t^{N,2}\cap E_t^{N,3}\cap E_t^{N,4}$, all particles are to the right of $L_{t-N^{-1/40}}-\epsilon$ at time $t$. On the events $E^{N,1}_t$ and $E^{N,2}_t$, there are at most $N^{19/20}$ red particles at time $t$. On $E^{N,4}_t$, all green particles never go to the left of $L_{t-N^{-1/40}}-\epsilon/3$ during the time interval $[t-N^{-1/40},t]$. Under $E^{N,3}_t$, during the time interval $[t-N^{-1/40},t]$, there are at least $N^{19/20}$ events in which a green particle branches and the leftmost particle is simultaneously deleted. Suppose that at time $t$ there is at least one red particle remaining which has not been deleted due to one of the $N^{19/20}$ branching events of a green particle in which the leftmost particle is deleted. This means it must not have been the leftmost particle at any of these branching times. Since green particles never go to the left of $t-N^{-1/40}-\epsilon/3$ during $[t-N^{-1/40},t]$, this red particle must have been, at some point in $[t-N^{-1/40},t]$, to the right of $L_{t-N^{-1/40}}-\epsilon/3$. Therefore, due to the event $E^{N,4}_t$, this red particle must always be right of $L_{t-N^{-1/40}}-\epsilon$ during $[t-N^{-1/40},t]$. Therefore by time $t$, any red particle which has not been deleted is to the right of $L_{t-N^{-1/40}}-\epsilon$ and all green particles are to the right of $L_{t-N^{-1/40}}-\epsilon/3$. Therefore we can conclude that at time $t$, all particles are to the right of $L_{t-N^{-1/40}}-\epsilon$.
    
    Then as $L$ is by assumption Lipschitz continuous in a small interval around $t$, with constant $K$, thus $|L_{t-N^{-1/40}}-L_t|<KN^{-1/40}$. For $N$ sufficiently large, $KN^{-1/40}<\epsilon$, so then $X_1^N(t)>L_t-2\epsilon$ for $N$ sufficiently large. Therefore by the union bound and 
    \begin{align*}\B{P}(X_1^N(t)<L_t&-\epsilon-KN^{-1/40})=\B{P}(F_t^N)\leq 1-\B{P}(E_t^{N,1}\cap E_t^{N,2}\cap E_t^{N,3}\cap E_t^{N,4})\\
    &= \B{P}((E_t^{N,1})^c\cup (E_t^{N,2}\cap E_t^{N,3})^c \cup (E_t^{N,4})^c)\\
    &\leq c_1N^{-11/10}+c_2N^{-41/40}+c_3\exp\left(-pN^{39/40}\right)+c_4N^{79/80}\exp\left(-\epsilon^2 N^{1/40}/18\right)\\
    &\leq \tilde{c}N^{-41/40}\end{align*}
    for some constant $\tilde{c}$. Therefore $\sum_{N=1}^\infty\B{P}(F_t^N)<\infty$, and therefore by the first Borel-Cantelli lemma, $X_1^N(t)<L_t-2\epsilon$ occurs only finitely often with probability $1$.

    Next we will prove that $X_1^N(t)>L_t+\epsilon$ finitely often with probability $1$. Since $u$ is continuous and positive on $(L_t,R_t)$, therefore $\ell:=\int_{L_t+\epsilon}^\infty u(y,t)dy < 1$, and so applying Theorem \ref{hydroLimitThm} with $B=1-\ell$ and $\beta=0$, for some constant $C$ and for $N$ sufficiently large:
    $$\B{P}(X_1^N(t)>L_t+\epsilon)\leq \B{P}\left(\left|\int_{L_t+\epsilon}^\infty \pi^N_t\rho(y)dy - \int_{L_t+\epsilon}^\infty u(y,t)dy\right|<1-\ell\right)\leq C(1-\ell)^{-9}N^{-2},$$
    therefore since $\sum_{N=1}^\infty C(1-\ell)^{-9}N^{-2}<\infty$, thus by the Borel-Cantelli lemma, the event $\{X_1^N(t)>L_t+\epsilon\}$ occurs only finitely often with probability $1$.

    Combining the upper and lower bounds, we have that $\{|X_1^N(t)-L_t|>2\epsilon\}$ occurs finitely often with probability $1$. Therefore since $\epsilon$ was arbitrary, we have that $X_1^N(t)\xrightarrow[N\to\infty]{a.s.}L_t$. By symmetry the result holds for $X_N^N(t)\xrightarrow[N\to\infty]{a.s.}R_t$.
\end{proof}

\vspace{1em}

\begin{propn} \label{boundariesConvergeE}
    Consider the solution $(u(x,t),L_t,R_t)_{t\geq 0,x\in \B{R}}$ of the FBP \eqref{fbpMain} with initial condition $\rho$. Let $(X^N(t))_{t\geq 0}$ be the $(N,p)$-BBM with each particle initially being i.i.d. with density $\rho$. Then for fixed $t$, $\B{E}[X_1^N(t)]\xrightarrow[N\to\infty]{}L_t$ and $\B{E}[X_N^N(t)]\xrightarrow[N\to\infty]{}R_t$.
\end{propn}

\begin{proof}
    For ease let $Z=Z^N_t:=|X_1^N(t)-L_t|$, and fix some $t\in (0,\infty)$ and $\delta>0$. We can write:
    \begin{equation}\label{expConvBoundZ}\B{E}[Z]=\B{E}\Big[Z\big(\is_{\{Z\leq \delta\}}+\is_{\{\delta < Z \leq (\log N)^2\}} + \is_{\{(\log N)^2 < Z\}}\big)\Big]\leq \delta + (\log N)^2 \B{P}(Z>\delta) + \B{E}[Z\is_{\{Z>(\log N)^2\}}].\end{equation}

    First let us recall that from the proof of Proposition \ref{asBdyConverge} that, for fixed $\epsilon>0$ and fixed $t$, we have that for $N$ sufficiently large $\B{P}(|X_1^N(t)- L_t|>2\epsilon)\leq c N^{-41/40}$ for some constant $c$
    therefore if we choose $\epsilon = \delta/2$ and $N$ such that $KN^{-1/40}<\epsilon$, then it follows that $(\log N)^2\B{P}(Z>\delta)\to 0$ as $N\to\infty$. 

    To show that $\B{E}[Z\is_{\{Z>(\log N)^2\}}]\xrightarrow{N\to\infty}0$, note that the $(N,p)$-BBM with initial density $\rho$ is embedded inside a branching Brownian motion (BBM) started from $N$ particles each with initial density $\rho$. Let the $i$\textsuperscript{th} particle of the BBM at time $0$ have, at time $t$, offspring $u\in \C{N}_i(t)$, with locations $B_u^i(t)$ for $u\in \C{N}_i(t)$. Therefore
    \begin{equation}\label{expConvergeBound}Z \leq \max_{i=1}^N \max_{u\in \C{N}_i(t)} |B_u^i(t)| + |L_t|.\end{equation}
    
    Then we observe that
    \begin{align}\label{term3ExpConvergeBound}Z\is_{\{Z>(\log N)^2\}}\leq \max_{i=1}^N \max_{i\in \C{N}_i(t)}(|B_u^i(t)|+|L_t|)\is_{\big\{\max_{i=1}^N\max_{i\in \C{N}_i(t)}|B_u^i(t)|>(\log N)^2-|L_t|\big\}}\\ \leq \sum_{i=1}^N \sum_{u\in \C{N}_i(t)}(|B_u^i(t)|+|L_t|)\is_{\{|B_u^i(t)|>(\log N)^2-|L_t|\}},\end{align}
    so that combining (\ref{expConvergeBound}) and (\ref{term3ExpConvergeBound}) and the many-to-one lemma, we get that:
    $$\B{E}[Z\is_{\{Z>(\log N)^2\}}]\leq Ne^t \B{E}[(|B(t)|+L_t)\is_{\{|B(t)|>(\log N)^2-L_t\}}],$$
    where $B$ is a standard Brownian motion started from an initial location $B(0)\sim \rho$, which has compact support, so by the tail behaviour of the normal distribution, this behaves asymptotically like a constant times $N\exp(-(\log N)^4)$, which converges to $0$ as $N\to\infty$. Thus as $\delta$ was arbitrary, then by equation (\ref{expConvBoundZ}), we have that $\B{E}[Z]\to 0$ as $N\to\infty$, and therefore $\lim_{N\to\infty}\B{E}[X_1^N(t)]\to L_t$, as required. By symmetry, $\lim_{N\to\infty}\B{E}[X_N^N(t)]\to R_t$. 
\end{proof}

\section{Convergence of Speeds} \label{speed}
In this section, we will prove Theorem \ref{speedThm}; that the asymptotic velocity $v_N$ of the $(N,p)$-BBM process converges to the travelling wave speed $c=\sqrt{\frac{2\log^2(p/(1-p))}{\log^2(p/(1-p))+\pi^2}}$ as $N\to\infty$.

To prove this result, we will use the fact that the $(N,p)$-BBM is monotonic, and, as viewed from its leftmost (resp. rightmost) particle is positive Harris recurrent (Propositions \ref{monotonicity} and \ref{harrisRecc}), the fact that the FBP (\ref{fbpMain}) has a unique travelling wave with wave speed $c$ (Proposition \ref{travellingWave}), and the fact that, for suitable initial conditions, the leftmost and rightmost particles of the $(N,p)$-BBM converge in expectation to the left and right free boundaries of the free boundary problem (Proposition \ref{boundariesConvergeE}). First we prove a technical lemma.

\vspace{1em}

\begin{lemma} \label{birkhoffLemma}
    Let $(\bar{X}^{N,p}(t))_{t\geq 0}:=(X^{N,p}(t)-X_1^{N,p}(t)\underbar{1})_{t\geq 0}$ be the $(N,p)$-BBM process $X^{N,p}$ as viewed from its leftmost particle, and let $\bar{\pi}^N$ denote its stationary distribution. Then for any $t\geq 0$, $\B{E}[X_1^{N,p}(t)|X^{N,p}(0)\sim \bar{\pi}^N]=v_Nt$, where $v_N:=\lim_{t\to\infty}X_1^{N,p}(t)/t$.
\end{lemma}

\begin{proof}
    Fix $T\geq 0$ and consider the $(N,p)$-BBM with initial distribution $X^{N,p}(0)\sim \bar{\pi}^N$. Now consider the random variable $X^{N,p}_1(t)-X^{N,p}_1(0)$ and say it has distribution $\xi^{N,p}$. Since the process $(X^{N,p}(t))_{t\geq 0}$ starts from the stationary distribution $\bar{\pi}^N$, thus the mapping which maps $X^{N,p}_1(it)-X^{N,p}_1(it-t)$ to $X^{N,p}_1(it+t)-X^{N,p}_1(it)$ is a $\xi^{N,p}$-preserving map. Therefore we can apply Birkhoff's ergodic Theorem (see Theorem 25.6, \cite{kallenberg}, for example) to show that:
    \begin{align*}\B{E}[X_1^{N,p}(t)|X^{N,p}(0)\sim \bar{\pi}^N]&=\B{E}[X_1^{N,p}(t)-X_1^{N,p}(0)|X^{N,p}(0)\sim \bar{\pi}^N]\\
    &=\lim_{n\to\infty}\frac1n\sum_{i=0}^{n-1}(X_1^{N,p}(it+t)-X_1^{N,p}(it)) = \lim_{n\to\infty}\frac1n X_1^{N,p}(nt)=v_Nt\end{align*}
\end{proof}

\begin{proof} (Of Theorem \ref{speedThm})
    Consider the travelling wave solution of the free boundary problem (\ref{fbpMain}), shifted so that the right barrier is at $0$ at time $0$. Proposition \ref{travellingWave} tells us that this travelling wave solution has a profile given by:
    $$\hat{u}(x):=\frac{2p}{\sqrt{2-c_p^2}}e^{-c_px}\sin\left(x\sqrt{2-c_p^2}\right)\is_{\{-R_0 \leq x\leq 0\}}$$
    where $c_p$ is the travelling wave speed $c_p=\sqrt{\frac{2\log^2(p/(1-p))}{\log^2(p/(1-p))+\pi^2}}$ for $p\geq 1/2$ and $c_p=-\sqrt{\frac{2\log^2(p/(1-p))}{\log^2(p/(1-p))+\pi^2}}$ for $p\leq 1/2$, and $R_0=\frac{\pi}{\sqrt{2-c_p^2}}$. 
    Then by monotonicity (Proposition \ref{monotonicity}), we have that $\B{E}_{\hat{u}}[X_1^{N,p}(t)]\leq \B{E}[X_1^{N,p}(t)|X^{N,p}(0)\sim \bar{\pi}^N]$ for any $t$ and $N$. Notice that the expectation $\B{E}_{\hat{u}}[X_1^{N,p}(t)]$ is the expectation conditional on each particle having i.i.d. initial density $\hat{u}$, whereas $\B{E}[X_1^{N,p}(t)|X^{N,p}(0)\sim \bar{\pi}^N]$ is the expectation conditional on the initial configuration being drawn from the probability distribution $\bar{\pi}^N$ on $\B{R}^N$. 

    Then by Proposition \ref{boundariesConvergeE}, $\B{E}_{\hat{u}}[X_1^{N,p}(t)]\to \bar{L}_t$ as $N\to\infty$, where $\bar{L}_t=c_pt-R_0$ is the position of the left free boundary in the FBP (\ref{fbpMain}) with initial condition $\hat{u}$. So fix $\epsilon>0$ and let $t \geq 1 \vee (2R_0/\epsilon)$. Moreover, let $N$ be sufficiently large so that $|\B{E}_{\hat{u}}[X_1^{N,p}(t)]-\bar{L}_t|<\epsilon/2$. So then $$v_N =\B{E}[X_1^{N,p}(t)|X^{N,p}(0)\sim \bar{\pi}^N]/t \geq \B{E}_{\hat{u}}[X_1^{N,p}(t)]/t \geq c_p-R_0/t-\epsilon/2 \geq c_p-\epsilon.$$

    Similarly, we can consider the stationary distribution of $X^{N,p}$ as viewed from its rightmost particle, which we will call $\hat{\pi}^N$, and the travelling wave solution of (\ref{fbpMain}), as given in Proposition \ref{travellingWave}, which we will call $\bar{u}$. Consider $N$ such that $|\B{E}_{\bar{u}}[X_N^{N,p}]-R_t|<\epsilon/2$, where $R_t=c_pt$ is the right-hand boundary of the travelling wave solution. Then by the same reasoning as before, for $t\geq 1\vee R_0/2\epsilon$, $$v_N=\B{E}[X_N^{N,p}(t)|X^{N,p}(0)\sim \hat{\pi}^N]/t \leq \B{E}_{\bar{u}}[X_N^{N,p}(t)]/t\leq c_p+R_0/t +\epsilon/2 \leq c_p+\epsilon.$$
    Therefore for $N$ sufficiently large, $|v_N - c_p|<\epsilon$, and hence as $\epsilon$ was arbitrary $\lim_{N\to\infty}v_N = c_p$. 
\end{proof}

\section{Appendix}\label{appendix}\begin{proof}
    \textit{(Of Lemma \ref{llnProbBound})}
    By Markov's inequality
    \begin{align*}\B{P}\left(\frac1N\left|\sum_{i=1}^N A_i - \sum_{i=1}^N \mu_i \right|>BN^{-\beta}\right) &= \B{P}\left(\left| \frac1N \sum_{i=1}^N (A_i - \mu_i)\right|^4 > B^4 N^{-4\beta}\right) \\ 
    &\leq B^{-4} N^{4\beta}\frac{1}{N^4}\B{E}\left[\left(\sum_{i=1}^N (A_i - \mu_i)\right)^4\right].\end{align*}
    Then we note that the only non-zero terms of the expansion of $\left(\sum_{i=1}^N (A_i - \mu_i)\right)^4$ are $$\sum_{i=1}^N \B{E}[(A_i - \mu_i)^4] + 6\sum_{i,j\in\{1,2,\ldots,N\}, i\neq j}\B{E}[(A_i-\mu_i)^2]\B{E}[(A_j - \mu_j)^2]$$
    therefore 
    \begin{align*}B^{-4}N^{4\beta}\frac{1}{N^4}\B{E}\left[\left(\sum_{i=1}^N (A_i - \mu)\right)^4\right]&=B^{-4}N^{4\beta - 2}\left(\frac{1}{N^2}\sum_{i=1}^N \B{E}[(A_i - \mu_i)^4]+\frac{6}{N^2}\sum_{i\neq j}\B{E}[(A_i - \mu_i)^2]\B{E}[(A_j-\mu_j)^2]\right)\\
    &\leq B^{-4}CkN^{4\beta-2}\end{align*}
    for sufficiently large $N$.
\end{proof}

\begin{proof}
    \textit{(Of Corollary \ref{probBoundAddition})}
    By the triangle inequality
    \begin{align*}
        \B{P}(|(A_N+A'_N)&-(a+a')|>BN^{-\beta})\leq \B{P}(|A_N-a|+|A'_N-a'|>BN^{-\beta}) \\
        &\leq \B{P}(|A_N-a|>BN^{-\beta}/2)+\B{P}(|A'_N-a'|>BN^{-\beta}/2) \\
        &\leq 16B^{-4}C_AN^{4\beta-2}+16B^{-4}C'_AN^{4\beta-2} 
    \end{align*}
\end{proof}

\begin{proof}
    \textit{(Of Lemma \ref{compoundPoissonLemma})}
    For $i=1,2,\ldots,n$, the random variable $Y_i:=\sum_{j=1}^{M^i_N}A_{i,j}$ is called a \textit{compound Poisson} random variable, and if we denote the characteristic function of $A_{i,1}$ by $\phi_{A_i}(t)$, then the characteristic function of $Y_i$ is $\phi_{Y_i}(t):=e^{\lambda N (\phi_{A_i}(t)-1)}$. Therefore we can calculate the moments of $Y_i$ by differentiating $\phi_{Y_i}(t)$. We temporarily omit the subscript $i$ for readability. Differentiating four times gives:
    \begin{align*}
        \phi'_Y(t)&=\lambda N\phi'_X(t)\phi_Y(t) \\
        \phi''_Y(t)&=\lambda N\phi''_X(t)\phi_Y(t) + (\lambda N\phi'_X(t))^2 \phi_Y(t) \\
        \phi^{(3)}_Y(t)&= \lambda N\phi^{(3)}_X(t)\phi_Y(t) + 3\lambda^2 N^2 \phi'_X(t)\phi''_X(t)\phi_Y(t) + (\lambda N \phi'_X(t))^3\phi_Y(t) \\
        \phi^{(4)}_Y(t)&=\lambda N \phi^{(4)}_X(t)\phi_Y(t) + 4\lambda^2 N^2\phi^{(3)}_X(t) \phi'_X(t) \phi_Y(t)+ 3\lambda^2 N^2 \phi''_X(t)^2\phi_Y(t) +  \\
        &+ 6\lambda^3 N^3 \phi'_X(t)^2 \phi''_X(t)\phi_Y(t) + (\lambda N\phi'_X(t))^4\phi_Y(t).
    \end{align*}
    Then evaluating at $t=0$, if we write $m_1,m_2,m_3,m_4$ for the first $4$ moments of $A$ respectively, then we can calculate
    \begin{align*}
        \B{E}[Y]&=\lambda N m_1 \\
        \B{E}[Y^2]&=\lambda N m_2 + \lambda^2 N^2 m_1^2 \\
        \B{E}[Y^3]&=\lambda N m_3 + 3\lambda^2 N^2 m_1 m_2 + \lambda^3 N^3 m_1^3\\
        \B{E}[Y^4]&=\lambda N m_4 + 4\lambda^2 N^2 m_1 m_3 + 3\lambda^2 N^2 m_2^2 + 6\lambda^3 N^3 m_1^2 m_2 + \lambda^4 N^4 m_1^4.
    \end{align*}
    Therefore expanding $\B{E}[(Y-N\mu)^4]$ and evaluating we yield 
    \begin{align*}
        \B{E}[(Y-\lambda N m_1)^2]=\B{E}[Y^2 - 2\lambda N m_1 Y + \lambda^2 N^2 m_1^2]=\lambda N m_2,
    \end{align*}
    and
    \begin{align*}
        \B{E}[(Y-\lambda Nm_1)^4]=\B{E}[Y^4 &-4\lambda Nm_1 Y^3 + 6\lambda^2 N^2 m_1^2 Y^2 -4 \lambda^3 N^3 m_1^3 Y + \lambda^4 N^4 m_1^4] \\
        =\lambda N m_4& + 4\lambda^2 N^2 m_1 m_3 + 3\lambda^2 N^2 m_2^2 + 6\lambda^3 N^3 m_1^2 m_2 + \lambda^4 N^4 m_1^4 \\ &-4\lambda^2N^2 m_1 m_3 -12 \lambda^3 N^3 m_1^2 m_2 -4\lambda^4 N^4 m_1^4 + 6\lambda^3 N^3 m_1^2 m_2 + 6\lambda^4 N^4 m_1^4 \\
        &- 4\lambda^4 N^4 m_1^4 + \lambda^4 N^4 m_1^4 \\
        =\lambda N m_4 & + 3\lambda^2 N^2 m_2^2.
    \end{align*}
    Then noting that $\B{E}[Y_i-\lambda N \mu_i]=0$ and using the expressions for $\B{E}[(Y_i - \lambda N \mu_i)^2]$ and $\B{E}[(Y_i - \lambda N \mu_i)^4]$ that we calculated above, we can expand:
    \begin{align*}
        \B{E}\left[\left(\sum_{i=1}^n (Y_i - N\lambda \mu_i)\right)^4\right] = \sum_{i=1}^n \B{E}[(Y_i - N\lambda \mu_i)^4] + 6\sum_{i,j=1, i\neq j}^n \B{E}[(Y_i - N\lambda \mu_i)^2]\B{E}[(Y_j - N\lambda \mu_j)^2] \\
        = \sum_{i=1}^n (N\lambda \B{E}[A_{i,1}^4] + 3N^2 \lambda^2 \B{E}[A_{i,1}^2]^2) + 6\sum_{i,j=1,i\neq j}^n N^2 \lambda^2 \B{E}[A_{i,1}^2]\B{E}[A_{j,1}^2]
    \end{align*}
    Then by Markov's inequality: 
    \begin{align*}
        \B{P}\left(\left|\sum_{i=1}^n\frac1N \sum_{i=1}^{M_N} A_i - \lambda \sum_{i=1}^n \mu_i\right|>BN^{-\beta}\right) = \B{P}\left(\left|\sum_{i=1}^n(Y_i-N\lambda \mu_i)\right|^4>B^4N^{4-4\beta}\right)\\
        \leq B^{-4}C\left(\sum_{i=1}^n \lambda^2 \B{E}[A_{i,1}^2]^2 + 6\sum_{i,j=1,i\neq j}^n \lambda^2 \B{E}[A_{i,1}^2]\B{E}[A_{j,1}^2]\right)N^{4\beta-2}
    \end{align*}
    for some constant $C$ and sufficiently large $N$.
\end{proof}

\begin{proof}\textit{(Of Lemma \ref{functionLims})} For all $x$, $(\int_x^\infty f_n(s)ds)_{n\in \B{N}}$ is monotonic increasing and bounded above by, therefore has a limit $\lim_{n\to\infty}\int_x^\infty f_n(s)ds:=\Psi(x)\leq \int_x^\infty b(s)ds$. In remains to prove that the convergence is uniform in $x$. First note that by sandwiching $\int_x^\infty f_1(s)ds \leq \Psi(x)\leq \int_x^\infty b(s)ds$, we have that $\Psi(x)\to 1$ and $x\to-\infty$ and $\Psi(x)\to 0$ as $x\to \infty$. Therefore fix $\epsilon>0$ and choose $M$ such that $x\geq M \implies \int_x^\infty b(s)ds < \epsilon$ and $x\leq - M \implies \int_x^\infty f_1(s)ds > 1-\epsilon$. Therefore, since $\int_x^\infty f_1(s)ds \leq \int_x^\infty f_n(s)ds \leq \Psi(x)\leq \int_x^\infty b(s)ds$, thus $|x|>M \implies \left|\int_x^\infty f_n(s) - \Psi(x)\right|\leq \epsilon$ for all $n\in \B{N}$. Therefore convergence is uniform outside of $[-M,M]$. By compactness, convergence is uniform inside $[-M,M]$. Putting this together gives the desired result. 
\end{proof}

\section*{Acknowledgements}
The author thanks Julien Berestycki for his careful supervision of the project, as well as helpful conversations with Sarah Penington and Oliver Tough and feedback from Matthias Winkel and Christina Goldschmidt on an earlier version of the work.

\bibliographystyle{plain}
\bibliography{bps}

@article{globalSolsNBBM,
author = {Berestycki, Julien and Brunet, Éric and Penington, Sarah},
copyright = {2019 IOP Publishing Ltd & London Mathematical Society},
issn = {0951-7715},
journal = {{Nonlinearity}},
language = {eng},
number = {10},
pages = {3912-3939},
publisher = {IOP Publishing},
title = {{Global existence for a free boundary problem of Fisher-KPP type}},
volume = {32},
year = {2019},
}

@article{lalleySellke,
issn = {0091-1798},
journal = {The Annals of probability},
pages = {1052--1061},
volume = {15},
publisher = {Institute of Mathematical Statistics},
number = {3},
year = {1987},
title = {{A Conditional Limit Theorem for the Frontier of a Branching Brownian Motion}},
copyright = {Copyright 1987 Institute of Mathematical Statistics},
language = {eng},
address = {Hayward, CA},
author = {Lalley, S. P. and Sellke, T.},
}

@article{beesBBNP,
issn = {0091-1798},
journal = {The Annals of probability},
pages = {2133--2177},
volume = {50},
publisher = {INST MATHEMATICAL STATISTICS-IMS},
number = {6},
year = {2022},
title = {BROWNIAN BEES IN THE INFINITE SWARM LIMIT},
language = {eng},
address = {CLEVELAND},
author = {Berestycki, Julien and Brunet, Eric and Nolen, James and Penington, Sarah},
}

@incollection{hydroNBBM,
author = {De Masi, Anna and Ferrari, Pablo A. and Pesutti, Errico and Soprano-Loto, Nahuel},
address = {Switzerland},
booktitle = {Stochastic Dynamics Out of Equilibrium},
copyright = {Springer Nature Switzerland AG 2019},
isbn = {9783030150952},
issn = {2194-1009},
language = {eng},
pages = {523-549},
publisher = {Springer International Publishing AG},
series = {Springer Proceedings in Mathematics \& Statistics},
title = {{Hydrodynamics of the N-BBM Process}},
volume = {282},
year = {2019},
}

@article{afree,
issn = {0002-9947},
journal = {Transactions of the American Mathematical Society},
pages = {6269--6329},
volume = {374},
publisher = {Amer Mathematical Soc},
number = {9},
year = {2021},
title = {A free boundary problem arising from branching Brownian motion with selection},
copyright = {Copyright 2021, American Mathematical Society},
language = {eng},
address = {PROVIDENCE},
author = {Berestycki, Julien and Brunet, Eric and Nolen, James and Penington, Sarah},
}

@book{todorovic,
author = {Todorovic, P.},
address = {New York},
booktitle = {An introduction to stochastic processes and their applications},
isbn = {038797783X},
language = {eng},
lccn = {91046692},
publisher = {Springer-Verlag},
series = {Springer series in statistics. Probability and its applications},
title = {An introduction to stochastic processes and their applications },
year = {1992},
}

@article{athreyaNey,
 ISSN = {00029947},
 URL = {http://www.jstor.org/stable/1998882},
 author = {K. B. Athreya and P. Ney},
 journal = {Transactions of the American Mathematical Society},
 pages = {493--501},
 publisher = {American Mathematical Society},
 title = {{A New Approach to the Limit Theory of Recurrent Markov Chains}},
 urldate = {2024-01-18},
 volume = {245},
 year = {1978}
}

@article{sipKingman,
    author = {Kingman, J.F.C},
    title = {{Subadditive Ergodic Theory}},
    journal = {{The Annals of Probability}},
    year = {1973},
    volume = {1},
    number = {6},
    pages = {883-909}
}

@article{kim,
language = {eng},
number = {2},
pages = {1315-},
publisher = {Institute of Mathematical Statistics},
title = {{The maximum of branching Brownian motion in $\B{R}^d$}},
volume = {33},
year = {2023},
author = {Kim, Yujin H. and Lubetzky, Eyal and Zeitouni, Ofer},
address = {Hayward},
copyright = {Copyright Institute of Mathematical Statistics Apr 2023},
issn = {1050-5164},
journal = {The Annals of applied probability},
}

@book{carinci,
author = {Carinci, Gioia and De Masi, Anna and Giardinà, Cristian and Presutti, Errico},
address = {Switzerland},
booktitle = {Free boundary problems in PDEs and particle systems},
isbn = {9783319333700},
language = {eng},
publisher = {Springer},
series = {SpringerBriefs in mathematical physics, volume 12},
title = {Free boundary problems in PDEs and particle systems },
year = {2016},
}

@article{strip,
author = {Harris, S. C. and Hesse, M. and Kyprianou, A. E.},
address = {Hayward},
copyright = {Copyright © 2016 Institute of Mathematical Statistics},
issn = {0091-1798},
journal = {The Annals of probability},
language = {eng},
number = {1},
pages = {235-275},
publisher = {Institute of Mathematical Statistics},
title = {BRANCHING BROWNIAN MOTION IN A STRIP: SURVIVAL NEAR CRITICALITY},
volume = {44},
year = {2016},
}

@article{bbs1,
author = {Berestycki, Julien and Berestycki, Nathanaël and Schweinsberg, Jason},
address = {Hayward},
copyright = {Copyright © 2013 Institute of Mathematical Statistics},
issn = {0091-1798},
journal = {The Annals of Probability},
language = {eng},
number = {2},
pages = {527-618},
publisher = {Institute of Mathematical Statistics},
title = {{The Genealogy of branching Brownian motion with absorption}},
volume = {41},
year = {2013},
}

@article{guoHu06,
 ISSN = {0033569X, 15524485},
 URL = {http://www.jstor.org/stable/43638734},
 author = {Jong-Shenq Guo and Bei Hu},
 journal = {Quarterly of Applied Mathematics},
 number = {3},
 pages = {413--431},
 publisher = {Brown University},
 title = {ON A TWO-POINT FREE BOUNDARY PROBLEM},
 urldate = {2025-03-25},
 volume = {64},
 year = {2006}
}

@article{guoHuKoh03,
    author = {Chang, Y.-L. and Guo, J.-S. and Kohsaka, Y.},
    journal = {Asymptotic Analysis},
    volume = {34},
    year = {2003},
    title = {On a two-point free boundary problem for a quasilinear parabolic equation},
    pages = {333--358}
}

@article{predPreyFBP,
author = {Wang, Mingxin and Zhao, Jingfu},
address = {New York},
copyright = {Springer Science+Business Media New York 2015},
issn = {1040-7294},
journal = {Journal of dynamics and differential equations},
number = {3},
pages = {957-979},
publisher = {Springer US},
title = {A Free Boundary Problem for the Predator–Prey Model with Double Free Boundaries},
volume = {29},
year = {2017},
}

@article{downey,
author = {Downey, Peter J.},
address = {Amsterdam},
copyright = {1990},
issn = {0167-6377},
journal = {Operations research letters},
language = {eng},
number = {3},
pages = {189-201},
publisher = {Elsevier B.V},
title = {Distribution-free bounds on the expectation of the maximum with scheduling applications},
volume = {9},
year = {1990},
}

@article {nonLocalNBBM,
    AUTHOR = {De Masi, Anna and Ferrari, Pablo A. and Presutti, Errico and
              Soprano-Loto, Nahuel},
     TITLE = {Non local branching {B}rownian motions with annihilation and
              free boundary problems},
   JOURNAL = {Electron. J. Probab.},
  FJOURNAL = {Electronic Journal of Probability},
    VOLUME = {24},
      YEAR = {2019},
     PAGES = {Paper No. 63, 30},
      ISSN = {1083-6489},
   MRCLASS = {60K35 (82C20)},
  MRNUMBER = {3978213},
MRREVIEWER = {Aernout\ C. D. van Enter},
       DOI = {10.1214/19-EJP324},
       URL = {https://doi.org/10.1214/19-EJP324},
}

@article {durrettRemenik,
    AUTHOR = {Durrett, Rick and Remenik, Daniel},
     TITLE = {Brunet-{D}errida particle systems, free boundary problems and
              {W}iener-{H}opf equations},
   JOURNAL = {Ann. Probab.},
  FJOURNAL = {The Annals of Probability},
    VOLUME = {39},
      YEAR = {2011},
    NUMBER = {6},
     PAGES = {2043--2078},
      ISSN = {0091-1798,2168-894X},
   MRCLASS = {60J80 (35K91 35R35 60F99)},
  MRNUMBER = {2932664},
MRREVIEWER = {Nicolas\ Fournier},
       DOI = {10.1214/10-AOP601},
       URL = {https://doi.org/10.1214/10-AOP601},
}

@misc{beresTough,
      title={{Selection principle for the $N$-BBM}}, 
      author={Berestycki, Julien and Tough, Oliver},
      year={2024},
      eprint={2407.05792},
      archivePrefix={arXiv},
      primaryClass={math.PR},
note={\textit{arXiv:2407.05792}},
      howpublished={\url{https://arxiv.org/abs/2407.05792}}, 
}

@book{yongZhou,
author = {Yong, Jiongmin and Zhou, Xun Yu},
address = {New York, NY},
booktitle = {Stochastic Controls : Hamiltonian Systems and HJB Equations},
edition = {1st ed. 1999.},
isbn = {1-4612-1466-1},
language = {eng},
publisher = {Springer New York},
series = {Stochastic Modelling and Applied Probability, 43},
title = {Stochastic Controls : Hamiltonian Systems and HJB Equations },
year = {1999},
}

@article{spreadVanDichot,
author = {Du, Yihong and Lin, Zhigui},
address = {Philadelphia},
copyright = {Copyright Society for Industrial and Applied Mathematics 2010},
issn = {0036-1410},
journal = {SIAM journal on mathematical analysis},
number = {1},
pages = {377-405},
publisher = {Society for Industrial and Applied Mathematics},
title = {Spreading-Vanishing Dichotomy in the Diffusive Logistic Model with a Free Boundary},
volume = {42},
year = {2010},
}

@article{shiftCutoff,
author = {Brunet, Eric and Derrida, Bernard},
issn = {1063-651X},
journal = {Physical review. E, Statistical physics, plasmas, fluids, and related interdisciplinary topics},
language = {eng},
number = {3},
pages = {2597-2604},
title = {Shift in the velocity of a front due to a cutoff},
volume = {56},
year = {1997},
}

@article{microModelsCutoff,
author = {Brunet, Eric and Derrida, Bernard},
copyright = {1999},
issn = {0010-4655},
journal = {Computer physics communications},
language = {eng},
pages = {376-381},
publisher = {Elsevier B.V},
title = {Microscopic models of traveling wave equations},
volume = {121-122},
year = {1999},
}

@article{effectOfMicroNoise,
issn = {0022-4715},
journal = {Journal of statistical physics},
pages = {269--282},
volume = {103},
publisher = {Springer Verlag},
number = {1-2},
year = {2001},
title = {Effect of Microscopic Noise on Front Propagation},
copyright = {Distributed under a Creative Commons Attribution 4.0 International License},
language = {eng},
author = {Brunet, Éric and Derrida, Bernard},
}

@misc{me,
      title={{Critical Drift for Brownian Bees and a Reflected Brownian Motion Invariance Principle}}, 
      author={Jacob Mercer},
      year={2024},
      eprint={2412.04527},
      archivePrefix={arXiv},
      primaryClass={math.PR},
      note={\textit{arXiv:2412.04527}},
      howpublished={\url{https://arxiv.org/abs/2412.04527}}, 
}

@article{anulova,
author = {Anulova, S. V.},
address = {Philadelphia},
copyright = {[Copyright] © Society for Industrial and Applied Mathematics},
issn = {0040-585X},
journal = {Theory of probability and its applications},
language = {eng},
number = {2},
pages = {362-366},
publisher = {Society for Industrial and Applied Mathematics},
title = {{On Markov Stopping Times with a Given Distribution for a Wiener Process}},
volume = {25},
year = {1981},
}

@article{cccs06,
author = {Cheng, Lan and Chen, Xinfu and Chadam, John and Saunders, David},
address = {Philadelphia},
copyright = {[Copyright] © 2006 Society for Industrial and Applied Mathematics},
issn = {0036-1410},
journal = {SIAM journal on mathematical analysis},
language = {eng},
number = {3},
pages = {845-873},
publisher = {Society for Industrial and Applied Mathematics},
title = {Analysis of an Inverse First Passage Problem from Risk Management},
volume = {38},
year = {2006},
}

@article{cccs11,
author = {Chen, Xinfu and Cheng, Lan and Chadam, John and Saunders, David},
address = {Hayward},
copyright = {Copyright © 2011 Institute of Mathematical Statistics},
issn = {1050-5164},
journal = {The Annals of applied probability},
language = {eng},
number = {5},
pages = {1663-1693},
publisher = {Institute of Mathematical Statistics},
title = {Existence and uniqueness of solutions to the inverse boundary crossing problem for diffusions},
volume = {21},
year = {2011},
}

@article{ccs22,
author = {Chen, Xinfu and Chadam, John and Saunders, David},
issn = {0036-1410},
journal = {SIAM journal on mathematical analysis},
language = {eng},
number = {4},
pages = {4695-4720},
title = {Higher-Order Regularity of the Free Boundary in the Inverse First-Passage Problem},
volume = {54},
year = {2022},
}

@article{ekstromJanson,
author = {Ekström, Erik and Janson, Svante},
address = {Hayward},
copyright = {Copyright © 2016 Institute of Mathematical Statistics},
issn = {1050-5164},
journal = {The Annals of applied probability},
language = {eng},
number = {5},
pages = {3154-3177},
publisher = {Institute of Mathematical Statistics},
title = {The inverse first-passage problem and optimal stopping},
volume = {26},
year = {2016},
}

@book {kallenberg,
    AUTHOR = {Kallenberg, Olav},
     TITLE = {Foundations of modern probability},
    SERIES = {Probability Theory and Stochastic Modelling},
    VOLUME = {99},
   EDITION = {Third},
 PUBLISHER = {Springer, Cham},
      YEAR = {[2021] \copyright 2021},
     PAGES = {xii+946},
      ISBN = {978-3-030-61871-1; 978-3-030-61870-4},
   MRCLASS = {60-01 (60A10 60G05)},
  MRNUMBER = {4226142},
MRREVIEWER = {Myron\ Hlynka},
       DOI = {10.1007/978-3-030-61871-1},
       URL = {https://doi.org/10.1007/978-3-030-61871-1},
}

\end{document}